%% file: arxiv.tex
\documentclass[10pt]{article} 
 
\usepackage[top=3cm, bottom=3cm, left=4cm, right=4cm]{geometry}
\usepackage{graphicx}

\usepackage{microtype}

\usepackage{hyperref}

\usepackage[parfill]{parskip}

\usepackage{url}

\usepackage{subfigure} 

\usepackage{amsthm}
\usepackage{amssymb}
\usepackage{amsmath}

\usepackage[nameinlink,capitalize]{cleveref}

\usepackage{bibunits}

\usepackage{algorithm}
\usepackage{algpseudocode}

\usepackage{booktabs}

 \usepackage{cite}
 
\usepackage{xcolor}
 \definecolor{dgreen}{rgb}{0.00,0.49,0.00}
\definecolor{dblue}{rgb}{0,0.08,0.75}
\hypersetup{
    colorlinks = true,
    citecolor=dgreen,
    urlcolor=dblue,
    linkcolor=dblue
}

\usepackage[shortlabels]{enumitem}

\usepackage{newtxtext}      

\newcommand{\vol}{\operatorname{vol}}
\newcommand{\range}{\operatorname{range}}
\newcommand{\ran}{\range}
\newcommand{\pdm}{\mathbb{S}_+}

\newcommand{\Cmul}{\mathsf{M}}

\newcommand{\Cdiff}{\mathsf{D}_m}

\input{macro}

\title{Finding Global Minima via Kernel Approximations}

\author{Alessandro Rudi \and Ulysse Marteau-Ferey \and Francis Bach}

\date{
		{\small 
			INRIA - Département d’Informatique de l’École Normale Supérieure, \\
		PSL Research University, 2 rue Simone Iff, 75012, Paris, France \\
		\texttt{\{alessandro.rudi, ulysse.marteau-ferey, francis.bach\}@inria.fr}
	} \\   
	
	$ $\\

	18 December 2020
}

\begin{document}


\maketitle

\begin{abstract}
We consider the global minimization of smooth functions based solely on function evaluations. Algorithms that achieve the optimal number of function evaluations for a given precision level typically rely on explicitly constructing an approximation of the function which is then minimized with algorithms that have exponential running-time complexity. In this paper, we consider an approach that jointly models the function to approximate and finds a global minimum. This is done by using infinite sums of square smooth functions and has strong links with polynomial sum-of-squares hierarchies. Leveraging recent representation properties of reproducing kernel Hilbert spaces, the infinite-dimensional optimization problem can be solved by subsampling in time polynomial in the number of function evaluations, and with theoretical guarantees on the obtained minimum. 
Given $n$ samples, the computational cost is $O(n^{3.5})$ in time, $O(n^2)$ in space, and we achieve a convergence rate to the global optimum that is $O(n^{-m/d + 1/2 + 3/d})$ where $m$ is the degree of differentiability of the function and $d$ the number of dimensions. The rate is nearly optimal in the case of Sobolev functions and more generally makes the proposed method particularly suitable for functions which have a large number of derivatives. Indeed, when $m$ is in the order of $d$, the convergence rate to the global optimum does not suffer from the curse of dimensionality, which affects only the worst case constants (that we track explicitly through the paper).
\end{abstract}

\section{Introduction}

We consider  the general problem of unconstrained optimization. Let $f:\R^d \to \R$ be a possibly non-convex function. Our goal is to solve the following problem
\eqal{\label{eq:prob-true}
\min_{x \in \R^d} f(x).
}
In particular, we will consider the setting where (a) the function is smooth, that is, $f \in C^m(\R^d)$ with $m \in \Np$ ($f$ $m$-times continuously differentiable), and (b) we are able to evaluate it on given points, without the need of computing the gradient. For this class of problems there are known lower-bounds \cite{novak2006deterministic,nesterov2013introductory} that show that it is not possible to achieve a global minimum with error $\eps$ with less than $O(\eps^{-d/m})$ function evaluations. In this paper, we want to achieve this lower bound in terms of function evaluations, while having an optimization algorithm which has a running-time which is polynomial in the underlying dimension and the number of function evaluations.

Several methods are available to solve this class of problems. For example, the function~$f$ can be approximated from its values at $n$ sampled points, and the approximation of the function globally minimized instead of $f$. If the approximation is good enough, then this can be optimal in terms of $n$, but computationally infeasible. Optimal approximations can be obtained by multivariate polynomials~\cite{ivanov1971optimum} or functions in Sobolev spaces~\cite{novak2008tractability}, with potentially adaptive ways of selecting points where the function is evaluated~(see, e.g.,~\cite{osborne2009gaussian} and references therein). Alternatively, when the function is itself a polynomial, algorithms based on the ``sum-of-squares'' paradigm can be used, but their computational complexity grows polynomially on $d^{r/2}$, where $r$ is in the most favorable situations the order of the polynomial, but potentially larger when so-called hierarchies are used~\cite{lasserre2001global,laurent2009sums,lasserre2010moments}.

It turns out that the analysis of lower-bounds on the number of function evaluations shows an intimate link between function interpolation and function minimization, i.e., the lower-bounds of one problem are the same for the other problem. However, existing methods consider a two-step approach where (1) the function is approximated optimally, and (2) the approximation is minimized.
In this paper, we consider a joint approach where approximation and optimization are done \emph{jointly}.

We derive an algorithm that cast the possibly non-convex problem in \cref{eq:prob-true} in terms of a simple convex problem based on a non-parametric representation of  non-negative functions via positive definite operators~\cite{marteau2020non}. As shown below, it can be considered as an infinite-dimensional counter-part to polynomial optimization with sums of squares, with two key differences: (1) the relaxation is always tight for the direct formulation, and (2) the computational cost does not depend on the dimension of the model (here infinite anyway), by using a subsampling algorithm and a  computational trick common in statistics and machine learning.

The resulting algorithm with $n$ sampled points will be able to achieve
an error of $\eps =$ $O ( n^{-m/d + 3/d + 1/2})$ as soon as $m\geq3+d/2$, with $n$ function evaluations to reach the global minimum with precision $\eps$, and a computational complexity of $O(n^{3.5} \log ( 1 / \eps) )$ (with explicit constants). This is still not the optimal complexity in terms of number of function evaluations (which is $\eps = O( n^{-m/d})$), but this is achieved with a polynomial-time algorithm in $n$. This is particularly interesting in the contexts where the function to be optimized is very smooth, i.e., $m \gg d$, possibly $C^\infty$ or a polynomial. For example, if the function is differentiable at least $d+3$ times, even if non-convex, the proposed algorithm finds the global minimum with error $O(n^{-1/2})$ and time $O(n^{3.5} \log n).$
Note that the (typically exponential) dependence on the dimensionality $d$ is only in the constants and tracked explicitly in the rest of the paper. 

Moreover the  algorithm  is based on simple interior-point methods for semidefinite programming, directly implementable and based only on function evaluations and matrix operations. It can thus leverage multiple GPU architectures to reach large values of $n$, which are needed when the dimension grows.

\section{Outline of contributions}\label{sec:sketch}

In this section, we present our framework, our algorithm and summarize the associated guarantees.

Denote by $\zeta \in \R^d$ a global minimizer of $f$ and assume to know a bounded open region $\Omega \subset \R^d$ that contains $\zeta$. We start with a straightforward and classical convex characterization of the problem in \cref{eq:prob-true}, with infinitely many constraints:
\eqal{\label{eq:prob-convex-version}
\max_{c \in \R} ~c \quad \mbox{ such that } \quad  \forall x \in \Omega, \ f(x) \geq c.
}
Note that the solution $c_*$ of the problem above corresponds to $c_* = f(\zeta) = f_\ast$, the global minimum of $f$. The problem above is convex, but typically intractable to solve, due to the dense set of inequalities that $c$ must satisfy. 

To solve \cref{eq:prob-convex-version} our main idea is to represent the dense set of inequalities in terms of a dense set of {\em equalities} and then to approximate them by subsampling.

\paragraph{Tight relaxation.} 
We start by introducing a quadratic form $\scal{\phi(x)}{A \phi(x)}$ with $A$ a self-adjoint positive semidefinite operator from $\hh$ to $\hh$, for a suitable map $\phi: \Omega \to \hh$ and an infinite-dimensional Hilbert space $\hh$, to define the following problem
\eqal{\label{eq:prob-intermediate}
\max_{c \in \R,\  A \in \pdm(\hh)} ~c \quad \mbox{ such that } \quad  \forall x \in \Omega, \ f(x) - c = \scal{\phi(x)}{A \phi(x)},
}
where $\pdm(\hh)$ is the set of bounded self-adjoint positive semi-definite operators on~$\hh$. 

The problem in \cref{eq:prob-intermediate} has a smaller optimized objective function than  the problem in \cref{eq:prob-convex-version} because we constrain $A$ to be positive semi-definite and any feasible point for  \cref{eq:prob-intermediate} is feasible for \cref{eq:prob-convex-version}. In fact, when $f$ is a polynomial and $\phi(x)$ is composed of monomials of degree less than half the degree of $f$ (and thus $\hh$ finite-dimensional), then we recover the classical ``sum-of-squares'' relaxation of polynomial optimization. In that situation, the relaxation is tight only if $f-f_\ast$ is itself a sum-of-squares, which is known to not always be the case. Then, to make the relaxation tight, several hierarchies of polynomial optimization problems have been considered using polynomials of increasing degrees~\cite{lasserre2001global,laurent2009sums,lasserre2010moments}.

In this paper, we consider a well-chosen  infinite-dimensional space $\hh$, and we prove that if~$f$ is  smooth enough (i.e., $m$-times differentiable with $m > 3+d/2$), under mild geometrical assumptions on $f$ then there always exists a map~$\phi$, and a finite rank $A_\ast \in \pdm(\hh)$ for which the problem in \cref{eq:prob-convex-version} and the one above are equivalent, that is, the relaxation is tight. 

Note that, the resulting $\phi$, despite being infinite-dimensional, has an explicit and easy-to-compute ($O(d)$ in memory and time) inner product $k(x,x') = \scal{\phi(x)}{\phi(x')}$ that will be the only quantity required to run the algorithm. We will thus use Hilbert spaces $\hh$ which are reproducing kernel Hilbert spaces~\cite{berlinet2011reproducing}, such as Sobolev spaces~\cite{adams2003sobolev}.

\paragraph{Subsampling.}
We approximate the problem above as follows. Given a finite set $\widehat{X} = \{x_1,\dots,x_n\}$ which is a subset of~$\Omega$, we restrict the inequality in \cref{eq:prob-intermediate} to only $x_1, \dots, x_n$.  

Without further assumptions, subsampling cannot work since, while the function $f$ is assumed smooth enough, the map $x \mapsto  \scal{\phi(x)}{A \phi(x)}$ needs to be regular enough so that satisfying the equality constraint on $\widehat{X}$ leads to a an approximate satisfaction on all of $\Omega$. We thus need to penalize $A$ in some way, and we consider the trace of $A$ and solve the problem
\eqal{\label{eq:prob-intermediate-trace}
\max_{c \in \R,\  A \in \pdm(\hh)} c - \lambda \tr(A) \quad \mbox{ such that } \quad \forall i \in \{1,\dots,n\}, \ f(x_i) - c = \scal{\phi(x_i)}{A \phi(x_i)},
}
for some positive $\lambda$ (with the implicit assumption that we optimize over operators $A$ with finite trace). We show in this paper that solving 
\cref{eq:prob-intermediate-trace} leads to an approximate optimum of the original problem in  \cref{eq:prob-convex-version}, when $n$ is large enough and $\lambda$ small enough. However it is   still formulated in an infinite-dimensional space.

\begin{algorithm}[t]
\caption{Global minimum. Given $f:\R^d\to \R$, $\Omega, n \in \Np, \la > 0, s > d/2$.}\label{alg:glm}
\begin{algorithmic}[1]
\State $\widehat{X} ~\gets~ \{x_1,\dots, x_n\}$\Comment{Sampled i.i.d. uniformly on $\Omega$}
\State $f_j ~\gets~ f(x_j), ~~\forall j \in [n]$ 
\Statex\vspace{-0.2cm}
\Statex \hspace{-0.5cm} Features computation
\State $K_{ij} ~\gets~ k(x_i,x_j) ~~ i,j \in [n]$ \Comment{$k$ Sobolev~kernel of smoothness $s$, \cref{eq:sobolev-kernel}}
\State $R ~\gets~ \texttt{cholesky}(K)$ \Comment{upper triangular Cholesky}
\State $\Phi_j = j$-th column of $R$,  $~~ \forall j \in [n]$
\Statex\vspace{-0.2cm}
\Statex \hspace{-0.5cm} Solution of the approximate problem (use any algorithm in~\cref{sec:algorithm})
\State $\hat{c} ~\gets~ \max_{c \in \R, B \in \pdm(\R^{n})} ~c - \la \tr(B) ~~~ \mbox{ such that }~~~\forall j \in [n],  f_j - c = \Phi_j^\top B \Phi_j$ \label{state:alg-glm-sdp}
\State \textbf{return} $\hat{c}$ 
\end{algorithmic}
\end{algorithm}

\paragraph{Finite-dimensional algorithm.}
We can now leverage the particular choice of penalty by the trace of $A$ and the choice of Hilbert space. Indeed, for reproducing kernel Hilbert spaces, then, following~\cite{marteau2020non}, we only need to solve the problem in the finite-dimensional Hilbert space spanned by $\phi(x_1),\dots, \phi(x_n)$, that is, we only need to look at $A$ of the form $A = \sum_{i,j=1}^n C_{ij} \phi(x_i) \otimes \phi(x_j)$ for some positive semi-definite \emph{matrix} $C \in \R^{n \times n}$. We can then write $\tr(A) = \tr ( C K )$, with $K \in \R^{n \times n}$ the matrix of dot-products with $K_{ij} = \langle \phi(x_i),\phi(x_j)\rangle = k(x_i,x_j)$, and $\langle\phi(x_i), A \phi(x_i)\rangle  = ( K C K )_{ii}$.

Consider the Cholesky decomposition   of $K$ as $K = R^\top R$, with $R \in \R^{n \times n}$ upper-triangular. We can directly solve for $B = R C R^\top$, noting that $KCK = R^\top B R$ and $\tr ( C K ) = \tr (B)$. We can thus use 
a representation in terms of finite-dimensional vectors $\Phi_1,\dots,\Phi_n \in \R^n$ defined as the columns of 
$  R  $.
 We thus study the following problem,
\eqal{\label{eq:prob-relax}
\max_{c \in \R, \ B \in \pdm(\R^n)} ~c - \la \tr(B) \quad \mbox{ such that } \quad \forall i \in \{1,\dots,n\}, \ f(x_i) - c = \Phi_i^\top B \Phi_i.
}
From an algorithmic viewpoint the problem above can be solved  efficiently since this is a semi-definite program.
We show in \cref{sec:algorithm} how we can apply Newton method and classical interior-point algorithms, leading to a computational complexity of $O(n^{3.5} \log(1/\eps))$ in time and $O(n^2)$ in space.

Note that in the context of sum-of-squares polynomials, the relationship with reproducing kernel Hilbert spaces had been explored for approximation purposes after a polynomial optimization algorithm is used~\cite{marx2020semi}. In this paper, we propose to leverage kernel methods \emph{within} the optimization algorithm. 

\paragraph{Why not simply subsampling the inequality?} One straightforward algorithm is to subsample the dense set of \emph{inequalities} in \cref{eq:prob-convex-version}. Doing this will simply lead to outputting $\min_{i \in \{1,\dots,n\} } f(x_i)$. Subsampling the dense set of \emph{equalities} in \cref{eq:prob-intermediate} allows to use smooth interpolation tools. When $\lambda=0$, the optimal value is also $\min_{i \in \{1,\dots,n\} } f(x_i)$ (if the kernel matrix is invertible, see \cref{sec:algorithm}), but for $\lambda > 0$, we can leverage smoothness as shown below.

\paragraph{Theoretical guarantees.}

From a theoretical viewpoint, denoting by $\hat{c}$ the minimizer of \cref{eq:prob-relax}, we provide upper bounds for $|f_* - \hat{c}|$ with explicit constants and that hold under mild geometrical assumptions on $f$. We prove that the bound depends on how the points in $\widehat{X} = \{x_1,\dots,x_n\}$ are chosen. In particular we prove that when they are chosen uniformly at random on $\Omega$, the problem in \cref{eq:prob-relax} achieves the global minimum with error $\eps$ with a precise dependence on $n$.

The results in this paper hold under the following assumptions.

\ba[Geometric properties on $\Omega$ and $f$]\label{asm:geom-f}
The following holds:
\begin{enumasm}
\item\label{asm:geom-f:a} Let $\Omega = \cup_{x \in S} B_r(x)$, where $S$ is a bounded subset of $\R^d$ and $B_r(x)$ is the open ball of radius $r > 0$, centered in $x$. 
\item\label{asm:geom-f:b} The function $f$ is in $C^2(\R^d)$. $\Omega$ contains at least one global minimizer. The minimizers in $\Omega$ are isolated points with strictly positive Hessian and their number is finite. There is no minimizer on the boundary of $\Omega$.


\end{enumasm}
\ea
Note that \cref{asm:geom-f:a} can be easily relaxed to $\Omega$ having locally Lipschitz-continuous boundaries \cite[Section 4.9]{adams2003sobolev}. \cref{asm:geom-f:b} is satisfied if all global minimizers of $f$ are in $\Omega$, and are second-order strict local minimizers.

\bt[Main result, informal]\label{thm:informal}
Let $\Omega \subset \R^d$ be a ball of radius $R > 0$.  Let $s > d/2$ and let $k$ be Sobolev kernel of smoothness $s$ (see \cref{ex:sobolev-kernel}). Let $f \in C^{s+3}(\R^d)$ and that satisfies \cref{asm:geom-f:b}. Let $\hat{c}$ be the result of \cref{alg:glm} executed with $n \in \Np$ points chosen uniformly at random in $\Omega$ and $\la > 0$. Let $\delta > 0$. There exist $n_{s,d,\delta}, C_{s,d} > 0$ such that, when $n > n_{s,d,\delta}$, and
$$\la ~~\geq~~ C_{s,d} ~ n^{- {s}/{d} +{1}/{2}} ~(\log\tfrac{n}{\delta})^{ {s}/{d}-{1}/{2}},$$
then, with probability at least $1-\delta$,
$$|\hat{c} - f_*| ~~\leq~~ 3\,\la\,\left(\tr(A_*) + |f|_{\Omega,\lceil s-{d}/{2} \rceil}\right),$$
where $A_\ast$ is any solution of \cref{eq:prob-intermediate}.
\et
Note that $A_\ast$ exists since $f \in C^{s+3}(\R^d)$ and it satisfies the geometrical mild condition in \cref{asm:geom-f:b} (as we prove in \cref{sec:inf-dim}), and that all constants can be made explicit (see \cref{thm:alg-glm}). From the result above, and with $m = s+3$, for $s > d/2$, we can achieve an error of order $\displaystyle {n^{-s/d+1/2}}$, which translates to   $\eps = O ( n^{-m/d + 3/d + 1/2})$ as soon as $m> d/2 + 3$. The rate for the class of functions $C^m(\Omega)$ is sub-optimal by a factor $1/2 + 3/d$. In the following remark we are going to show that our algorithm achieves nearly-optimal convergence rates when the function to optimize is in a Sobolev space. Denote by $W^s_2(\Omega)$ the Sobolev space of squared-integrable functions of smoothness $s > 0$, i.e., the space of functions whose weak derivatives up to order $s$ are square-integrable on $\Omega$, (see \cite{adams2003sobolev}). 

\br[\bf Nearly optimal rates for Sobolev spaces.]\label{rm:sobolev-nearly-optimal} If $\Omega$ satisfies \cref{asm:geom-f:a}, $f$ satisfies \cref{asm:geom-f:b} and $f \in W^{s}_2(\Omega)$, with $s > d/2 + 3$, then \cref{alg:glm} with Sobolev kernel of smoothnes $s-3$ achieves the convergence rate 
$$
O\big(n^{-{s}/{d} + {1}/{2} + {3}/{d}}\big),
$$
modulo logarithmic factors, as proven in \cref{thm:alg-glm}. When $d$ is large, then the error exponent is asymptotically optimal, since the term $3/d$ becomes negligible, leading to the optimal exponent $-{s}/{d} + {1}/{2}$ (see, e.g., \cite[Prop.~1.3.11]{novak2008tractability}).
\er

\paragraph{Finding the global minimizer.}
In \cref{sec:minimizer} we derive an extension of the problem in  \cref{eq:prob-relax}, with the goal of finding the global minimizer. Under the additional assumption that the minimizer is unique we obtain the similar rate as \cref{thm:bound-prob-sampled} for the localization of the global minimizer.

\paragraph{Warm restart scheme for linear rates.} Applying a simple warm restart scheme, we prove, in \cref{sec:warm-start}, that when $f$ has a unique global minimum, then it is possible to achieve it with error~$\eps$, with a number of observations that is only logarithmic in $\eps$
$$n = O(C_{d,m} \log(1/\eps)),$$
for some constant $C_{d,m}$ that can be exponential in $d$ (note that the added assumption of unique minimizer makes this result not contradict the lower bound in $\eps^{-d/m}$).

\paragraph{Rates for functions with low smoothness $m \leq d/2$ or functions that are not in $\hh$.} In \cref{sec:low-smoothness} we study a variation of the problem in  \cref{eq:prob-relax} that allows to have some error $\tau > 0$ on the constraints. When $f \in C^{m+2}(\Omega)$, by tuning $\tau$ appropriately with respect to $\la$, we show that \cref{alg:glm} applied on this different formulation achieves an error in the order 
$$O\big(\,n^{-\frac{m}{2d}(1-(d - m)/(2r - m))}\,\big),$$
where $r$  is now the index of the Sobolev kernel and can be chosen arbitrarily large. The exponent of the rate above matches the optimal one for $C^{m+2}$ functions (that is $-(m+2)/d$) up to a multiplicative factor of $\frac{1}{2} \frac{1}{1+2/m}$.
 
\paragraph{Relationship to polynomial optimization.}
When $f$ is a polynomial of degree $2r$, then it is natural to consider $\phi(x)$ composed of all monomials of degree less than $r$, leading to a space~$\hh$ of dimension ${d+r \choose r}$. All polynomials can be represented as $f(x) = c + \phi(x)^\top A \phi(x)$ for some symmetric matrix $A$. When $A \succcurlyeq 0$, by using its eigendecomposition, we can see that the polynomial $x \mapsto  \phi(x)^\top A \phi(x)$ is a sum-of-squares polynomial.

However, in general $A$ may not be positive semi-definite, as non-negative polynomials are not all sum-of-squares. Moreover, even when there exists a matrix $A \succcurlyeq 0$, the corresponding~$c$ may not be the minimum of $f$ (it only needs to be a lower bound)---see, e.g.,~\cite{lasserre2001global} and references therein.

If $f(x) - f_\ast$ is a sum of squares, then, with $\lambda = 0 $ and $n = {d+2r \choose 2r}$ points (to ensure that subsampling is exact), we exactly get the minimum of $f$, as we are solving \emph{exactly} the usual optimization problem.

When $f(x)-f_\ast$ is not a sum of squares, then a variety of hierarchies have been designed, that augment the problem dimensionality to reach global convergence\cite{lasserre2001global,laurent2009sums,lasserre2010moments}. In \cref{sec:polynomial}, we show how our framework fits with one these hierarchies, and also can provide computational gains.

Note that our framework, by looking directly at an infinite-dimensional space circumvents the need for hierarchies, and solves a single optimization problem. The difficulty is that it requires sampling. Moreover by using only kernel evaluations, we circumvent the explicit construction of a basis for $\hh$, which is computationally cumbersome when $d$ grows.

\paragraph{Organization of the paper.}
The paper is organized as follows: in \cref{sec:setting}, we present the kernel setting our paper relies on; then, in \cref{sec:inf-dim}, we analyze the infinite-dimensioal problem and show its equivalence with global minimization. Then, in \cref{sec:finite-dim}, we present
our theoretical guarantee for the finite-dimensional algorithm, as summarized in \cref{thm:informal}. In \cref{sec:algorithm} we present the dual algorithm based on self-concordant barriers and the damped Newton algorithm. In \cref{sec:minimizer}, we present our extension to find the global minimizer, while in \cref{sec:certificate}, we provide certificates of optimality for potentially inexactly solved problems. 
In \cref{sec:polynomial}, we discuss further relationships with polynomial hierarchies, and provide illustrative experiments in \cref{sec:experiments}. We conclude in \cref{sec:discussion} with a discussion opening up to many future problems.

\section{Setting}
\label{sec:setting}

In this section, we first introduce some definitions and notation about {\em reproducing Kernel Hilbert spaces} in \cref{sec:intro-rkhs} (for more details, see \cite{aronszajn1950theory,paulsen2016introduction}), and present our detailed assumptions in \cref{sec:rkhs-assumptions}. In \cref{sec:inf-dim} we show how our infinite-dimensional  sum-of-squares representation can be built, and in  \cref{sec:finite-dim} we provide guarantees on subsampling.

\subsection{Definitions and notation} \label{sec:notation}
\label{sec:intro-rkhs}
In this section we denote by $u \cdot v$, $a \circ v$ respectively the pointwise multiplication between the functions $u$ and $v$, and the composition between the functions $a$ and $v$. We denote by $\Nz$ the set of natural numbers including $0$, by $\Np$ the set $\Np = \Nz \setminus \{0\}$ and $[n]$ the set $\{1,\dots,n\}$ for $n \in \Np$. We will always consider $\R^d$ endowed with the Euclidean norm $\|\cdot\|$ if not specified otherwise. Moreover we denote by $B_r(z)$ the open ball $B_r(z) = \{x \in \R^d~|~ \|x-z\| < r\}$.
Let $\Omega \subseteq \R^d$ be an open set. Let $\alpha \in \Nz^d$. We introduce the following {\em multi-index notation} $|\alpha| = \alpha_1 + \dots + \alpha_d$ and $\partial_x^\alpha = \frac{\partial^{|\alpha|} }{\partial x^{\alpha_1}_1 \dots \partial x^{\alpha_d}_d}$ \cite{adams2003sobolev}. For $m \in \Nz$, and $\Omega$ an open set of $\R^d$, denote by $C^m(\Omega)$ the set of $m$-times differentiable functions on $\Omega$ with continuous $m$-th derivatives. For any function $u$ defined on a superset of $\Omega$ and $m$ times differentiable on $\Omega$, define the following semi norm.
\begin{equation}\label{eq:df_seminorm}
   |u|_{\Omega,m} ~~= ~~ \max_{|\alpha| = m} ~ \sup_{x \in \Omega} \big|\partial^\alpha u(x)\big|.
\end{equation}

\paragraph{Positive definite matrices and operators.} Let $\hh$ be a Hilbert space, endowed with the inner product $\scal{\cdot}{\cdot}$. Let $A:\hh \to \hh$ be a linear operator and denote by $A^*$ the adjoint operator, by $\tr(A)$ the trace of $A$ and by $\|\cdot\|_F$ the Hilbert-Schmidt norm $\|A\|^2_F = \tr(A^*A)$. We always endow $\R^p$ with the standard inner product $x^\top y = \sum_{i=1}^p x_i y_i$ for any $x,y \in \R^p$.
In the case $\hh = \R^p$, with the standard inner product, then $A \in \R^{p \times p}$ is a matrix and the Hilbert-Schmidt norm corresponds to the Frobenius norm. We say that $A \succeq 0$ or $A$ is a {\em positive operator} (positive matrix if $\hh$ is finite dimensional), when $A$ is bounded, self-adjoint, and
$ \scal{u}{A u} \geq 0, \  \forall u \in \hh.$
We denote by $\pdm(\hh)$ the space of positive operators on $\hh$. Moreover, we denote by $A \succ 0$, or $A$ strictly positive operator, the case $\scal{u}{A u} > 0$ for all  $u \in \hh$ such that $u \neq 0$.

\paragraph{Kernels and reproducing kernel Hilbert spaces.} For this section we refer to \cite{aronszajn1950theory,steinwart2008support,paulsen2016introduction}, for more details (see also \cref{sec:rkhs-appedix}, \cpageref{sec:rkhs-appedix}).
Let $\Omega$ be a set. A function $k:\Omega \times \Omega \to \R$ is called a {\em positive definite kernel} if all matrices of pairwise evaluations are positive semi-definite, that is, if it satisfies the following equation
$$ \sum_{i,j=1}^n \alpha_i \alpha_j k(x_i,x_j) \geq 0, \qquad \forall n \in \Nz, \alpha_1,\dots,\alpha_n \in \R, x_1,\dots, x_n \in \Omega.$$
Given a kernel $k$, the {\em reproducing kernel Hilbert space} (RKHS) $\hh$,  with the associated inner product $\scal{\cdot}{\cdot} $, is a space of real functions with domain $\Omega$, with the following properties.
\begin{enumerate}[(a)]
    \item The function $k_x = k(x,\cdot)$ satisfies $k_x \in \hh $ for any $x \in \Omega$.
    \item The inner product satisfies $\scal{f}{k_{x}} = f(x)$ for all $f \in \hh $, $x \in \Omega$. In particular $\scal{k_{x'}}{k_{x}} = k(x',x)$ for all $x,x' \in \Omega$.
\end{enumerate}
In other words, function evaluations are uniformly bounded and continuous linear forms and the $k_x$ are the evaluation functionals.
The norm associated to $\hh$ is the one induced by the inner product, i.e., $\|f\|^2  = \scal{f}{f}$. We remark that given a kernel on $\Omega$ there exists a unique associated RKHS on $\Omega$ \cite{berlinet2011reproducing}. 
Moreover, the kernel admits a characterization in terms of a {\em feature map} $\phi$,
$$\phi:\Omega \to \hh, \quad \textrm{defined as} \quad \phi(x) = k(x, \cdot) = k_x, ~~ \forall x \in \Omega.$$
Indeed according to the point (b) above, we have $k(x,x') = \scal{\phi(x)}{\phi(x')}$ for all $x,x' \in \Omega$. We will conclude the section with an example of RKHS that will be useful in the rest of the paper.

\bex[Sobolev kernel \cite{wendland2004scattered}]\label{ex:sobolev-kernel} 
Let $s > d/2$, with $d \in \Np$, and $\Omega$ be a bounded open set. Let 
\eqal{\label{eq:sobolev-kernel}
k_s(x,x') = c_s \|x-x'\|^{s-d/2} {\cal K}_{s-d/2}(\|x-x'\|), \quad \forall x,x' \in \Omega,
}
where ${\cal K}:\R_+ \to \R$ the Bessel function of the second kind (see, e.g., 5.10  in \cite{wendland2004scattered}) and $c_s = \frac{2^{1 + d/2 -s}}{\Gamma(s-d/2)}$. The constant $c_s$ is chosen such that $k_s(x,x) = 1$ for any $x \in \Omega$. In the particular case of $s = d/2 + 1/2$, we have $k(x,x') = \exp( - \| x - x'\|)$. Note that a scale factor is often added as $k(x,x') = \exp( - \| x - x'\|/ \sigma)$ in this last example. In such case, all bounds that we derive in this paper would then have extra factors proportional to powers of $\sigma$.
To conclude, when $\Omega$ has locally Lipschitz boundary (a sufficient condition is \cref{asm:geom-f:a}) then $\hh = W^s_2(\Omega)$,
where $W^s_2(\Omega)$ is the Sobolev space of functions whose weak-derivatives up to order $s$ are square-integrable \cite{adams2003sobolev}. Moreover, in this case $\|\cdot\|_{\hh}$ is equivalent to $\|\cdot\|_{W^s_2(\Omega)}$. 
\eex

Reproducing kernel Hilbert spaces are classically used in fitting problems, such as appearing  in statistics and machine learning, because of function evaluations $f \mapsto f(x)$ are bounded operators for any $x$, and optimization problems involving $f$ only through function evaluations at a finite number of points $x_1,\dots,x_n$, and penalized with the norm $\| f\|$, can be solved by looking only a $f$ of the form $f(x) = \sum_{i=1}^n \alpha_i k(x,x_i)$~\cite{aronszajn1950theory,paulsen2016introduction}. We will use an extension of this classical ``representer theorem'' to operators and spectral norms in \cref{sec:finite-dim}.

\subsection{Precise assumptions on reproducing kernel Hilbert space}
\label{sec:rkhs-assumptions}

On top of \cref{asm:geom-f} (made on the function $f$ and the set $\Omega$), we make the following assumptions on the space $\hh$ and the associated kernel $k$.

\ba[Properties of the space $\hh$]\label{asm:H-rich}
Given a bounded open set $\Omega \subset \R^d$, let $\hh$ be a space of functions on $\Omega$ with norm $\|\cdot\|_\hh$, satisfying the following conditions
\begin{enumasm}
\item\label{asm:H-rich:a}  $w|_\Omega \in \hh , ~\forall w \in C^\infty(\R^d)$. Moreover there exists $\Cmul \geq 1$ such that
$$\|u \cdot v\|_{\hh } \leq \Cmul \|u\|_{\hh } \|v\|_{\hh }, \quad \forall u,v \in {\hh }.$$ 
\item\label{asm:H-rich:b} $a \circ v \in \hh $, for any $a \in C^\infty(\R^p)$, $v = (v_1, \dots, v_p)$, $v_j \in \hh $, $j \in [p]$.
\item\label{asm:H-rich:c} Let $z \in \R^d, r > 0$ s.t. the ball $B_r(z)$ is in $\Omega$. For any $u \in \hh$, there exists $g_{r,z} \in \hh$ s.t.
$$g_{r,z}(x) = \int_0^1(1-t) u(z + t(x-z)) dt, \quad  \forall x \in B_r(z).$$
\item\label{asm:H-rich:d} $\hh$ is a RKHS with associated kernel $k$. For some $m \in \Np$ and some $\Cdiff \geq 1$, the kernel $k$ satisfies
$$\max_{|\alpha| = m} \sup_{x,y \in \Omega} |\partial^\alpha_x \partial^\alpha_y k(x,y)| \leq  \Cdiff^2 < \infty.$$
\end{enumasm}
\ea
\cref{asm:H-rich:a,asm:H-rich:b,asm:H-rich:c} above require essentially that $\hh$ contains functions in $\hh$ can be multiplied by other functions in $\hh$, or by infinitely smooth functions, that can be composed with infinitely smooth functions, or integrated, and still be in $\hh$. Moreover \cref{asm:H-rich:d} requires that $\hh$ is a RKHS with a kernel that is $m$-times differentiable. An interesting consequence of \cref{asm:H-rich:d} is the following remark (for more details, see for example \cite[Corollary 4.36]{steinwart2008support}).
\br\label{rm:k-differentiable}
 \cref{asm:H-rich:d} guarantees that $\hh  \subseteq C^m(\Omega)$ and 
$|u|_{\Omega,m} \leq \Cdiff \|u\|_{\hh }$. 
\er
Note that \cref{asm:H-rich:a,asm:H-rich:b,asm:H-rich:c} are the only required in \cref{sec:inf-dim} to prove the crucial decomposition in \cref{thm:g-good-has-A} and are satisfied by notable spaces (that are not necessarily RKHS) like $C^s(\Omega)$ or Sobolev spaces $W^s_p(\Omega)$ with $s > d/p$ and $p \in [1,\infty]$. Instead, \cref{asm:H-rich:d} is required for the analysis of the finite-dimensional problem and in particular \cref{thm:inequality-scattered-data,thm:bound-prob-sampled}. In the following proposition we show that $W^s_2(\Omega)$ with $s > d/2$ and $\Omega$ satisfying \cref{asm:geom-f:a} satisfy the whole of \cref{asm:H-rich}.

\bp[Sobolev kernels satisfy \cref{asm:H-rich}]\label{prop:sobolev-kernel} 
Let $\Omega$ be a bounded open set of $\R^d$. The Sobolev kernel with $s > d/2$ recalled in \cref{ex:sobolev-kernel} satisfies \cref{asm:H-rich} for any $m \in \Np, m < s - \tfrac{d}{2}$ and
$$\Cmul  =(2\pi)^{d/2}  2^{s+1/2} ,  \qquad \Cdiff = (2\pi)^{d/4}\sqrt{\frac{ \Gamma(m+d/2)\Gamma(s-d/2-m)}{\Gamma(s-d/2)\Gamma(d/2)}}.$$
\ep

The proof of proposition above is in \cref{proof:prop:sobolev-kernel}, \cpageref{proof:prop:sobolev-kernel}. We make a last assumption regarding the differentiability of $f$, namely that $f$ and its second-derivatives are in $\hh$.

\ba[Analytic properties of $f$]\label{asm:analytic-f}
The function $f$ satisfies $f|_\Omega \in C^2(\Omega) \cap \hh$ and $\frac{\partial^2 f}{\partial x_i \partial x_j}|_\Omega \in \hh $ for all $i,j \in [d]$.
\ea

\section{Equivalence of the infinite-dimensional problem}\label{sec:inf-dim}

In \cref{thm:g-good-has-A} and \cref{cor:g-good-has-A}, we provide  a representation of $f-f_\ast$ in terms of an infinite-dimensional, {\em but finite-rank}, positive operator, under basic geometric conditions on $f$ and algebraic properties of~$\hh$. In \cref{thm:Hilbert-isgood} we use this operator to prove that \cref{eq:prob-intermediate} achieves the global minimum of $f$. In this section we analyze the conditions under which the problem in \eqref{eq:prob-intermediate} has the same solution as the one in \cref{eq:prob-convex-version}.

The proof follows by explicitly constructing a bounded positive operator $A_*$ (which will have finite trace) that satisfy $f(x) - f_\ast  = \scal{\phi(x)}{A_\ast \phi(x)}$ for all $x \in \Omega$. Note that, by construction $f - f_\ast$ is a non-negative function. If $w := \sqrt{f - f_\ast} \in \hh$ then $A_* = w \otimes w$ would suffice. However, denoting by $\zeta \in \Omega$ a global minimizer, note that 
$f(\zeta) - f_\ast = 0$ and the smoothness of $\sqrt{f - f_\ast}$ may degrade around $\zeta$, making $\sqrt{f - f_\ast} \notin \hh$ even if $f - f_\ast \in \hh$.

 Here we follow a different approach. In \cref{lm:repr-conv-g} we provide a decomposition that represents the function $f - f_\ast$ locally around each global optimum using the fact that it is locally strongly convex around the minimizers. In the proof of \cref{thm:g-good-has-A} we provide a decomposition of the function far from the optimal points; we then glue these different decompositions via bump functions.

\blm\label{lm:repr-conv-g}
Let $\hh$ be a space of functions on $\Omega$ that satisfy \cref{asm:H-rich:a,asm:H-rich:b,asm:H-rich:c}. Let $\zeta \in \Omega$ and $r,\gamma > 0$. Let $B_r(\zeta) \subset \Omega$ be a ball centered in $\zeta$ of radius $r$ and $g \in C^2(\Omega)$ satisfy $g(\zeta) = 0$,  $\nabla^2 g (x) \succcurlyeq \gamma I$ for $x \in B_r(\zeta)$ and $\tfrac{\partial^2}{\partial x_i \partial x_j}g \in \hh $ for $i,j \in [d]$. Then, there exists $w_j \in \hh , j \in [d]$ such that
\eqal{\label{eq:repr-conv-g}
g(x) = \sum_{j=1}^d w_j(x)^2, \quad \forall x \in B_r(\zeta).
}
\elm
\bpr
Let $x \in B_r(\zeta)$ and consider the function $h(t) = g(\zeta + t(x-\zeta))$ on $[0,1]$. Note that $h(0) = g(\zeta)$ and $h(1) = g(x)$. Taking the Taylor expansion of $h$ of order 1, we have $h(1) = h(0) + h'(0) + \int_{0}^1 (1-t) h''(t)dt$, with  $h(0) = g(\zeta)$, $h'(0) = (x-\zeta)^\top \nabla g(\zeta)$ and $h''(t) = (x-\zeta)^\top \nabla^2 g(\zeta + t(x-\zeta))(x-\zeta)$. Since $g(\zeta) = 0$ by construction and $\nabla g(\zeta) = 0$ since $\zeta$ is a local minimizer of $g$, we have $h(0) = h'(0) = 0$ leading to
\eqal{\label{eq:repr-g}
g(x) = (x-\zeta)^\top R(x) (x-\zeta), \quad R(x) = \int_0^1 (1-t) \nabla^2 g(\zeta + t(x-\zeta)) dt.
}
Note that for $x \in B_{r}(\zeta)$ we have $\nabla^2 g(x) \succcurlyeq \gamma I$ and so $R(x) \succcurlyeq\gamma I$. In particular, this implies that for any $x \in B_r(\zeta)$, $S(x) = \sqrt{R(x)}$ is well defined ($\sqrt{\cdot} : \pdm(\R^d) \rightarrow \pdm(\R^d)$ is the spectral square root, where  for any $M \in \pdm(\R^d)$ and any eigen-decomposition $M = \sum_{j=1}^d \la_j u_j u_j^\top$, $\sqrt{M} =  \sum_{j=1}^d \sqrt{\la_j} u_j u_j^\top$). Thus,
\[\forall x \in B_r(\zeta),~ g(x) = (x-\zeta)^\top S(x)S(x)(x-\zeta) = \sum_{i=1}^d{\left(e_i^\top S(x)(x-\zeta)\right)^2}.\]
The following steps prove the existence of $w_i \in \hh$ such that $w_i|_{B_r(\zeta)} =e_i^\top S(\cdot)(\cdot-\zeta)$. Let $(e_1,...,e_d)$ be the canonical basis of $\R^d$ and $\mathbb{S}(\R^d)$ be the set of symmetric matrices on $\R^d$ endowed with Frobenius norm, in the rest of the proof we identify it with the isometric space $\R^{d(d+1)/2}$ (corresponding of taking the upper triangular part of the matrix and reshaping it in form of a vector).

\noindent{\bf Step 1. } \textit{There exists a function $\overline{R} : \Omega \rightarrow \mathbb{S}(\R^d)$, such that 
$$ \forall i,j \in [d],~ e_i^\top \overline{R}e_j \in \hh \text { and }  \overline{R}|_{B_r(\zeta)} = R.$$}
This is a direct consequence of the fact that $\frac{\partial^2}{\partial x_i \partial x_j}g \in \hh $ for all $i \leq j \in [d]$,  of \cref{asm:H-rich:c} and of the definition of $R$ in \cref{eq:repr-g}.

\noindent{\bf Step 2. } \textit{There exists a function $\overline{S} : \Omega \rightarrow \mathbb{S}(\R^d)$ such that 
$$ \forall i,j \in [d],~  e_i^\top \overline{S}e_j \in \hh \text { and } \forall x \in B_r(\zeta),~ \overline{S} (x)=  \sqrt{R}(x).$$}
Let $\tau := \sup_{x \in B_{r}(\zeta)} \|R(x)\|_{\rm op}= \|\overline{R}(x)\|_{\rm op}$, which is well defined because $R$ is continuous since $g \in C^2(\Omega)$.
Define the compact set $K = \{T \in \mathbb{S}(\R^d)~|~\gamma I \preceq T \preceq \tau I\}$ and the open set $U = \{ T \in \mathbb{S}(\R^d)~|~ \tfrac{\gamma}{2} I \prec T \prec 2\tau I\}$. Note that $ K \subset U \subset \mathbb{S}(\R^d)$.

Fix $i, j \in [d]$ and consider the function $\theta_{i,j} : U \rightarrow \R$ defined by $\theta_{i,j}(M) = e_i^\top \sqrt{M} e_j$. Since the square root $\sqrt{\cdot} : \pdm(\R^d) \rightarrow \pdm(\R^d)$ is infinitely differentiable (see e.g. the explicit construction in \cite{DELMORAL2018259} Thm.~1.1) and $U \subset \pdm(\R^d)$ then $\theta_{i,j}$ is infinitely differentiable on $U$, i.e.,  $\theta_{i,j} \in \C^{\infty}(U)$. By \cref{prop:Cinfty-extension}, since $K$ is a compact set in $U$, there exists $\overline{\theta}_{i,j} \in C_0^{\infty}(\mathbb{S}(\R^d))$ such that $\forall T \in K,~ \overline{\theta}_{i,j}(T) = \theta_{i,j}(T)$.

Define $\overline{S}(x)  = \sum_{i,j \in [d]} {(\overline{\theta}_{i,j}\circ \overline{R})(x) e_ie_j^\top}$ for any $x \in \Omega$. Applying \cref{asm:H-rich:b}, $e_i^\top\overline{S}e_j = \overline{\theta}_{i,j}\circ \overline{R} \in \hh$ since the $\overline{R}_{k,l} \in \hh,~k,l \in [d]$ and $\overline{\theta}_{i,j}$ is in $C_0^{\infty}(\mathbb{S}(\R^d))$. Moreover, by construction, for any $x \in B_r(\zeta)$, we have $\overline{R}(x) = R(x) \in K$ and so
$$\overline{S}_{i,j}(x) = \overline{\theta}_{i,j}(\overline{R}(x)) = \theta_{i,j}(R(x)) = e_i^\top \sqrt{R(x)}e_j. $$
Note that here, we have applied \cref{prop:Cinfty-extension} and \cref{asm:H-rich:b} to $\mathbb{S}(\R^d)$ and not to $\R^{d(d+1)/2}$; this can be made formal by using the linear isomorphism between $\mathbb{S}(\R^d)$ endowed with the Frobenius norm and $\R^{d(d+1)/2}$ endowed with the Euclidean norm.

\noindent{\bf Step 3. } \textit{There exists a function $\overline{h} = (\overline{h}_j)_{j \in [d]}: \Omega \rightarrow \R^d$ such that 
$$\forall j \in [d],\overline{h}_j \in \hh \text { and } \forall x \in B_r(\zeta),~ \overline{h}(x) = x-\zeta. $$}
Fix $j \in [n]$. Define $\overline{B}_r(\zeta) = K \subset U = B_{2r}(\zeta)$ and apply proposition \cref{prop:Cinfty-extension} to $x \in U \mapsto e_j^\top(x-\zeta)$ to get $h_j \in C_0^\infty(\R^d)$ which coincides with $e_j^\top(\cdot-\zeta)$ on $K$ hence on $B_r(\zeta)$. Applying \cref{asm:H-rich:a}, the restriction $\overline{h}_j = h_j |_{\Omega}$ is in $\hh$, and hence $\overline{h} = \sum_{j \in [d]}\overline{h}_j e_j$ satisfies the desired property. 

\noindent{\bf Step 4. } \textit{The $w_i = e_i^\top \overline{S}~ \overline{h},~ i \in [d]$ have the desired property. }

It is clear that the $w_i$ are in $\hh$ as a linear combination of products of functions in $\hh$ (see \cref{asm:H-rich:a}), since $w_i = \sum_{j \in [d]} \overline{S}_{ij}(x) \overline{h}_j(x)$ for any $x \in \Omega$. Moreover, 
\[\sum_{i \in [d]}{w_i^2}  = \overline{h}^\top \overline{S}^\top ~\left(\sum_{i =1}^d{e_i e_i^\top}\right) \overline{S}~ \overline{h}  =\overline{h}^\top \overline{S}^2\overline{h}.\]
Using the previous points, 
\[\forall x \in B_r(\zeta),~ \sum_{i \in [d]}{w_i^2(x)} =\overline{h}^\top(x) \overline{S}^2(x)\overline{h}(x)  = (x-\zeta)^\top R(x)(x-\zeta) = g(x).\]
\epr

Now we are going to use the local representations provided by the lemma above to build a global representation in terms of a finite-rank positive operator. Indeed far from the global optima the function $f - f_\ast$ is strictly positive and so we can take a smooth extension of the square root to represent it and glue it with the local representations around the global optima via bump functions as follows. 

\bt\label{thm:g-good-has-A} 
Let $\Omega$ be a bounded open set and let $\hh$ be a space of functions on $\Omega$ that satisfy \cref{asm:H-rich:a,asm:H-rich:b,asm:H-rich:c}.  Let $f$ satisfy \cref{asm:geom-f:b,asm:analytic-f}. Then there exist $w_1,\dots,w_q \in \hh$ with $q \leq d p + 1$ and $p \in \Np$ the number of minimizers in $\Omega$, such that
\eqal{\label{eq:g-good-A-statement}
f(x) - f_* = \sum_{j \in [q]} w_j(x)^2, \qquad \forall~ x \in \Omega.
}
\et
\bpr
Let $Z = \{\zeta_1,\dots,\zeta_p\}$, $p \in \Np$ be the \textit{non-empty} set of global minima of $f$, according to \cref{asm:geom-f:b}. Denote by $f_* = \min_{x \in \Omega} f(x)$ the global minimum of $f$, and by $g : \Omega \rightarrow \R$ the function $g = f|_\Omega - f_* \mathbf{1}|_\Omega$ where $\mathbf{1}$ is the function $\mathbf{1}(x) = 1$ for any $x \in \R^d$.
\cref{asm:analytic-f} implies that $\nabla^2g = \nabla^2f|_\Omega$ is continuous, an that $\frac{\partial^2 g}{\partial x_i \partial x_j} \in \hh$ for any $i,j \in [d]$.
Moreover, $g \in \hh$. Indeed, by construction $f_*\mathbf{1}$ is in $C^\infty(\R^d)$, and since $\hh$ satisfies \cref{asm:H-rich:a}, $f_*\mathbf{1}
|_\Omega \in \hh$. Since $f|_\Omega \in \hh$ by \cref{asm:analytic-f}, then $g \in \hh$. 

\noindent{\bf Step 1.} \textit{There exists $r>0$ and $\alpha > 0$ such that (i) the $B_r(\zeta_l) ,~ l \in [p]$ are included in $\Omega$ and (ii) for any $x \in \bigcup_{l \in [p]} B_r(\zeta_l)$, it holds $\nabla^2 g(x) \succeq \alpha I$.}

By \cref{asm:geom-f:b}, for all $\zeta \in Z$, $\nabla^2 g(\zeta) \succ 0$. Since $\nabla^2 g$ is continuous, $Z$ is a finite set, and $\Omega$ is an open set, there exists a radius $r > 0$ and $\alpha >0$ such that for all $l \in [p]$, $B_{r}(\zeta_l) \subset \Omega$ and $\nabla^2 g|_{B_r(\zeta_l)} \succeq \alpha I$. 
For the rest of the proof, fix $r,\alpha$ satisfying this property. For any $X \subset \Omega$ denote with $\boldsymbol{1}_X$ the indicator function of a $X$ in $\Omega$. We define $\chi_0 = \boldsymbol{1}_{\Omega \setminus{\bigcup_{l \in [p]}B_{r/2}(\zeta_l)}}$, and $\chi_l = \boldsymbol{1}_{B_r(\zeta_l)},~ l \in [p]$.

\noindent{\bf Step 2.} \textit{ There exists $w_0 \in \hh$ s.t. $w_0^2 \chi_0 =g \chi_0$.} 

 $\Omega$ is bounded and by \cref{asm:geom-f:b}, the set of global minimizers of $f$ included in $\Omega$ is finite and there is no mimimizer of $f$ on the boundary, i.e., there exists $m_1 > 0$ and a compact $K \subset \Omega$ such that $\forall x \in \Omega \setminus{K},~ g(x) \geq m_1$.

Moreover, $f$ has no global optima on the compact $K \setminus \bigcup_{\zeta \in Z}B_{r/2}(\zeta)$ since the set of global optima is $Z$, hence the existence of $m_2 > 0$ such that $\forall x \in K \setminus{\bigcup_{l \in [p] }B_{r/2}(\zeta_l)},~ g(x) \geq m_2$. Taking $m = \min(m_1,m_2)$, it holds $\forall x \in \Omega \setminus{\bigcup_{l \in [p] }B_{r/2}(\zeta_l)},~g(x) \geq m > 0$. Since $f \in C^2(\Omega)$, $f$ is also bounded above on $\Omega$ hence  the existence of $M >0$ such that $g \leq M$. Thus
\[\forall x \in \Omega \setminus{\bigcup_{l \in [p] }B_{r/2}(\zeta_l)},~ g(x) \in I \subset (m/2,2M),\qquad I = [m,M].\]
Since $\sqrt{\cdot} \in C^\infty((m/2,2M))$, $(m/2,2M)$ is an open subset of $\R$ and $I$ is compact, applying \cref{prop:Cinfty-extension}, there exists a smooth extension $s_I \in C^\infty_0(\R)$ such that $ s_I(t) = \sqrt{t}$ for any $t \in I$. Now since $g \in \hh$ and $s_I \in C^\infty_0(\R)$, by  \cref{asm:H-rich:b}, $w_0 := s_I \circ g \in \hh$. Since $\forall x \in \Omega \setminus{\bigcup_{l \in [p] }B_{r/2}(\zeta_l)},~g \in I$, this shows $g \chi_0 = w_0^2\chi_0$.

\noindent{\bf Step 3. }\textit{For all $l \in [p]$, there exists $(w_{l,j})_{j \in [d]} \in \hh^d$ s.t.   $g(x)\chi_l = \sum_{j=1}^d{w_{l,j}^2~\chi_l}$.}

This is an immediate consequence of \cref{lm:repr-conv-g} noting that  $\nabla g(x) \geq \alpha I$ on $B_r(\zeta_l)$.

\noindent{\bf Step 4.}\textit{ There exists $b_l \in C^\infty(\R^d)$ s.t.  $b_l  =  b_l ~\chi_l$ for all $l \in \{0,1,\dots,p\}$ and $\sum_{l=0}^p{b_l^2} = 1$.}

This corresponds to \cref{lm:exists-b2}, \cref{sec:c-infty-results}, \cpageref{sec:c-infty-results} applied to the balls $B_r(\zeta_l),~l \in [p]$.

\noindent{\bf Step 5.} Using all the previous steps
\begin{align*} 
g = \sum_{l = 0}^{p} g~ b_l^2 &= \sum_{l = 0}^{p} g  (\chi_l ~ b_l)^2 =\sum_{l = 0}^{p} { (\chi_l~g ) ~(\chi_l  b_l^2 )}\\
&= (\chi_0 w_0^2)~( \chi_0 ~b_0^2 ) + \sum_{l=1}^p{\big( \chi_l \sum_{j=1}^d{w_{l,j}^2}\big)~\chi_l~b_l^2} \\
& = ([b_0 ~\chi_0]~ w_0)^2 + \sum_{l=1}^p{ \sum_{j=1}^d{([b_l ~\chi_l] ~w_{l,j})^2}} = (b_0 ~ w_0)^2 + \sum_{l=1}^p{ \sum_{j=1}^d{(b_l  ~w_{l,j})^2}}.
\end{align*}

Applying \cref{asm:H-rich:a} to each function inside the squares in the previous expressions yields the result.

\epr

A direct corollary of the theorem above is the existence of $A_* \in \pdm(\hh)$ when $\hh$ is a reproducing kernel Hilbert space satisfying the assumptions of \cref{thm:g-good-has-A}.
\bcor\label{cor:g-good-has-A}
Let $k$ be a kernel whose associated RKHS $\hh$ satisfies \cref{asm:H-rich:a,asm:H-rich:b,asm:H-rich:c} and let $f$ satisfy \cref{asm:geom-f:b,asm:analytic-f}, then there exists $A_* \in \pdm(\hh)$ with $\rank(A_*) \leq d|Z| + 1$ such that $f(x) - f^* = \scal{\phi(x)}{A_* \phi(x)}$ for all $x \in \Omega$.
\ecor
\bpr
By \cref{thm:g-good-has-A} we know that if $f$ satisfies \cref{asm:geom-f:b,asm:analytic-f} w.r.t. a space $\hh$ that satisfies \cref{asm:H-rich:a,asm:H-rich:b,asm:H-rich:c}, there exists $w_1,\dots,w_q \in \hh$ with $q \leq d|Z| + 1$ such that $f(x) - f^* = \sum_{j \in [q]} w_j^2(x)$ for any $x \in \Omega$. 
Since $\hh$ is a reproducing kernel Hilbert space, for any $h \in \hh,~x \in \Omega$ we have $h(x) = \scal{\phi(x)}{h}_{\hh}$. Moreover, by the properties of the outer product in Hilbert spaces, for any $h,v \in \hh$, it holds $(\scal{h}{v}_{\hh})^2 = \scal{h}{(v \otimes_{\hh} v) h}$.

Thus, for any $x \in \Omega, j \in [q]$, it holds $w_j(x)^2 = \scal{\phi(x)}{(w_j \otimes w_j) \phi(x)}$ and hence
$$\forall x \in \Omega,~f(x) - f^* = \scal{\phi(x)}{A_* \phi(x)}, \qquad A_* = \sum_{j \in [q]} w_j \otimes w_j.$$
\epr

To conclude the section we prove the problem in \cref{eq:prob-intermediate} admits a maximizer whose non-negative operator is of rank at most $d|Z|+1$.

\bt\label{thm:Hilbert-isgood}
Let $\Omega \subset \R^d$ be an open set, $k$ be a kernel, $\hh$ the associated RKHS, and $f:\R^d \to \R$. Under  \cref{asm:geom-f,asm:H-rich,asm:analytic-f}, the problem in \cref{eq:prob-intermediate} admits an optimal solution $(c_*, A_*)$ with $c_* = f_*$,  and $A_*$ a positive operator on $\hh $ with rank at most $d|Z|+1$.
\et




\bpr
Let $p_0$ be the maximum of \cref{eq:prob-convex-version}.
Since $A \succeq 0$ implies $\scal{\phi(x)}{ A \phi(x)} \geq 0$ for all $x \in \Omega$, the problem in \cref{eq:prob-convex-version} is a relaxation of \cref{eq:prob-intermediate}, where the constraint $f(x) - c = \scal{\phi(x)}{ A \phi(x)}$ is substituted by $f(x) - c \geq 0, \forall x \in \Omega$.
Then $p_0 \geq p_*$ if a maximum $p_*$ exists for \cref{eq:prob-intermediate}. Moreover if there exists $A$ that satisfies the constraints in \cref{eq:prob-intermediate} for the value $c_* = f_*$, then $p_0 = p_*$ and $(c_*, A)$ is a maximizer of \cref{eq:prob-intermediate}.
The proof is concluded by  applying \cref{cor:g-good-has-A} that shows that there exists $A$ satisfying the constraints in \cref{eq:prob-intermediate} for the value $c =  f_*$. 
\epr

In \cref{cor:g-good-has-A} and \cref{thm:Hilbert-isgood} we proved the existence of an infinite-dimensional trace-class positive operator $A_*$ that satisfies $\scal{\phi(x)}{A_* \phi(x)} = f(x) - f_*$ for any $x \in \Omega$ and maximizing \cref{eq:prob-intermediate}. The proof is quite general, requiring some geometric properties on $f$, the fact that $f$ and its second derivatives belong to $\hh$ and some algebraic properties of the space~$\hh$, in particular to be closed to multiplication with a $C^\infty$ function, to integration, and  to composition with a $C^\infty$ map. The generality of the proof does not allow to derive an easy characterization of the trace of $A_*$.

\section{Properties of the finite-dimensional problem}\label{sec:finite-dim}

In the previous section we proved that there exists a finite rank positive operator $A_*$ minimizing \cref{eq:prob-intermediate}.
In this section we study the effect of the discretization of \cref{eq:prob-intermediate} on a given a set of distinct points $\widehat{X} = \{x_1,\dots,x_n\}$.
First, we derive \cref{thm:inequality-scattered-data} which is fundamental to prove \cref{thm:bound-prob-sampled}, and is our main technical result (we believe it can have a broader impact beyond the use in this paper as discussed in \cref{sec:technical-contr}). Given a smooth function $g$ on $\Omega$, in \cref{thm:inequality-scattered-data} we prove that if there exists a matrix $B \in \pdm(\R^n)$ such that $g(x_i) = \Phi_i^\top B \Phi_i$ for $i \in [n]$ (the vectors $\Phi_j \in \R^n$ are defined before \cref{eq:prob-relax}), then the inequality $g(x) \geq - \eps$ holds for any $x \in \Omega$ for an $\eps$ depending on the smoothness of the kernel, the smoothness of $g$ and how well the points in $\widehat{X}$ cover $\Omega$. We denote by $h_{\widehat{X},\Omega}$ the {\em fill distance} \cite{wendland2004scattered}, 
\eqal{\label{eq:fill-distance}
h_{\widehat{X},\Omega} = \sup_{x \in \Omega} \min_{i \in [n]} \|x - x_i\|,
}
corresponding to the maximum distance between a point in $\Omega$ and the set $\widehat{X}$.
In particular, if the kernel and $g$ are $m$-times differentiable,  \cref{thm:inequality-scattered-data} proves that $g(x) \geq -\eps$ holds with $\eps = O(h^m_{\widehat{X},\Omega})$ which is an improvement when $m \gg 2$ with respect to standard discretization results that guarantee exponents of only $1$ or $2$.
Then in \cref{lm:exists-B} we show that there exists a finite-dimensional positive definite matrix $B \in \pdm(\R^n)$ such that $\tr(B) \leq \tr(A_*)$ and  $\Phi_i^\top B \Phi_i = \scal{\phi(x_i)}{A_* \phi(x_i)}$ for all $i \in [n]$.
Finally, in \cref{thm:bound-prob-sampled}, we combine \cref{lm:exists-B} with \cref{thm:inequality-scattered-data}, to show that the problem in \cref{eq:prob-relax} provides a solution that is only $O(h^m_{\widehat{X},\Omega})$ distant from the solution of the infinite dimensional problem in \cref{eq:prob-intermediate}.

To start we recall some basic properties of $\Phi_i$ and $\phi(x_i)$, for $i \in [n]$, already sketched in \cref{sec:sketch}. In particular, the next proposition shows that, by construction, $\Phi_i^\top \Phi_j = \phi(x_i)^\top \phi(x_j)$ for any $i,j \in [n]$ and more generally that the map $V$ that maps $f \in \hh  \mapsto R^{-\top}(\scal{\phi(x_1)}{f},\dots,\scal{\phi(x_n)}{f}) \in \R^n$ is a partial isometry and that $\Phi_i = V\phi(x_i)$. The map $V$ will be crucial to characterize the properties of the finite dimensional version of the operator $A_*$


\blm[Characterizing $\Phi_j$ in terms of $\phi$]\label{prop:properties-V}
Let $k$ be a kernel satisfying \cref{asm:H-rich:a}. There exists a linear operator $V:\hh  \to \R^n$ such that
$$\Phi_i = V \phi(x_i), \qquad \forall i \in [n].$$
Moreover $V$ is a partial isometry: $VV^*$ is the identity on $\R^n$, $P = V^*V$ is a rank $n$ projection operator satisfying $P \phi(x_i) = \phi(x_i), \forall i \in [n]$.
\elm

The proof of \cref{prop:properties-V} is given in \cref{proof:prop:properties-V} in \cpageref{proof:prop:properties-V} and is based on the fact that the kernel matrix $K$ is positive definite and invertible when $k$ is {\em universal} \cite{steinwart2008support}, property that is implied by \cref{asm:H-rich:a}, and that $R$ is an invertible matrix that satisfies $K = R^\top R$.

\subsection{Uniform inequality from scattered constraints}\label{sec:inequality}

In this section we derive \cref{thm:inequality-scattered-data}. Here we want to guarantee that a function $g$ satisfies $g(x) \geq - \eps$ on $\Omega$, by imposing some constraints on $g(x_i)$ for $i \in [n]$. If we use the most natural discretization, that consists in the constraints $g(x_i) \geq 0$, by Lipschitzianity of $g$ we can guarantee only $\eps = |g|_{\Omega,1} h_{\widehat{X},\Omega}$ (recall the definition of $|\cdot|_{\Omega,m}$ for $m \in \Nz$ from \cref{eq:df_seminorm}).  In the case of {\em equality constraints}, instead, standard results for {\em functions with scattered zeros} \cite{wendland2004scattered} (recalled in \cref{sec:scattered-data}) guarantee for all $x \in \Omega$
$$|u(x)| \leq \eps, \quad \eps = C h_{\widehat{X},\Omega}^m |u|_{\Omega,m},$$ 
when $u$ is $m$-times differentiable and satisfies $u(x_i) = 0$ for any $i \in [n]$ (see \cite{wendland2004scattered,narcowich2003refined} or \cref{thm:zeros-scattered-data} for more details). Thus, in this case the discretization leverages the degree of smoothness of $u$, requiring much less points to achieve a given $\eps$ than in the inequality case. 

The goal here is to derive a guarantee for {\em inequality constraints} that is as strong as the one for the equality constraints. In particular, given a function $g$ defined on $\Omega$ and that satisfies $g(x_i) - \Phi_i B \Phi_i = 0$ on $\widehat{X}$, with $B \succeq 0$, we first derive a function $u$ defined on the whole $\Omega$ and matching $g(x_i) - \Phi_i B \Phi_i$ on $\widehat{X}$. This is possible since we know that $\Phi_i = V \phi(x_i)$, by \cref{prop:properties-V}, then $u(x) = g(x) - \scal{\phi(x)}{V^*BV\, \phi(x)}$ satisfies $u(x_i) = g(x_i) - \Phi_i B \Phi_i$ for any $i \in [n]$. Finally, we apply the results for functions with scattered zeros on $u$. 
The desired result is obtained by noting that, since $\scal{\phi(x)}{V^*BV\, \phi(x)} \geq 0$ for any $x \in \Omega$, by construction, then for all $x \in \Omega$
$$-g(x) \leq - g(x) + \scal{\phi(x)}{V^*BV\, \phi(x)} \leq |g(x) - \scal{\phi(x)}{V^*BV\, \phi(x)}| = |u(x)| \leq \eps,$$
i.e., $g(x) \geq - \eps$ for all $x \in \Omega$ with $\eps = C h_{\widehat{X},\Omega}^m |u|_{\Omega,m}$. In the following theorem we provide a slightly more general result, that allows for $|g(x_i) - \Phi_i B \Phi_i| \leq \tau$ with $\tau \geq 0$.

\bt[Uniform inequality from scattered constraints]\label{thm:inequality-scattered-data}
Let $\Omega$ satisfy \cref{asm:geom-f:a} for some $r > 0$. Let $k$ be a kernel satisfying \cref{asm:H-rich:a,asm:H-rich:d} for some $m \in \Np$.
Let $\widehat{X} = \{x_1,\dots, x_n\} \subset \Omega$ with $n \in \Np$ such that $h_{\widehat{X},\Omega} \leq r\min(1, \frac{1}{18(m-1)^2})$.
Let $g \in C^m(\Omega)$ and assume there exists $B \in \pdm(\R^n)$ and $\tau \geq 0$ such that 
\eqal{\label{eq:empirical-inequality}
|g(x_i) - \Phi_i^\top B \Phi_i| \leq \tau, \quad \forall i \in [n],
}
where the $\Phi_i$'s are defined in \cref{sec:sketch}. The following statement holds:
\eqal{
g(x) \geq -(\eps + 2\tau) ~~ \forall x \in \Omega, \quad \textrm{where} \quad \eps = C  h_{\widehat{X},\Omega}^m,
}
and $C = C_0(|g|_{\Omega,m}+\Cmul\Cdiff\tr(B))$ with $C_0 = 3\tfrac{\max(\sqrt{d},3\sqrt{2d}(m-1))^{2m}}{m!}$. The constants $m, \Cmul, \Cdiff$, defined in \cref{asm:H-rich:a,asm:H-rich:d}, do not depend on $n, \widehat{X}, h_{\widehat{X},\Omega}, B$ or $g$.
\et
\bpr
Let the partial isometry $V:\hh  \to \R^n$ and the projection operator $P = V^*V$ be defined as in \cref{prop:properties-V}. Given $B \in \pdm(\R^n)$ satisfying \cref{eq:empirical-inequality},
define the operator $A \in \pdm(\hh )$ as $A = V^* B V$ and the functions $u,r_A:\Omega \to \R$ as follows
$$ r_A(x) =\scal{\phi(x)}{A\phi(x)},\quad u(x) = g(x) - r_A(x), \quad \forall x \in \Omega.$$
Since $\Phi_i = V\phi(x_i)$ for all $i \in [n]$, then for all $i \in [n]$:
\eqals{
r_A(x_i) = \big<\phi(x_i), V^* B V\phi(x_i)\big>  = (V\phi(x_i))^\top B (V\phi(x_i)) = \Phi_i^\top B \Phi_i,
}
and hence $u(x_i) = g(x_i) - \Phi_i^\top B \Phi_i$. Thus, $|u(x_i)| \leq \tau$ for any $i \in [n]$. This allows to apply one of the classical results on functions with scattered zeros \cite{narcowich2003refined,wendland2004scattered} to bound $\sup_{x \in \Omega}|u(x)|$, which we derived again in \cref{thm:zeros-scattered-data} to obtain explicit constants. Since we have assumed $h_{\widehat{X},\Omega} \leq r/\max(1,18(m-1)^2)$, applying \cref{thm:zeros-scattered-data}, the following holds
$$\sup_{x \in \Omega} |u(x)| \leq 2\tau + \eps, \qquad  \eps = c~ R_m(u) \,\, h_{\widehat{X},\Omega}^m,$$
where $c = 3 \max(1,18(m-1)^2)^{m}$ and $R_m(v) = \sum_{|\alpha| = m}\frac{1}{\alpha!}\sup_{x \in \Omega}|\partial^\alpha v(x)|$ for any $v \in C^m(\Omega)$ using the multi-index notation (recalled in \cref{sec:notation}).
Since $r_A(x) = \scal{\phi(x)}{A\phi(x)} \geq 0$ for any $x \in \Omega$ as $A \in \pdm(\hh)$, it holds :
\begin{equation}
\label{eq:boundscattered} g(x) \geq  g(x) - r_A(x) = u(x) \geq -|u(x)| \geq -(2\tau + \eps), \qquad \forall x \in \Omega.
\end{equation}
The last step is bounding $R_m(u)$. Recall the definition of $|\cdot|_{\Omega,m}$ from \cref{eq:df_seminorm}.
First, note that $A = V^*BV$ is finite rank (hence trace-class). Applying the cyclicity of the trace and the fact that $VV^*$ is the identity on $\R^n$, it holds
$$\tr(A) = \tr(V^*BV) = \tr(BVV^*) = \tr(B).$$
Since $k$ satisfies \cref{asm:H-rich:a}, by \cref{lm:H-norm-tr-norm}, \cpageref{lm:H-norm-tr-norm}, $r_A \in \hh $ and $\|r_A\|_{\hh } \leq \Cmul \tr(A) = \Cmul \tr(B)$ where $\Cmul$ is fixed in \cref{asm:H-rich:a}.  Moreover, since the kernel $k$ satisfies \cref{asm:H-rich:d} with $m$ and $\Cdiff$, then $|v|_{\Omega,m} \leq \Cdiff \|v\|_{\hh }$, for any $v \in \hh$ as recalled in \cref{rm:k-differentiable}. In particular, this implies $|r_A|_{\Omega,m} \leq \Cdiff \|r_A\|_\hh \leq \Cdiff \Cmul \tr(B)$.
To conclude, note that, by the multinomial theorem,
\eqals{
R_m(u) & =  \sum_{|\alpha| = m}\frac{1}{\alpha!} \sup_{x \in \Omega} \big|\partial^\alpha u(x)\big| \leq \sum_{|\alpha| = m}\frac{1}{\alpha!} ~  |u|_{\Omega,m} ~=~ \frac{d^m}{m!} ~|u|_{\Omega,m}.
}
Since $|u|_{\Omega,m} \leq |g|_{\Omega,m} + |r_A|_{\Omega,m}$, combining all the previous bounds, it holds 
\[\eps \leq C_0 ~(|g|_{\Omega,m} + \Cdiff \Cmul \tr(B))~ h_{\widehat{X},\Omega}^m,\qquad C_0 = 3 \frac{d^m\max(1,18(m-1)^2)^m}{m!}.\]
The proof is concluded by bounding $\eps$ in \cref{eq:boundscattered} with the inequality above.
\epr

In the theorem above we used a domain satisfying \cref{asm:geom-f:a} and a version of a bound for functions with scattered zeros (that we derived in \cref{thm:zeros-scattered-data} following the analysis in \cite{wendland2004scattered}), to have explicit and relatively small constants. However, by using different bounds for functions with scattered zeros, we can obtain the same result as \cref{thm:inequality-scattered-data}, but with different assumptions on $\Omega$ (and different constants). For example, we can use Corollary 6.4 in \cite{narcowich2003refined}  to obtain a result that holds for $\Omega = [-1, 1]^d$ or Theorem 11.32 with $p=q=\infty, m=0$ in \cite{wendland2004scattered} to obtain a result that holds for $\Omega$ with locally Lipschitz-continuous boundary.

\subsection{Convergence properties of the finite-dimensional problem}

Now we use \cref{thm:inequality-scattered-data} to bound the error of \cref{eq:prob-relax}. First, to apply \cref{thm:inequality-scattered-data} we need to prove the existence of at least one finite-dimensional $B \succeq 0$ that satisfies the constraints of \cref{eq:prob-relax} and such that the trace of $B$ is independent of $n$ and $h_{\widehat{X},\Omega}$. This is possible since we proved in \cref{thm:Hilbert-isgood} that there exists at least one finite rank operator $A$ that solves \cref{eq:prob-intermediate} and thus satisfies its constraints, of which the ones in \cref{eq:prob-relax} constitute a subset. In the next lemma we construct $\overline{B} \in \pdm(\R^n)$, such that $\scal{\phi(x_i)}{A \phi(x_i)} = \Phi_i^\top \overline{B} \Phi_i$. In particular, 
$\overline{B} = V A_* V^* = R^{-\top} C R^{-1}$, with $C_{i,j} = \scal{\phi(x_i)}{A_* \phi(x_j)}$ for $i,j \in [n]$, where $A_*$ is one solution of \cref{eq:prob-intermediate} with minimum trace-norm, since the bound in \cref{thm:inequality-scattered-data} depends on the trace of the resulting matrix.

\blm\label{lm:exists-B}
Let $\Omega$ be an open set and $\{x_1,\dots,x_n\} \subset \Omega$ with $n \in \Np$. Let $g: \Omega \rightarrow \R$ and $k$ be a kernel on $\Omega$. Denote by $\hh$ the associated RKHS and by $\phi$ the associated canonical feature map. Let $A \in \pdm(\hh)$ satisfy $\tr(A) < \infty$ and $\scal{\phi(x)}{A\phi(x)} = g(x),~ x \in \Omega$. Then there exists $\overline{B} \in \pdm(\R^n)$ such that $\tr(\overline{B}) \leq \tr(A)$ and $g(x_i) = \Phi_i^\top \overline{B} \Phi_i, ~ \forall i \in [n]$.
\elm

\bpr
Let $V:\hh  \to \R^n$ be the partial isometry defined in \cref{prop:properties-V} and $P = V^*V$ be the associated projection operator. Define $\overline{B} \in \R^{n\times n}$ as $\overline{B} = V A V^*$. Since by \cref{prop:properties-V}, $\Phi_i = V \phi(x_i)$ and $P$ satisfies $P\phi(x_i) = \phi(x_i)$ for $i \in [n]$,
\eqals{
\Phi_i^\top \overline{B} \Phi_i &= (V \phi(x_i))^\top(V A V^*)(V \phi(x_i)) = \big<V^*V\phi(x_i), A V^*V \phi(x_i)\big>\\ & = \big<P\phi(x_i), A P \phi(x_i)\big> = \scal{\phi(x_i)}{A\phi(x_i)}\quad \forall i \in [n].
}
Note that $\overline{B}$ satisfies: (a) $\overline{B} \in \pdm(\R^n)$, by construction; (b) the requirement $\Phi_i^\top \overline{B} \Phi_i = g(x_i)$, indeed $\Phi_i^\top \overline{B} \Phi_i = \scal{\phi(x_i)}{A \phi(x_i)}$ and $\scal{\phi(x)}{A \phi(x)} = g(x)$ for any $x \in \Omega$; (c) $\tr(\overline{B}) \leq \tr(A)$, indeed, by the cyclicity of the trace, 
\eqals{
\tr(\overline{B}) = \tr(V A V^*) = \tr(AV^*V) = \tr(AP).
}
The proof is concluded by noting that, since $A \succeq 0$ and $\|P\|_{\rm op} \leq 1$ because $P$ is a projection, then $\tr(A P) \leq \|P\|_{\rm op}\tr(|A|) = \|P\|_{\rm op}\tr(A) \leq \tr(A)$.
\epr

We are now ready to prove the convergence rates of \cref{eq:prob-relax} to the global minimum. We will use the bound for the inequality on scattered data that we derived \cref{thm:inequality-scattered-data} and the fact that there exists $\overline{B} \succeq 0$ that satisfies the constraints of \cref{eq:prob-relax} with a trace bounded by $\tr(A_*)$ as we proved in the lemma above (that is in turn bounded by the the trace of the operator explicitly constructed in \cref{thm:g-good-has-A}). The proof is organized as follows. We will first show that \cref{eq:prob-relax} admits a minimizer, that we denote by $(\hat{c}, \hat{B})$. The existence of $\overline{B}$ allows to derive a lower-bound on $\hat{c} - f_\ast$. Using \cref{thm:inequality-scattered-data} on the constraints of \cref{eq:prob-relax} and evaluating the resulting inequality in one minimizer $\zeta$ of $f$ allows to find an upper bound on $\hat{c} - f_\ast$ and an upper bound for $\tr(\hat{B})$.

\bt[Convergence rates of \cref{eq:prob-relax} to the global minimum]\label{thm:bound-prob-sampled} 
Let $\Omega$ be a set satisfying \cref{asm:geom-f:a} for some $r > 0$. Let $n\in \Np$ and $\widehat{X} = \{x_1,\dots,x_n\} \subset \Omega$ with fill distance $h_{\widehat{X},\Omega}$. Let $k$ be a kernel and $\hh$ the associated RKHS satisfying \cref{asm:H-rich} for some $m \in \Np$. Let $f$ be a function satisfying \cref{asm:geom-f:b} and \cref{asm:analytic-f} for $\hh$. The problem in \cref{eq:prob-relax} admits a solution.
Let $(\hat{c}, \hat{B})$ be any solution of \cref{eq:prob-relax}, for a given $\la > 0$. The following holds
\eqal{\label{eq:appr-c}
|\hat{c} - f_*| ~\leq~ 2\eta\,|f|_{\Omega,m} ~+~ \la\,\tr(A_*), \qquad \eta ~=~ C_0 \,h^m_{\widehat{X},\Omega},
}
when $h_{\widehat{X},\Omega} \leq r\min(1, \frac{1}{18(m-1)^2})$ and $\la \geq 2 \Cmul\Cdiff \eta $. Here
 $C_0 = 3\tfrac{\max(\sqrt{d},3\sqrt{2d}(m-1))^{2m}}{m!}$, $\Cdiff, \Cmul$ are defined in \cref{asm:H-rich} and $A_*$ is given by \cref{thm:Hilbert-isgood}. Moreover, under the same conditions
\eqal{\label{eq:appr-A}
\tr(\hat{B}) ~\leq~ 2\,\tr(A_*) ~+~ 2\tfrac{\eta}{\la}\,|f|_{\Omega,m}.
}
\et
\bpr
We divide the proof in few steps. 

\noindent{\bf Step 0. Problem \cref{eq:prob-relax} admits always a solution.} 

\noindent (a) On the one hand, \textit{$c$ cannot be larger than $c_0 = \min_{i \in [n]} f(x_i)$}, otherwise there would be a point $x_j$ for which $f(x_j) - c < 0$ and so the constraint $\Phi_j^\top B \Phi_j = f(x_j) - c$ would be violated, since does not exist any positive semi-definite matrix for which $\Phi_j^\top B \Phi_j < 0$.

\noindent (b) On the other, \textit{there exists an admissible point.} Indeed let $(c_*, A_*)$ be the solution of \cref{eq:prob-intermediate} such that $A_*$ has minimum trace norm. By \cref{thm:Hilbert-isgood}, we know that this solution exists with $c_* = f_\ast$, under \cref{asm:geom-f,asm:H-rich,asm:analytic-f}. Then, by \cref{lm:exists-B} applied to $g(x) = f(x)-c_*$ and $A = A_*$, given $\widehat{X} = \{x_1,\dots,x_n\}$ we know that there exists $\overline{B} \in \pdm(\R^n)$ satisfying $\tr(\overline{B}) \leq \tr(A_*)$ such that the constraints of \cref{eq:prob-relax} are satisfied for $c = c_*$. Then $(c_*, \overline{B})$ is admissible for the problem in \cref{eq:prob-relax}. 

\noindent Thus, since there exists an admissible point for the constraints of \cref{eq:prob-relax} and its functional cannot be larger than $c_0$ without violating one constraint, the SDP problem in \cref{eq:prob-relax} admits a solution~(see \cite{boyd2004convex}).

\noindent{\bf Step 1. Consequences of existence of $A_*$.} 
Let $(\hat{c}, \hat{B})$ be one minimizer of \cref{eq:prob-relax}. The existence of the admissible point $(c_*, \overline{B})$ proven in the step above implies that
$$
\hat{c} - \la \tr(\hat{B}) \geq c_* - \la\tr(\overline{B}) \geq f_* - \la\tr(A_*),
$$
from which we derive,
\eqal{\label{eq:lb-cc}
\la\tr(\hat{B}) - \la\tr(A_*) \leq \Delta, \quad \Delta := \hat{c} - f_*.
}

\noindent{\bf Step 2. $\boldsymbol{f|_\Omega \in C^{m+2}(\Omega)}$.} \cref{asm:analytic-f} guarantees that $f|_\Omega \in C^2(\Omega)$ and that for all $i,j \in [d]$, $\tfrac{\partial}{\partial x_i \partial x_j} f|_\Omega \in \hh$. Since under \cref{asm:H-rich:d}, $\hh \subset C^m(\Omega)$ by \cref{rm:k-differentiable}, we see that $\tfrac{\partial}{\partial x_i \partial x_j} f|_\Omega \in C^m(\Omega)$ for all $i,j \in [d]$ and hence $f|_\Omega \in C^{m+2}(\Omega)$. 

\noindent{\bf Step 3. $L^\infty$ bound due to the scattered zeros.} Let $(\hat{c}, \hat{B})$ be one minimizer of \cref{eq:prob-relax} and define $\hat{g}(x) = f(x) - \hat{c}$ for all $x \in \Omega$.
Note that $\hat{g}(x_i) = \Phi_i^\top \hat{B} \Phi_i$ for $i \in [n]$. Moreover, $\hat{g} \in C^m(\Omega)$ because $f \in C^m(\Omega)$ and $\hat{c}$ is a constant. Considering that $h_{\widehat{X},\Omega} \leq \frac{r}{\max(1,18(m-1)^2)}$, by assumption, then all the conditions in \cref{thm:inequality-scattered-data} are satisfied for $g = \hat{g}$, $\tau = 0$ and $B = \hat{B}$. Applying \cref{thm:inequality-scattered-data}, we obtain, 
$$\forall x \in \Omega,~ f(x) - \hat{c} = \hat{g}(x) \geq - \eta (|\hat{g}|_{\Omega,m} + \Cmul\Cdiff\tr(\hat{B})), \qquad \eta = C_0 h^m_{\widehat{X},\Omega},$$
where $C_0$ is defined in \cref{thm:inequality-scattered-data}. Since the inequality above holds for any $x \in \Omega$, by evaluating it in one global minimizer $\zeta \in \Omega$, we have $f(\zeta) = f_*$ and hence
$$-\Delta = f_* - \hat{c} = f(\zeta) - \hat{c} = \hat{g}(\zeta) \geq -\eta (|\hat{g}|_{\Omega,m} + \Cmul\Cdiff\tr(\hat{B})).$$
Since $\hat{g} = f - \hat{c} \boldsymbol{1}_\Omega$, and since for any $m \in \Np$, $|\boldsymbol{1}_{\Omega}|_{\Omega,m} = 0$,  we have $|\hat{g}|_{\Omega,m} \leq |f|_{\Omega,m} + |\boldsymbol{1}_{\Omega}|_{\Omega,m} = |f|_{\Omega,m}$. Injecting this in the previous equation yields
\eqal{\label{eq:key-inequality}
\Delta \leq \eta |f|_{\Omega,m}  + \eta \Cmul \Cdiff \tr(\hat{B}).
}
\noindent{\bf Conclusion.}
Combining \cref{eq:key-inequality} with \cref{eq:lb-cc}, and since $\la \geq 2\Cmul \Cdiff \eta$ by assumption, 
$$\tfrac{\la}{2} \tr(\hat{B}) \leq (\la - \Cmul \Cdiff \eta) \tr(\hat{B}) \leq  \eta |f|_{\Omega,m} + \la\tr(A_*).$$
Note that \cref{eq:appr-A} is obtained from the one above, by dividing by $\tfrac{\lambda}{2}$. Finally the inequality \cref{eq:appr-c} is derived by bounding $\Delta$ from below as $\Delta \geq -\la \tr(A_*)$ by \cref{eq:lb-cc}, since $\tr(\hat{B}) \geq 0$ by construction, and bounding it from above as
$$\Delta  \leq 2\eta |f|_{\Omega,m} + \la\tr(A_*),$$
obtained by combining \cref{eq:key-inequality} with \cref{eq:appr-A} and with the assumption $\Cmul \Cdiff \eta \leq \frac{\la}{2}$.
\epr

The result above holds for any kernel satisfying \cref{asm:H-rich} and any function $f, \Omega$ satisfying the geometric conditions in \cref{asm:geom-f} and with $f \in C^2(\Omega)$ and $\frac{\partial^2 f}{\partial x_i \partial x_j} \in \hh$ for $i,j \in [d]$. The latter requirement is quite easy to verify for example when $\hh$ contains $C^s(\Omega)$ and $f \in C^{s+2}(\Omega)$ for some $s > 0$ as in the case of $\hh$ being a Sobolev space with $s > d/2$. Moreover the proposed result holds for any discretization $\widehat{X}$ (random, or deterministic). We would like to conclude with the following remark on the sufficiency of the assumptions on $f$.

\br[Sufficiency of \cref{asm:geom-f:b,asm:analytic-f}]\label{rem:sufficiency-asm-f}
Assumptions \cref{asm:geom-f:b,asm:analytic-f} are sufficient for \cref{thm:Hilbert-isgood,thm:bound-prob-sampled} to hold. However, by inspecting their proof it is clear that they hold by requiring only the existence of a trace-class operator $A_* \in \pdm(\hh)$ such that $f(x) - f_* = \scal{\phi(x)}{A_* \phi(x)}$ for any $x \in \Omega$, where $f_* = \inf_{x \in \Omega} f(x)$. Note that this is implied by \cref{asm:geom-f:b,asm:analytic-f} via \cref{cor:g-good-has-A}.
\er

In the next subsection we are going to apply the theorem above to the specific setting of \cref{alg:glm}.

\subsection{Result for Sobolev kernels and discussion}

In this we are going to apply \cref{thm:bound-prob-sampled} to \cref{alg:glm} which corresponds to $\hh$ be the Sobolev space of smoothness $s$ and the points $\widehat{X}$ selected independently and uniformly at random. First, in the next lemma we bound in high probability the fill distance $h_{\widehat{X},\Omega}$ with respect to the number of points $n$ that we sample, i.e.,  the cardinality of $\widehat{X}$.

\blm[Random sets of points]\label{thm:fill-distance-random-points}
Let $\Omega \subset \R^d$ be a bounded set with diameter $2R$, for some $R > 0$, and satisfying \cref{asm:geom-f:a} for a given $r > 0$. Let $\widehat{X} = \{x_1,\dots, x_n\}$ independent points sampled from the uniform distribution on $\Omega$. When $n \geq 2(\frac{6R}{r})^d \left(\log\frac{2}{\delta} + 2d \log\frac{4R}{r}\right)$, then the following holds with probability at least $1-\delta$:
$$h_{\widehat{X},\Omega} ~~~\leq ~~~ 11R ~ n^{-\frac{1}{d}} ~ (\log \tfrac{n}{\delta} + d \log\tfrac{2R}{r})^{{1}/{d}}.$$
\elm


The proof of \cref{thm:fill-distance-random-points} is in \cref{proof:thm:fill-distance-random-points},  \cpageref{proof:thm:fill-distance-random-points} and is a simpler version (with explicit constants) of more general results \cite[ Thm. 13.7]{penrose2003random}. In the next theorem we apply the bound in the lemma above with the explicit constants for Sobolev spaces derived in \cref{prop:sobolev-kernel} to \cref{thm:bound-prob-sampled}. The derivation of the theorem below is in \cref{proof:thm:alg-glm}, \cpageref{proof:thm:alg-glm}.

\bt[Convergence rates of \cref{alg:glm} to the global minimum]\label{thm:alg-glm} 
Let $\Omega \subset \R^d$ be a bounded set with diameter $2R$, for some $R > 0$, and satisfying \cref{asm:geom-f:a} for a given $r \in (0,R]$ (e.g. if $\Omega$ is a ball with radius $R$, then $r = R$).  Let $s$ satisfying $s > d/2$. Let $k$ be Sobolev kernel of smoothness $s$ (see \cref{ex:sobolev-kernel}). Assume that $f$ satisfies \cref{asm:geom-f:b} and that $f|_\Omega \in W^{s+2}_2(\Omega)$. Let $\hat{c}$ be the result of \cref{alg:glm} executed with $n \in \Np$ points chosen uniformly at random in $\Omega$ and $\la > 0$.
Let $\delta \in (0,1]$. When $m \in \Np$ satisfies $m < s-d/2$ and $n \geq \max(4,15(m-1))^{2d}\left(\tfrac{R}{r}\right)^d\left(2 \log \frac{2}{\delta} +  4d \log \tfrac{20 R~m }{r}\right)$ choose any $\la$ satisfying
$$\la ~~~\geq~~~ n^{-\tfrac{m}{d}} ~~ (\log \tfrac{2^d n}{\delta})^{\tfrac{m}{d}}~R^m C_{m,s,d},$$
where $C_{m,s,d} =  11^m C_0 \max(1,\Cmul \Cdiff)$ with $C_0$ defined in \cref{thm:bound-prob-sampled} and $\Cmul \Cdiff$ defined in \cref{prop:sobolev-kernel}. Note that $C_{m,s,d}$ is explicitely bounded in the proof in \cpageref{proof:thm:alg-glm} in terms of $s,m,d$.
Then, with probability at least $1-\delta$, the following holds
$$|\hat{c} - f_*| ~~\leq~~ 3\,\la\,\left(\tr(A_*) + |f|_{\Omega,{m}}\right).$$
\et


A direct consequence of the theorem above, already stated in \cref{rm:sobolev-nearly-optimal}, is the nearly-optimality of \cref{alg:glm} for the cases of Sobolev functions. Indeed by applying \cref{thm:alg-glm} with $m$ equal to the largest integer strictly smaller than $s-d/2$ we have that $m \geq s-d/2-1$, and so \cref{alg:glm} achieves the global minimum with a rate that is $O(n^{-\frac{s}{d}+\frac{1}{2}+\frac{1}{d}})$. The lower bounds from information based complexity state that, by observing the functions in $n$ points, it is not possible to find the minimum with error smaller than $n^{-\frac{s}{d} + \frac{1}{2}}$ for functions in $W^s_2(\Omega)$  (see, e.g., \cite{novak2006deterministic}, Prop.~1.3.11, page 36). Since in \cref{thm:alg-glm} we assume $f$   belongs to $W^{s+2}_2(\Omega)$, the optimal rate would be $n^{-\frac{s}{d} + \frac{1}{2} - \frac{2}{d}}$ so we are a factor $n^{3/d}$ slower than the optimal rate. Note that this factor is negligible if the function is very smooth, i.e., $s \gg d$, or $d$ is very large. An interesting corollary that corresponds to \cref{thm:informal}, can be derived considering that $C^{s+2}(\Omega) \subseteq W^{s+2}_2(\Omega)$, since $\Omega$ is bounded.

\section{Algorithm}
\label{sec:algorithm}

We need to solve the following optimization problem:
$$
\max_{ B \succcurlyeq 0,  c \in \R}
c - \ {\lambda}\tr(B) ~~~\mbox{ such that }~~~   f(x_i)   -  c - \Phi_i ^\top B \Phi_i = 0,~~\forall i \in [n].
$$
This is a semi-definite programming problem with $n$ constraints and a semi-definite constraint of size $n$. It can thus be solved with precistion $\eps$ in time $O(n^{3.5} \log ( 1/\eps))$ and memory $O(n^2)$ by standard software packages~\cite{boyd2004convex}. However, to allow applications to $n=1000$ or more, and on parallel architectures, we provide a simple Newton algorithm, which relies on penalization by a self-concordant barrier, that is, we aim to solve
$$
\max_{ B \succcurlyeq 0,  c \in \R}
c - \ {\lambda}\tr(B) + \frac{\eps}{n} \log \det (B) ~~~\mbox{ such that }~~~ f(x_i)   -  c - \Phi_i ^\top B \Phi_i = 0, ~~ \forall i \in [n],
$$
for which we know that at optimum, the deviation with the optimal value is at most $\eps$~\cite[Sec.~4.4]{nemirovski2004interior}.
By standard Lagrangian duality, we get, with $\Phi \in \R^{n \times n}$ the matrix with rows $\Phi_1,\dots,\Phi_n$, so that $\Phi \Phi^\top = K$:
\begin{eqnarray*}
 &  & \sup_{ B \succcurlyeq 0, c} \inf_{\alpha \in \mathbb{R}^n}
c  + \sum_{i=1}^n \alpha_i \big( f(x_i)  - c - \Phi_i^\top B \Phi_i \big)   - {\lambda}  \tr(B)
 + \frac{\eps}{n} \log \det (B)\\
& = &  \inf_{\alpha \in \mathbb{R}^n}
 \sum_{i=1}^n \alpha_i     f(x_i)  - \frac{\eps}{n} \log \det \big(\Phi^\top {\rm Diag}(\alpha) \Phi + \lambda I 
 \big)  + \frac{\eps}{n} \log \frac{  \eps}{n} - \eps  \ \mbox{ s. t. } \  \alpha^\top 1_n = 1   .
\end{eqnarray*} 
With the barrier term, this thus defines a dual function $H(\alpha)$, and we get the following gradient
$$
H'(\alpha)_i  = f_i - \frac{\eps}{n}\ \Phi_i^\top \Big( \Phi^\top {\rm Diag}(\alpha) \Phi + \lambda I   \Big)^{-1} \Phi_i = f_i - \frac{\eps}{n \alpha_i} \Big[
K(K+\lambda {\rm Diag}(\alpha)^{-1})^{-1}
\big]_{ii},$$
and Hessian
$$
H''(\alpha)_{ij} = \frac{\eps}{n} \big[ \Phi_i^\top \Big( \Phi^\top {\rm Diag}(\alpha) \Phi + \lambda I  \Big)^{-1} \Phi_j \big]^2
,$$
which can be rewritten
$$ H''(\alpha)_{ij} = \frac{\eps}{n  \alpha_j \alpha_i} \Big[
K (K+\lambda {\rm Diag}(\alpha)^{-1})^{-1}
\big]_{ij} \Big[
K (K+\lambda {\rm Diag}(\alpha)^{-1})^{-1}
\big]_{ji}.
$$
We can then compute
 the step for the Damped Newton algorithm:
$\alpha^+ =  \alpha - \tfrac{1}{1+ \sqrt{\frac{n}{\eps}}\lambda(\alpha)} \Delta$,
where $\Delta = H''(\alpha)^{-1} H'(\alpha) - \tfrac{1_n^\top H''(\alpha)^{-1}H'(\alpha)}{1_n^\top H''(\alpha)^{-1} 1_n } H''(\alpha)^{-1}1_n$ and $\lambda(\alpha)^2 = \Delta^\top H''(\alpha) \Delta $ is the Newton decrement (which can serve as a stopping criterion). 
Note that the algorithm is always feasible, without a need for any eigenvalue decomposition. The overall complexity is $O(n^3)$ per iteration due to matrix inversions and linear systems. Note that the conditioning of these linear systems is at least as bad as the conditioning of the kernel matrix $K$. Fortunately, for the $s$-th Sobolev kernels in dimension $d$, the $m$-th eigenvalue of the kernel matrix typically decay as $m^{-2s/d}$~\cite[Sec.~2.3]{bach2017breaking}.

\paragraph{Retrieving $c$ and $B$.} From an optimal $\alpha$, we can recover   $B = \frac{\eps}{n} \big( \Phi^\top {\rm Diag}(\alpha) \Phi + \lambda I   \big)^{-1} = \frac{\eps}{n \lambda} \big(
I - \Phi^\top ( \Phi \Phi^\top + \lambda {\rm  Diag}(\alpha)^{-1})^{-1} \Phi
\big)$ and $c = \frac{1}{n} H'(\alpha)^\top 1_n$ (since $c$ is the Lagrange multiplier for the constraint $\alpha^\top 1_n = 1 $). Thus, computing the model for a test point, can be done as $\frac{\eps}{n \lambda} \big(
k(x,x) - q(x)^\top ( K + \lambda {\rm  Diag}(\alpha)^{-1})^{-1} q(x)\big)$, where $q(x)_i = k(x,x_i)$. Alternatively, when $\Phi$ is invertible, we can use $q(x)^\top \Phi^{-\top} B \Phi^{-1} q(x)$.

\paragraph{Retrieving a minimizer.} Given the dual solution, based on our localizing arguments presented in \cref{sec:minimizer}, a good candidate solution will be 
\eqal{\label{eq:primal-candidate}
\hat{z} = \sum_{i=1}^n \alpha_i x_i
}
A more principled way to find a minimizer is provided in \cref{sec:minimizer}, of which the equation above corresponds to the limit solution of \cref{eq:prob-relax-z} for $\nu \to 0$ (see \cref{sec:minimizer-algorithm}).

\paragraph{Number of iterations.}
In order to reach a Newton decrement 
$n^{1/2}\eps^{-1/2}\lambda(\alpha) \leqslant \kappa$, a number of steps equal to a universal constant times $ \frac{n}{\eps} [ H(\alpha_0) - H(\alpha_\ast)]  +  \log \log \frac{1}{\kappa}$ is sufficient.~\cite{nemirovski2004interior}.

When initializing with $\alpha_0 = \frac{1}{n}1_n$, we have $H(\alpha_0) = \frac{1}{n} \sum_{i=1}^n f_i 
  - \frac{\eps}{n} \log \det \big(  K   + n \lambda I 
 \big)  + \frac{\eps}{n} \log  \eps - \eps,$ and $H(\alpha_\ast) \geqslant c_\ast - \lambda \tr (A_\ast) - \eps$. This leads to a number of Newton steps less than 
 $$
 \frac{n}{\eps}
 \big[ \langle f \rangle - \inf f \big] + 
 \log \det \big(  K   + n \lambda I 
 \big) + \frac{n}{\eps} \lambda \tr (A_\ast) + \log {\eps} +  \log \log \frac{1}{\kappa}.
 $$
 In our experiments, we do not perform path following (that would lead the classical interior-point method) and instead fixed value $\eps = 10^{-3}$, and a few hundred Newton steps. 

\paragraph{Behavior for $\lambda=0$.} If the kernel matrix $K$ is invertible (which is the case almost surely for Sobolev kernels and points sampled independently from a distribution with a density with respect to the Lebesgue measure), then we show that for $\lambda=0$, the optimal value of of the finite-dimensional problem in \cref{eq:prob-relax} is equal to $\min_{i \in [n]} f(x_i)$. Since $f(x_i) \geq c + \Phi_i^\top B \Phi_i$ implies $f(x_i) \geq c$, the optimal value has to be less than $\min_{i \in [n]} f(x_i)$. We therefore just need to find a feasible $B$ that achieves it. Since $K$ is assumed invertible (and thus its Cholesky factor as well), we can simply take $B = R^{-\top} {\rm Diag}[ (f(x_j) - \min_{i \in [n]} f(x_i))_j] R^{-1}$.

\section{Finding the global minimizer}
\label{sec:minimizer}

In this section we provide and study the problem in \cref{eq:prob-relax-z}, that is a variation of the problem in \cref{eq:prob-relax}, and allows to find also the minimizer of $f$ as we prove in \cref{thm:bound-prob-sampled-z}. As in \cref{sec:sketch} we start from a convex representation of the optimization problem and then we derive our sampled version, passing by an intermediate infinite-dimensional problem that is useful to derive the theoretical properties of the method. While the problem in \cref{eq:prob-convex-version} can be seen as finding the largest constant $c$ such that $f - c$ is still non-negative, in the problem below we find the parabola of the form $p_{z,\gamma}(x) = \frac{\nu}{2}\|x\|^2 - \nu x^\top z +  c = \frac{\nu}{2} \| x - z\|^2 + c - \frac{\nu}{2} \| z\|^2$ with the highest vertex such that $f - p_{z,c}$ is still non-negative. Since the height of the vertex of $p_{z,c}$ corresponds to $c - \frac{\nu}{2}\|z\|^2$, the resulting optimization problem is the following,
\eqal{\label{eq:prob-conv-z}
\max_{c \in \R, z \in \R^d} & c - \tfrac{\nu}{2} \|z\|^2  ~~~\mbox{ such that } ~~~ f(x) - \tfrac{\nu}{2}\|x\|^2 + \nu x^\top z - c \geq 0 \quad \forall x \in \Omega.
}
It is easy to see that if $f \in C^2(\R^d)$ has a unique minimizer $\zeta$ that belongs to $\Omega$ and is locally strongly convex around $\zeta$ then there exists a $\nu > 0$ such that the problem above achieves an optimum $(c_*,z_*)$ with $z_* = \zeta$ and $c_* = f_* + \tfrac{\nu}{2} \|\zeta\|^2$. In particular, to characterize $\nu$ explicitly we introduce the stronger assumption below.
\ba[Geometric assumption to find global minimizer]\label{asm:unique-minimizer}
The function $f: \R^d \to \R$ has a unique global minimizer in $\Omega$. 
\ea
If $f$ satisfies \cref{asm:geom-f:b,asm:unique-minimizer}, denote with $\zeta$ the unique minimizer of $f$ in $\Omega$ and with $f_* = f(\zeta)$ the corresponding minimum. 
\br\label{rm:existence-beta}
Under \cref{asm:geom-f:b,asm:unique-minimizer} $f$ can be lower bounded by a parabola with value $f_*$ at $\zeta$, i.e., there exists $\beta > 0$ such that 
\eqal{\label{eq:beta}
\forall x \in \Omega,~f(x) - f_* \geq \tfrac{\beta}{2} \|x - \zeta\|^2.
}
\er
The remark above is derived in  \cref{proof:rm:existence-beta}, \cpageref{proof:rm:existence-beta}.
In what follows, whenever $f$ satisfies \cref{asm:geom-f:b,asm:unique-minimizer}, then $\beta$ will be assumed to be the supremum among the value satisfying \cref{eq:beta}.  Now we are ready to summarize the reasoning above on the fact that \cref{eq:prob-conv-z} achieves the minimizer of $f$. 
\blm\label{lm:conv-z-is-ok}
Suppose $f$ satisfies \cref{asm:geom-f,asm:unique-minimizer}. Let $\zeta$ be the unique minimizer of $f$ in $\Omega$ and $f_* = f(\zeta)$ be the corresponding minimum. Let $\beta >0$ such that  \cref{eq:beta} holds.  If $\nu < \beta$ then the problem in \cref{eq:prob-conv-z} has a unique solution $(c_*, z_*)$ such that $z_* = \zeta$ and $c_* = f_* + \tfrac{\nu}{2} \|\zeta\|^2$. 
\elm
The lemma above guarantees that the problem in \cref{eq:prob-conv-z} achieves the global minimum and the global minimizer of $f$, when $f$ satisfies the geometric conditions \cref{asm:geom-f,asm:unique-minimizer}. Now, as we did for \cref{eq:prob-convex-version}, we consider the following problem of which \cref{eq:prob-conv-z} is a tight relaxation. 
\eqal{\label{eq:prob-intermediate-z}
\begin{split}
\max_{c \in \R, z \in \R^d, A \in \pdm(\hh )} & c - \tfrac{\nu}{2} \|z\|^2  \\
~~~\mbox{such that}~~~ &  f(x) - \tfrac{\nu}{2}\|x\|^2 + \nu x^\top z - c ~  = \scal{\phi(x)}{A\phi(x)} \quad \forall x \in \Omega.
\end{split}
}
Indeed, since $\scal{\phi(x)}{A\phi(x)} \geq 0$ for any $x \in \Omega$ and $A \in \pdm(\hh)$, for any triplet $(c,z,A)$ satisfying the constraints in the problem above, the couple $(c,z)$ satisfies the constraints in \cref{eq:prob-conv-z}. The contrary may be not true in general. In the next theorem we prove that when $\hh$ satisfies \cref{asm:H-rich} and $\Omega, f$ satisfy \cref{asm:geom-f,asm:analytic-f,asm:unique-minimizer}, then the relaxation is tight and, in particular, when $\nu < \beta$, there exists a finite rank operator $A_*$ such that the triplet $(f_* + \tfrac{\nu}{2}\|\zeta\|^2, \zeta, A_*)$ is optimal.
\bt\label{thm:Hilbert-z-isgood} 
Let $\Omega \subset \R^d$ be an open set, $k$ be a kernel, $\hh$ the associated RKHS, and $f:\R^d \to \R$ satisfying  \cref{asm:geom-f,asm:H-rich,asm:analytic-f}, and \cref{asm:unique-minimizer}. Let $\beta$ satisfying \cref{eq:beta}. For any $\nu < \beta$, the problem in \cref{eq:prob-intermediate-z} admits an optimal solution $(c_*, z_*, A_*)$ with $c_* = f_* + \frac{\nu}{2}\|\zeta\|^2$, $z_* = \zeta$,  and $A_*$ a positive semi-definite operator on $\hh$ with rank at most $d+1$.
\et
The proof of the theorem above is essentially the same of \cref{thm:Hilbert-isgood} and is reported for completeness in \cref{proof:thm:Hilbert-z-isgood}, \cpageref{proof:thm:Hilbert-z-isgood}. In particular, to prove the existence of $A_*$ we applied \cref{cor:g-good-has-A} to the function $f(x) - \frac{\nu}{2}\|x-\zeta\|^2$ that still satisfies \cref{asm:geom-f,asm:analytic-f} when $f$ does and $\nu < \beta$. Now we are ready to consider the finite-dimensional version of \cref{eq:prob-intermediate-z}. Given a set of points $\widehat{X} = \{x_1,\dots, x_n\}$ with $n \in \Np$, 
\eqal{\label{eq:prob-relax-z}
\begin{split}
\max_{c \in \R, z \in \R^d, B \in \pdm(\R^n)} & ~c - \tfrac{\nu}{2}\|z\|^2 - \la \tr(B)\\
 \mbox{such that} & \quad\forall i \in [n], \ f(x_i) - \tfrac{\nu}{2}\|x_i\|^2 + \nu x_i^\top z - c = \Phi_i^\top B \Phi_i.
 \end{split}
}
For the problem above we can derive similar convergence guarantees as for \cref{eq:prob-relax} and also a convergence of the estimated minimizer $z$ to $\zeta$, as reported in the following theorem.
\bt[Convergence rates of \cref{eq:prob-relax-z} to the global minimizer]\label{thm:bound-prob-sampled-z}
Let $\Omega$ be a set satisfying \cref{asm:geom-f:a} for some $r > 0$. Let $\widehat{X} = \{x_1,\dots,x_n\} \subset \Omega$ with fill distance $h_{\widehat{X},\Omega}$. Let $k$ be a kernel satisfying \cref{asm:H-rich} for some $m \geq 2$ and $f$ satisfying \cref{asm:geom-f,asm:analytic-f,asm:unique-minimizer}. The problem in \cref{eq:prob-relax-z} admits a solution. Denote by $(\hat{c}, \hat{z}, \hat{B})$ any solution of \cref{eq:prob-relax-z}, for a given $\la > 0$. Then
\eqal{\label{eq:appr-z}
\tfrac{\nu}{2}\|\hat{z} - \zeta\|^2 \leq 3\eta (|f|_{\Omega,m} + \nu)  + 2\la \tr(A_*), \qquad \eta ~=~ C \,h^m_{\widehat{X},\Omega},
}
when $h_{\widehat{X},\Omega} \leq \frac{r}{18(m-1)^2}$ and $\la \geq 2\Cmul\Cdiff \eta$. Here $C = 3\frac{(3\sqrt{2d}(m-1))^{2m}}{m!}$ and $\Cdiff, \Cmul$ are defined in \cref{asm:H-rich}. $A_*$ is from \cref{thm:Hilbert-z-isgood}. Moreover under the same conditions
\eqal{
\label{eq:appr-c-z} |\hat{c} - \tfrac{\nu}{2}\|\hat{z}\|^2 - f_*| ~&\leq~ 2\eta\,|f|_{\Omega,m} ~+~ \la\,\tr(A_*) ~+~ 2\eta \nu, \\
\label{eq:appr-A-z} \tr(\hat{B}) ~& \leq~ 2\,\tr(A_*) ~+~ 2\tfrac{\eta}{\la}\,|f|_{\Omega,m} ~+~ 2\nu\tfrac{\eta}{\la}. 
}
\et
The proof of the theorem above is similar to the one of \cref{thm:bound-prob-sampled} and it is stated for completeness in \cref{proof:thm:bound-prob-sampled-z}, \cpageref{proof:thm:bound-prob-sampled-z}. The same comments to \cref{thm:bound-prob-sampled} that we reported in the related section and the rates for Sobolev functions, apply also in this case. In the next section we describe the algorithm to solve the problem in \cref{eq:prob-relax-z}.

\subsection{Algorithm}\label{sec:minimizer-algorithm}
We can use the same dual technique as presented in \cref{sec:algorithm}, and obtain a dual problem to \cref{eq:prob-relax-z} with the additional penalty $ \frac{\eps}{n} \log \det B$. The dual problem can readily be obtained as (up to constants)
$$ \inf_{\alpha \in \mathbb{R}^n}
 \sum_{i=1}^n \alpha_i     f(x_i)  - \frac{\eps}{n} \log \det \big(\Phi^\top {\rm Diag}(\alpha) \Phi + \lambda I 
 \big) + \frac{\nu}{2} \Big( -\sum_{i=1}^n \alpha_i \|x_i\|_2^2  + \Big\| \sum_{i=1}^n \alpha_i x_i \Big\|_2^2 \Big), $$
such that $\alpha^\top 1_n = 1$, with the optimal $z$ that can be recovered as $z = \sum_{i=1}^n \alpha_i x_i$. We note that when $\nu$ tends to zero, we recover the dual problem from \cref{sec:algorithm}, and we keep the candidate above in $\Omega$ even when $\nu=0$.

\subsection{Warm restart scheme for linear rates} \label{sec:warm-start}
It is worth noting that \cref{thm:bound-prob-sampled-z} provides strong guarantees on the distance $\|\hat{z} - \zeta\|$ where $\hat{z}$ is the solution of the problem \cref{eq:prob-relax-z} and $\zeta$ the global optimum of $f$. This suggests that we can implement a warm restart scheme that leverage the additional knowledge of the position of $\zeta$. Assume indeed that $\Omega$ is a ball of radius $R$ centered in $z_0$. For $t  = 1,\dots,T$ with $T = \lceil \log \frac{1}{\eps} \rceil$, we apply \cref{eq:prob-relax-z} to a set $\widehat{X}_t$ that contains enough points sampled uniformly at random in the ball $B_{r_{t-1}}(z_{t-1})$  such that \cref{thm:bound-prob-sampled-z} guarantees that $\|z_t - \zeta\| \leq r_{t-1}/e$ where $z_t$ is the solution of \cref{eq:prob-relax-z}. The cycle is repeated with $r_t = r_{t-1}/e$ and the new center be $z_t$. By plugging the estimate of \cref{thm:fill-distance-random-points} for $h_{\widehat{X}_t,B_{r_{t-1}}(z_{t-1})}$ in \cref{thm:bound-prob-sampled-z} for each step $t$, we obtain a total number of points $n$ to achieve $\|z_T - \zeta\| \leq \eps$ with probability $1 - T \delta$, that is
$$n = O\Big(C_{d,m}^{{d}/{m}} \Big(\frac{\cal F}{\nu}\Big)^{{d}/{m}} R^{d} \log \frac{1}{\eps}\Big)$$
modulo logarithmic terms in $n$ and $\delta$, where $C_{d,m} = 3^m C \Cmul \Cdiff$ with $C$ defined in \cref{thm:bound-prob-sampled-z} and ${\cal F} = |f|_{\Omega,m} + \nu + \tr(A_*)$. This means that under the additional assumption of a unique minimizer in $\Omega$, we achieve a convergence rate that is only logarithmic in $\eps$, moreover when $m \gg d$ also the dependence with respect to $C_{d,m}$ (which is exponential in $m$ and $d$ in the case of the Sobolev kernel) and ${\cal F}$ improves, since $d/m$ tends to $0$.

\section{Extensions}
\label{sec:certificate}
In this section we deal with two aspects: (a) the effect of solving approximately the problem in \cref{eq:prob-relax}, and (b) how can we certify explicitly (no dependence on quantities of theoretical interest as $\tr(A_*)$) how close is a given (approximate) solution to the optimum; (c) we will also cover the case when the function $f$ does not have a positive definite representer $A_*$ in $\pdm(\hh)$ but in a larger space. This allows to cover the cases of $f \in C^s(\R^d)$ with $s \leq d/2+2$.

\subsection{Approximate solutions}
In this section we extend \cref{thm:bound-prob-sampled} to consider the case when we solve \cref{eq:prob-relax} in an approximate way. In particular, let $\la > 0, n \in \Np$ and $\widehat{X} = \{x_1,\dots,x_n\}$. Denote by $p_{\la, n}$ the optimal value achieved by \cref{eq:prob-relax} for such $\la, n$. We say that $(\tilde{c}, \tilde{B})$ is an {\em approximate solution} of \cref{eq:prob-relax} with parameters $\theta_1, \theta_2, \tau_1, \tau_2 \geq 0$ if it satisfies the following inequalities
\eqal{
\label{eq:appr-functional} p_{\la,n} ~-~ \tilde{c} + \la \tr(\tilde{B})  &\leq \theta_1 + \theta_2 \tr(\tilde{B}),\\
\label{eq:appr-constraints} |f(x_i) - \tilde{c} - \Phi_i^\top \tilde{B} \Phi_i| & \leq \tau_1 + \tau_2\tr(\tilde{B}), \quad \forall i \in [n].
}
\bt[Error of approximate solutions of \cref{eq:prob-relax}]\label{thm:appr-solution}
Let $(\tilde{c}, \tilde{B})$ be an approximate solution of \cref{eq:prob-relax} for a given $n \in \Np, \la > 0$ as defined in \cref{eq:appr-functional,eq:appr-constraints} w.r.t. $\tau_1, \tau_2, \theta_1, \theta_2 \geq 0$. 
Under the same assumptions and notation of \cref{thm:bound-prob-sampled} and \cref{rem:sufficiency-asm-f}, when $\tau_2, \theta_2 \leq \frac{\la}{8}$
\eqal{
\label{eq:appr-c-approx}
|\tilde{c} - f_*| &~\leq~ 7(2 \tau_1 + \eta |f|_{\Omega,m}) + 6(\theta_1 + \lambda \tr(A_*)) ,\\
\label{eq:appr-A-approx}
\tr(\tilde{B}) &~\leq~ 8\,\tr(A_*) ~+~ 8\tfrac{\eta}{\la}\,|f|_{\Omega,m} ~+~ 8 \tfrac{\theta_1 + 2\tau_1}{\lambda}.
}
\et
The proof of the theorem above is reported for completeness in \cref{proof:thm:appr-solution}, \cpageref{proof:thm:appr-solution}, and is a variation of the one of \cref{thm:bound-prob-sampled} where we used \cref{thm:inequality-scattered-data} with $\tau = \tau_1 + \tau_2\tr(\tilde{B})$ and we further bound $p_{\la,n}$ via \cref{eq:appr-functional}. From a practical side, the theorem above allows to use a wide range of methods and techniques to approximate the solution of \cref{eq:prob-relax}. In particular, it is possible to use lower dimensional approximations of $\Phi_1,\dots, \Phi_n$ and algorithms based on early stopping as described in \cref{sec:discussion}, since $\tau_1,\tau_2,\theta_1,\theta_2$ will take into account the error incurred in the approximations. An interesting application of the theorem above, from a theoretical side is that it allows also to deal with situations where $f$ does not have a representer $A_*$ in $\pdm(\hh)$ as we are going to discuss in the next section.

\subsection{Rates for $f$ with low smoothness}\label{sec:low-smoothness}
When $f \in C^{s+2}(\R^d)$ with $s \in \Nz$, but with a low smoothness, i.e., $s \leq d/2$, we can still apply our method to find the global minimum and obtain almost optimal convergence rates, as soon as it satisfies the geometric conditions in \cref{asm:geom-f:b}, as we are going to show in \cref{thm:global-min-low-smoothness} and the following discussion. 

\noindent In this section, for any function $u$ defined on a super-set of $\Omega$ and $s$ times differentiable on $\Omega$, we define the following semi norm :
\begin{equation}\label{eq:df_norm}
   \|u\|_{\Omega,s} ~~= ~~ \max_{|\alpha| \leq s} ~ \sup_{x \in \Omega} \big|\partial^\alpha u(x)\big|.
\end{equation}

\noindent We consider the following  variation of the problem in \cref{eq:prob-relax}:
\eqal{\label{eq:prob-relax-tau}
\max_{c \in \R, \ B \in \pdm(\R^n)} ~c - \la \tr(B) \quad \mbox{ such that } \quad \forall i \in [n], \ |f(x_i) - c - \Phi_i^\top B \Phi_i| \leq \tau.
}
The idea is that $f$, under the geometric conditions in \cref{asm:geom-f:b}, still admits a decomposition in the form $f(x) = \sum_{j \in [p]} w_j(x)^2$, $p \in \Np$ for any $x \in \Omega$, but now with respect to functions with low smoothness $w_1,\dots, w_p \in C^s(\R^d)$. To prove this we follow the same proof of \cref{sec:inf-dim} noting that the assumptions to apply \cref{lm:repr-conv-g} and \cref{thm:g-good-has-A} are that $f$ belongs to a normed vector space space that satisfy the algebraic properties in \cref{asm:H-rich:a,asm:H-rich:b,asm:H-rich:c} which does not have necessarily to be a RKHS. In particular, note that the space $\widetilde{\hh} = \{f|_{\Omega}~:~f \in C^s(\R^d) \}$ of restriction to $\Omega$ of functions in $\C^s(\R^d)$, endowed with the norm $\|\cdot\|_{\Omega,s} $ defined in \cref{eq:df_norm} (and is always finite on $\widetilde{H}$ since $\Omega$ is bounded) satisfies such assumptions. The reasoning above leads to the following corollary of \cref{thm:g-good-has-A} (the details can be found in \cref{proof:cor:g-in-Cs-has-A}  \cpageref{proof:cor:g-in-Cs-has-A}). 
\bcor\label{cor:g-in-Cs-has-A}
Let $\Omega$ be a bounded open set and $f \in C^{s+2}(\R^d)$, $s \in \Nz$, satisfying \cref{asm:geom-f:b}. Then there exist $w_1,\dots,w_p \in C^s(\R^d)$, $p \in \Np$, such that 
$$ \forall x \in \Omega,~f(x) - f_* = \sum_{j \in [p]} w_j^2(x).$$
\ecor
By using the decomposition above, when the kernel satisfies \cref{asm:H-rich:a}, we build an operator $A_\eps \in \pdm(\hh)$ that approximates $f$ with error $O(\epsilon^s)$ for any $\epsilon > 0$. First note that, for any bounded open set $\Omega \subset \R^d$ and any $s \leq r$, there exists $C_1$ and $C_2$ depending only on $r,s,\Omega$ such that for any $g \in C^s(\R^d)$ and $\eps > 0$ there exists a smooth approximation $g_\epsilon \in C^\infty(\R^d)$ such that $\sup_{x \in \Omega} |g(x) - g_\eps(x)| \leq C_1 \eps^s \|g\|_{\Omega,s}$ and such that
$\|g_\eps\|_{\Omega,r} \leq C_2 \eps^{-(r-s)} \|g\|_{\Omega,s}$ (see Thm. 5.33 of \cite{adams2003sobolev} for the more general case of Sobolev spaces, or \cite[Chapter 21]{cheney2009course} for explicit construction in terms of convolutions with smooth functions). 
Denote by $w^\eps_j$ the smooth approximation of $w_j$ on $\Omega$ for any $j \in [p]$. Since we consider kernels rich enough that the associated RKHS $\hh$ contains smooth functions (see \cref{asm:H-rich:a}), then we have that $w^\eps_j|_\Omega \in \hh$ for any $j \in [p]$. Then 
$$A_\eps = \sum_{j \in [p]} w^\eps_j|_\Omega \otimes w^\eps_j|_\Omega ~~~ \in ~~~\pdm(\hh).$$
The reasoning above is formalized in the next theorem (the proof is in \cref{proof:thm:Aeps-exists}, \cpageref{proof:thm:Aeps-exists}).
\bt\label{thm:Aeps-exists}
Let $d,p,s \in \Nz$. Let $\Omega$ satisfy \cref{asm:geom-f:a} and $f(x) = \sum_{j \in [p]} w_j^2(x),~x \in \Omega$ with $w_j \in C^s(\R^d)$ for $j \in [p]$.  Let $k_r$ be the Sobolev kernel of smoothness $r > \max(s,\frac{d}{2})$ and let $\hh$ be the associated RKHS. Then, for any $\eps \in (0,1]$ there exist $A_\eps \in \pdm(\hh)$ such that 
\eqal{\label{eq:Aeps-approx}
\tr(A_\eps) ~~\leq~~ C\eps^{-2(r-s)}, \qquad \sup_{x \in \Omega} |f(x) - f_* - \scal{\phi(x)}{A_\eps\phi(x)}| ~~\leq~~ C'\eps^{s},
}
where $C= p q w^2, C'= p q' w^2$, and $w = \max_{j \in [p]}\|w_j\|_{\Omega,s}$ and $q, q'$ are constants that depend only on $s,r,d,\Omega$ and are defined in the proof.
\et

Denote now by $(\tilde{c}, \tilde{B})$ one minimizer of \cref{eq:prob-relax-tau}, and consider the problem in \cref{eq:prob-relax} with respect to $f_\eps(x) = \scal{\phi(x)}{A_\eps \phi(x)} + f_*$, i.e.,
\eqal{\label{eq:prob-relax-eps}
\max_{c \in \R, \ B \in \pdm(\R^n)} ~c - \la \tr(B) \quad \mbox{ such that } \quad \forall i \in [n], \ f_\eps(x_i) - c = \Phi_i^\top B \Phi_i,
}
and denote by $p^\eps_{\la,n}$ its optimum.
Since $f_\eps(x_i) - c = \Phi_i^\top B \Phi_i$ implies $|f(x_i) - c - \Phi_i^\top B \Phi_i| \leq \tau$ when $\tau \geq \sup_{x \in \Omega}|f(x) - f_\eps(x)|$, then in this case \cref{eq:prob-relax-tau} is a relaxation of \cref{eq:prob-relax-eps} and we have that $p^\eps_{\la, n} - \tilde{c} - \la \tr(\tilde{B}) \leq 0$. Then, to obtain guarantees on $(\tilde{c}, \tilde{B})$ (the solution of \cref{eq:prob-relax-tau}) we can apply \cref{thm:appr-solution} to the problem in \cref{eq:prob-relax-eps} with $\theta_1,\theta_2, \tau_2 = 0$ and $\tau_1 = \tau$ with the requirement $\tau \geq \sup_{x \in \Omega}|f(x) - f_\eps(x)|$. The reasoning above is formalized in the following theorem and the complete proof is reported in \cref{proof:thm:global-min-low-smoothness}, \cpageref{proof:thm:global-min-low-smoothness}.
\bt[Global minimum for functions with low smoothness]\label{thm:global-min-low-smoothness}
Let $s \in \N$. Let $k_r$ be a Sobolev kernel with smoothness $r \geq s, r > d/2$ and $\hh$ be the associated RKHS. Let $\Omega \subset \R^d$ satisfying \cref{asm:geom-f:a} and $f \in C^{s+2}(\R^d)$, satisfying \cref{asm:geom-f:b}. The problem in \cref{eq:prob-relax-tau} admits a minimizer. Denote by $(\tilde{c}, \tilde{B})$ any of its minimizers for a given $\la > 0, \tau > 0$. With the same notation and the same conditions on $\la$ of \cref{thm:bound-prob-sampled}, when $\tau = \la^{s/(2r-s)}$
\eqals{
|\tilde{c} - f_*| ~\leq~ C_{1,f} (\la + \la^{\frac{s}{2r-s}}), \qquad \tr(\tilde{B}) &~\leq~ C_{2,f} (1+\la^{-(1-\frac{s}{2r-s})}).
}
with $C_{1,f}, C_{2,f}$ defined in the proof and depending only on $f$ and $r,s,d,\Omega$.
\et

The result above allows to derive the following estimate on \cref{alg:glm} applied on the problem in \cref{eq:prob-relax-tau} in the case of a function $f$ with low smoothness. Consider the application \cref{alg:glm} to the problem in \cref{eq:prob-relax-tau} to a function $f \in C^{s+2}(\Omega)$ satisfying \cref{asm:geom-f:b}, with a Sobolev kernel $k_r$, $r \geq s, r > d/2$, and with $\tau = \la^{s/(2r-s)}$, $\la = O(n^{-\frac{r}{d} + 1/2})$ on a set of $n$ points sampled independently and uniformly at random from $\Omega = B_1(0)$, the unit ball of $\R^d$. By combining the result of \cref{thm:global-min-low-smoothness} with the condition on $\la$ in \cref{thm:bound-prob-sampled} and with the upper bound on the fill distance in the case of points sampled uniformly at random in \cref{thm:fill-distance-random-points}, we have that
\eqals{
|\tilde{c} - f_*| ~~=~~ O\Big(\,n^{-\frac{s}{2d}\left(1-\frac{d - s}{2r - s}\right)}\,\Big),
}
modulo logarithmic factors, where $\tilde{c}$ is the solution of \cref{eq:prob-relax-tau}. The rate above must be compared with the optimal rates for global minimization of functions in $C^{s+2}(\Omega)$ via function evaluations, that is $n^{-\frac{s+2}{d}}$ for any $s \in \N$ (Prop. 1.3.9, pag. 34 of \cite{novak2006deterministic}). In the low smoothness setting, i.e., $s \leq d/2$ when we choose $r \gg d/2$, then the term $1-\frac{d - s}{2r - s} \to 1$ and so the exponent of the rate above differs from the optimal one by a multiplicative factor $1/2+\frac{1}{s}$, leading essentially to a rate of $O(n^{-s/(2d)})$. However, the choice of a large $r$ will impact the hidden constants that are not tracked in the analysis above. Then for a fixed $n$ there is a trade-off in $r$ between the constants and the exponent of the rate. So in practice it would be useful to select $r$ by parameter tuning.

\subsection{Certificate of optimality}

While in \cref{thm:bound-prob-sampled} we provide a bound on the convergence of \cref{eq:prob-relax} {\em a priori}, i.e.,  only depending on properties of $f, \Omega, \hh$, in this section we provide a bound {\em a posteriori}, that is a {\em certificate of optimality}. Indeed, the next theorem quantifies $f(z) - f^*$ for a candidate minimizer $z$, in terms of only $(\hat{c}, \hat{B})$, an (approximate) solution of \cref{eq:prob-relax} and $|f|_{\Omega,m}$. A candidate minimizer based on \cref{eq:prob-relax} is provided in \cref{eq:primal-candidate}. In section \cref{sec:minimizer} we study a different algorithm \cref{eq:prob-relax-z} that explicitly provides a minimizer and whose certificate is studied in \cref{sec:certificate-prob-relax-z}.

\bt[Certificate of optimality a minimizer from \cref{eq:prob-relax}]\label{thm:certificate-only-minimum}
Let $\Omega$ satisfy \cref{asm:geom-f:a} for some $r > 0$. Let $k$ be a kernel satisfying \cref{asm:H-rich:a,asm:H-rich:d} for some $m \in \Np$.
Let $\widehat{X} = \{x_1,\dots, x_n\} \subset \Omega$ with $n \in \Np$ such that $h_{\widehat{X},\Omega} \leq \frac{r}{18(m-1)^2}$.
Let $f \in C^m(\Omega)$ and let $\hat{c} \in \R, \hat{B} \in \pdm(\R^n)$ and $\tau \geq 0$ satisfying
\eqal{\label{eq:emp-in-cert-min}
|f(x_i) - \hat{c} ~-~ \Phi_i^\top \hat{B} \Phi_i | \leq \tau, \quad i \in [n],
}
where the $\Phi_i$'s are defined in \cref{sec:sketch}. Let $f_* = \min_{x \in \Omega} f(x)$. Then the following holds
\eqal{
|f(z) -f_*| &\leq f(z) - \hat{c} + \eps + 2\tau, ~~ \forall z \in \Omega, \quad \textrm{where} \quad \eps = C  h_{\widehat{X},\Omega}^m,
}
and $C = C_0(|f|_{\Omega,m}+\Cmul\Cdiff\tr(\hat{B}))$. The constants $C_0$, defined in \cref{thm:inequality-scattered-data}, and $m, \Cmul, \Cdiff$, defined in \cref{asm:H-rich:a,asm:H-rich:d}, do not depend on $n, \widehat{X}$, $h_{\widehat{X},\Omega}, \hat{c}, \hat{B}$ or $f$.
\et
\bpr
By applying \cref{thm:inequality-scattered-data} with $g(x) = f(x) - \hat{c}$, we have $f(x) - \hat{c} \geq - \eps - 2\tau$ for any $x \in \Omega$. In particular this implies that $f(\zeta) - \hat{c} \geq - \eps - \tau$. The proof is concluded by noting that $f(z) \geq f_*$ by definition of $f_*$.
\epr

\section{Relationship with polynomial hierarchies}
\label{sec:polynomial}

The formulation as an infinite-dimensional sum-of-squares bears some strong similarities with polynomial hierarchies. There are several such hierarchies allowing to solve any polynomial optimization problem~\cite{lasserre2001global,lasserre2007sum,lasserre2011new}, but one has a clear relationship to ours. The goal of the following discussion is to shed light on the benefits in terms of condition number and dimensionality of the problem, deriving by using an infinite dimensional feature map in the finite dimensional problem, instead of an explicit finite-dimensional polynomial map as in the case considered by the papers cited above.

\paragraph{Adding small perturbations.} We start this discussion from the following result from Lasserre~\cite{lasserre2007sum}, that is, for any multivariate non-negative polynomial $f$ on $\R^d$, and for any $\eta>0$, there exists a degree $r(f,\eta)$ such that the function 
$$ f_\eta(x) =  f(x) + \eta \sum_{k=0}^{r(f,\eps)} \frac{1}{k!} \sum_{j=1}^d x_j^{2k}
$$
is a sum of squares, and such that the $\ell_1$-norm between the coefficients of $f$ and $f_\eta$ tends to zero (here this $\ell_1$-norm is equal to $ \eta d \sum_{k=0}^{r(f,\eps)} \frac{1}{k!} \leqslant  \eta d e$).

This implies that for the kernel $k_r(x,y) = \sum_{k=0}^r \frac{(x^\top y)^k}{k!}$, with feature map $\phi_r(x)$ composed of all weighted monomials of degree less than $r$, the function
$$
 f(x) + \eta \| \phi_r(x) \|_2^2  = f(x) + \eta k_r(x,x)
$$
is a sum of squares, for any $r \geqslant r(f,\eta)$, with $\eta$ arbitrarily close to zero (this can be obtained by adding the required squares to go from $\sum_{j=1}^d x_j^{2k}$ to $\|x\|^{2k} = (\sum_{j=1}^d x_j^2)^k$).  This result implies that  minimizing $f$ arbitrarily precisely over any compact set $K$ (such that 
$\sup_{x \in K} k_r(x,x)$ is finite), can be done by minimizing $f(x)+\eta k(x,x)$, with sum-of-squares polynomials of sufficiently large degree.
We already showed that in this paper that if $f$ satisfies the geometric condition in \cref{asm:geom-f:b}, our framework is able to find the global minimum by the finite dimensional problem in \cref{eq:prob-relax}, which, in turn, is based on a kernel associated to an infinite dimensional space (as the Sobolev kernel, see \cref{ex:sobolev-kernel}).
We now show how our framework can provide approximation guarantees and potentially efficient algorithms for the problem above even when \cref{asm:geom-f:b} may not hold and we use a polynomial kernel of degree $r$ (with $r$ that may not be large enough). However, in this case the resulting problem would suffer of a possibly infinite condition number and a larger dimensionality than the one achievable with an infinite dimensional feature map.

\paragraph{Modified optimization problem.} Given the representation of $x \mapsto f(x) - f_\ast + \eta \| \phi_r(x) \|_2^2 $ as a sum-of-squares, we can explicitly model the function as
$$
f(x) - c  + \eta \| \phi_r(x) \|_2^2 = \langle \phi_r(x),A \phi_r(x) \rangle
$$
with $A$ positive definite and $\eta \geq 0$. Note that if $r$ is greater than twice the degree of $f$ this problem is always feasible by taking $\eta$ sufficiently large.
Moreover, for feasible $(c,\eta,A)$, we have for any $x \in \Omega$, 
$$f(x) \geq c - \eta \| \phi_r(x)\|^2 \geq c - \eta \sup_{y \in \Omega} \| \phi_r(y)\|_2^2.$$ Thus, a relaxation of the optimization problem is
$$
\sup_{c \in \R, A \succcurlyeq 0 , \eta \geqslant 0} \ \ c - \eta \sup_{y \in \Omega} \| \phi_r(y)\|_2^2 \ \ \mbox{ s. t. } \ \ \forall x\in \Omega, \ f(x) = c +  \phi_r(x)^\top A\phi_r(x) - \eta \| \phi(x)\|_2^2.
$$
Moreover, if we choose $r$ larger than $r(f-f_\ast, \eta)$, we know that there exists a feasible $A$ which is positive semi-definite, with $c = f_\ast - \eta \sup_{y \in \Omega} \| \phi_r(y)\|_2^2$, and thus the objective value is greater than $f_\ast - \eta  \sup_{y \in \Omega} \| \phi_r(y)\|_2^2$. Thus, the objective value of the problem above converges to $f_\ast$, when $\eta$ go to zero (and thus $r(f-f_\ast, \eta)$ goes to infinity), while always providing a lower bound.
Note that  if $f-f_\ast$ is a sum of squares, then the optimal value $\eta$ can be taken to be zero, and we recover the initial problem.

\paragraph{Subsampling and regularization.}
At this point, since $r$ is finite, subsampling ${d \choose 2r}$ points leads to an equivalent finite-dimensional problem. 
We can also add some regularization to sub-sample the problem and avoiding such a large number of points. Note here that the kernel matrix will probably be ill-conditioned, and the problem computationally harder to solve and difficult ro regularize.

\paragraph{Infinite-degree polynomials.}
In the approach outlined above, we need to let $r$ increase to converge to the optimal value. We can directly take $r=\infty$, since  $k_r(x,y) = \sum_{k=0}^r \frac{(x^\top y)^k}{k!}$ tends to the kernel $\exp( x^\top y)$, and here use subsampling. Again, it may lead to numerical difficuties. However, we can use Sobolev kernels (with guarantees on performance and controlled conditioning of kernel matrices), on the function $f(x) + \eta e^{\|x\|_2^2}$ for which we now there exists a sum of squares representation as soon as $f$ is a polynomial.

\section{Experiments}
\label{sec:experiments}

In this section, we illustrate our results with experiments on synthetic data.

\paragraph{Finding hyperparameters.}
Given a function to minimize and a chosen kernel, there are three types of hyperparameters: (a) the number $n$ of sample points, (b) the regularization parameter $\lambda$, and (c) the kernel parameters. Since $n$ drives the running time complexity of the method, we will always set it manually, while we will estimate the other parameters (regularization and kernel), by ``cross-validation'' (i.e.,   selecting the parameters of the algorithm that lead to the minimum value of $f$ at the candidate optimum, among a logarithmic range of parameters). This adds a few function evaluations, but allows to choose good parameters.

\begin{figure}[t]
\begin{center}
\includegraphics[scale=.3]{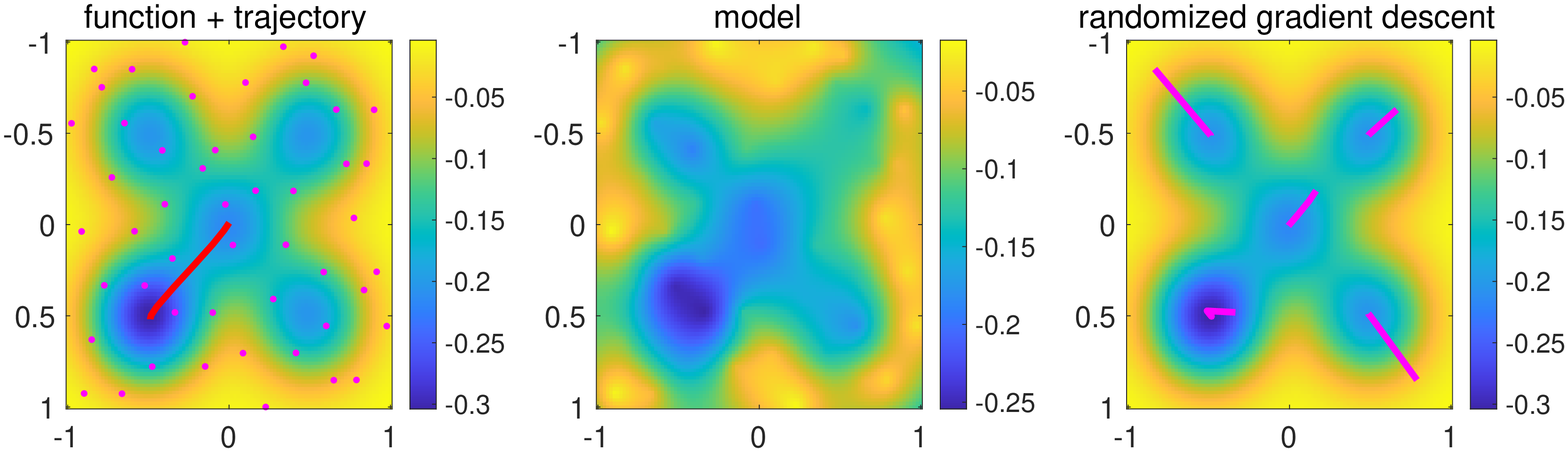} 
\includegraphics[scale=.3]{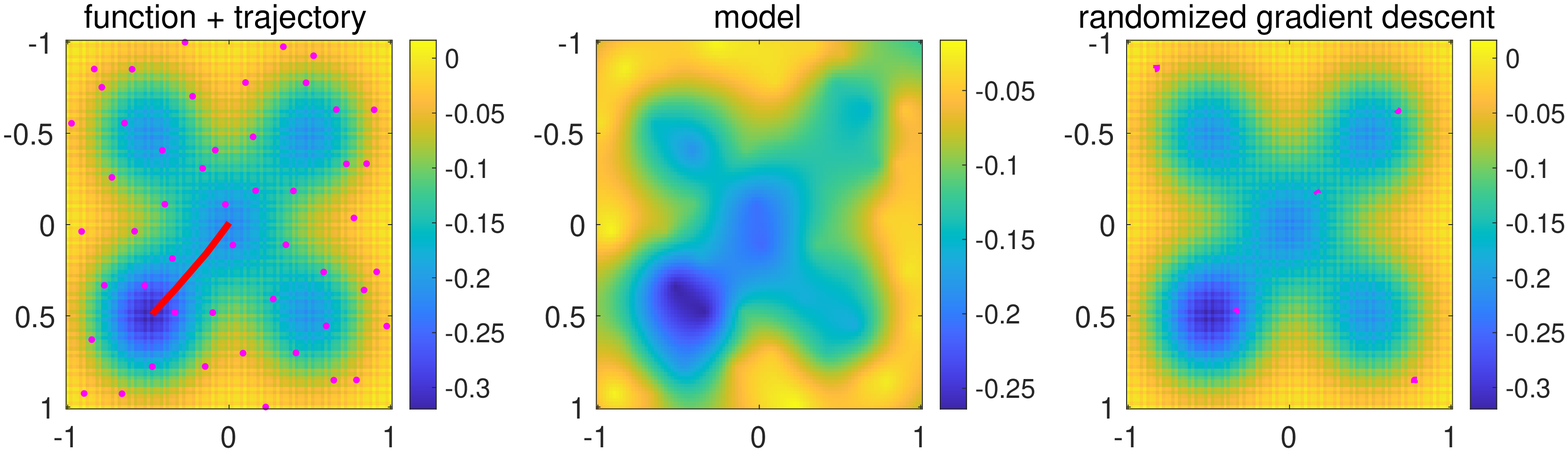} 
\end{center}
\caption{\label{fig:2D}  Top: 2D function without small-amplitude high-frequency components. Bottom: 2D function with small-amplitude high-frequency components. Left: sampled points and the trajectory of the proposed algorithm. Center: model reconstructed by the algorithm (see \cref{eq:appr-model}). Right: the trajectory of gradient descent starting from random points. As it is possible to see, even a small local non-convexity prevents the random+GD algorithms to converge properly, while the proposed method is quite robust to it.
}
\end{figure}

\paragraph{Functions to minimize.} We consider first a simple functions defined in $\R^2$ with their global minimimizer on $[-1,1]^d$, which is minus the sum of Gaussian bumps (see \cref{fig:2D}). To go to higher even dimensions with the possibility of computing the global minimum with high precision by grid search, we consider functions of the form $f(x) = f(x_1,x_2) + f(x_3,x_4) + \cdots + f(x_{d-1},x_d)$. We also consider adding a high-frequency cosine on the coordinate directions representing a more general scenario for a non-convex function. Note that in this second setting the gradient based methods cannot work properly (while ours can) as we are going to see in the simulations.

All results are reported by normalizing function values so that the range of values is 1, that is, $\max_{x \in [-1,1]^d} f(x) = 1$ and $\min_{x \in [-1,1]^d} f(x) = 0$.
 
\paragraph{Baseline algorithms.} We compare our algorithm with the exponential kernel and points sampled from a quasi-random sequence in $[-1,1]^d$, such as the Halton sequence~\cite{niederreiter1992random}, to:
\begin{itemize}
    \item Random search: select a quasi-random sequence in $[-1,1]^d$ and take the point with minimal function value.
    \item Random search with gradient descent: starting gradient descent for a certain number of iterations from quasi-random points, with a number of initialization divided by $d+1$ and the number of gradient steps, to account for gradient evaluations based on $d+1$ function evaluations (by finite-difference). The step-size for gradient descent is taken constant, but its values is optimized for smallest final value while providing a descent algorithm.
\end{itemize}

\paragraph{Illustration in two dimensions.} We show in \cref{fig:2D} a function in two dimensions, with sampled point in purple, the trajectory of the candidate optimum along Newton iterations in red, and the final model of the function.  We also compare to gradient descent with random starting points. We consider two functions below, one without extra high-frequency component (top), and one with (bottom). We can make the following observations:
\begin{itemize}
    \item Our algorithm outperforms random search, that is, it improves on the function values of the sampled points.
    \item For the smoother function, gradient descent performs quite well, but is not robust when high-frequency components are added.
\end{itemize}
Note that the proposed algorithm provides also a model of the function reconstructed starting from its evaluation on the sampled points. In particular, if $(\hat{c}, \hat{B})$ is a solution of the algorithm, the approximate function $\hat{g} \approx f - f^*$ corresponds to 
\eqal{
\label{eq:appr-model}
\hat{g}(x) = \scal{\phi(x)}{V^* B V \phi(x)} = v(x)^\top R^{-1} \hat{B} R^{-\top} v(x), \quad \forall x \in \Omega
}
with $v(x) = (k(x_1,x),\dots, k(x_n,x))$ for $x \in \Omega$ and where $V:\hh \to \R^n$ is in \cref{sec:finite-dim}.

\begin{figure}[t]
\begin{center}
\includegraphics[scale=.45]{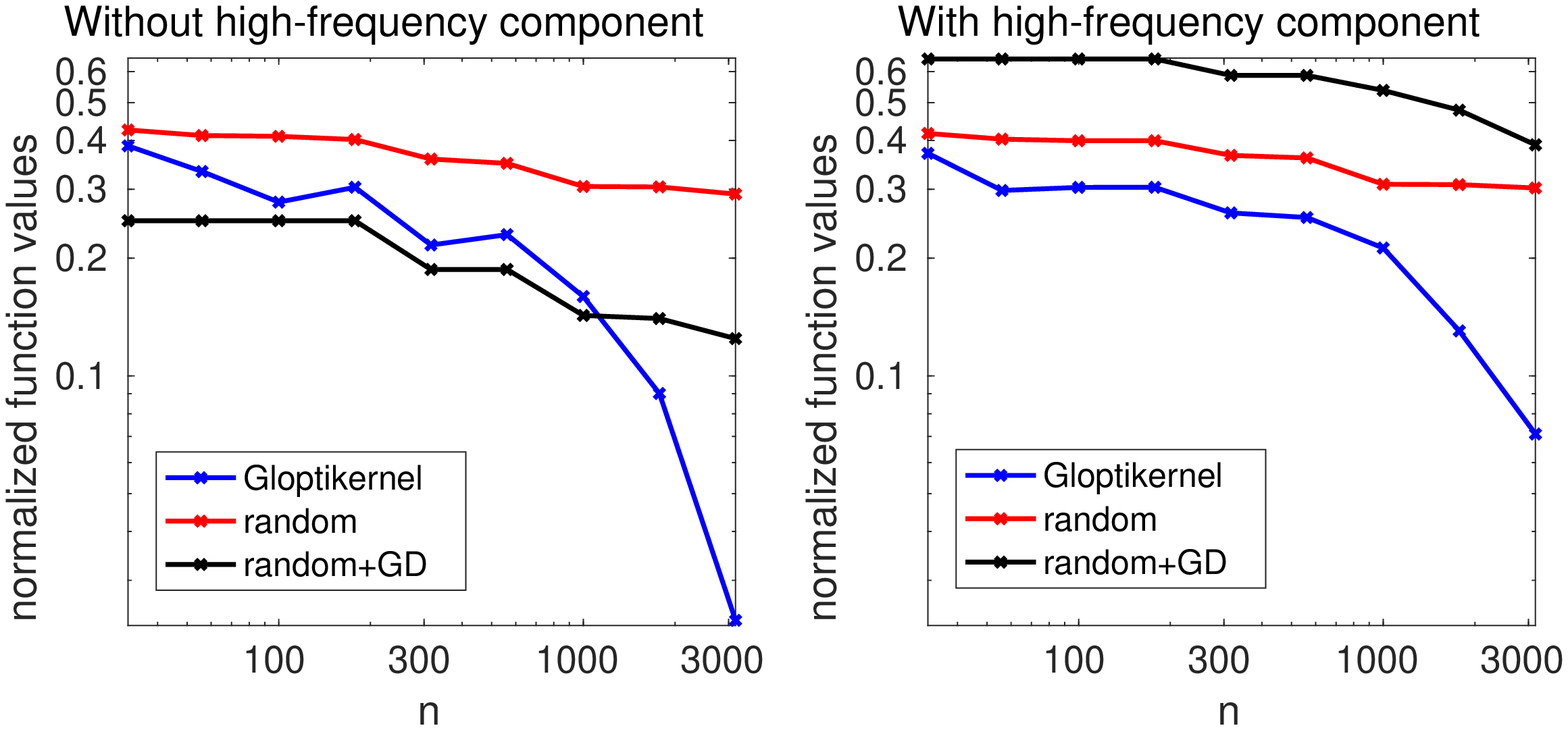} 
\end{center}
\caption{\label{fig:multi-D} Multivariate case $d=8$. Minimization error of our algorithm (gloptikernel) compared with random evaluations or random evaluations + GD. The function considered is built as described at the beginning of this section with domain $[-1,1]^d$ and shifted and rescaled to have minimum in $0$ and output in $[0,1]$. Left: function without small-amplitude high-frequency components. Right: function with small-amplitude high-frequency components.}
\end{figure}

\paragraph{Higher dimensions.} We compare the algorithms on a problem in dimension $d=8$, as $n$ increases, in order to assess how we approach the global optimum. We perform 4 replications with different random seeds for the sampling of points in $[-1,1]^d$. The function to be minimized is built as described at the beginning of this section and is shifted and rescaled to have output in $[0,1]$ and the minimum in $0$. 
We can see that as $n$ gets large, the performance of the proposed algorithm improves, and that with high frequency components, gradient descent with random restarts has worse performance and seem to show a slower rate overall, even in the case of the function without high-frequency components.

\section{Discussion}
\label{sec:discussion}

In this section, we discuss our results and propose a series of extensions.

\paragraph{Main technical contribution and extensions.}
\label{sec:technical-contr}

We see that from \cref{eq:prob-convex-version}, the problem of minimization can be easily written in terms of an infinite set of inequality constraints on $u(x) = f(x) - c$ that must hold for every $x \in \Omega$.
While it is well known how to approximate efficiently an infinite set of equality constraints via a finite subset (e.g. via {\em bounds on functions with scattered zeros} \cite{wendland2004scattered} from the field of approximation theory), leading to optimal rates for the approximation problem, the situation is more difficult in the case of an infinite set of inequality constraints.  The main technical contribution of this paper, on which the whole result of the paper is based, is \cref{thm:inequality-scattered-data}, that allows to deal with an infinite set of inequality constraints as efficiently as in the equality case as discussed in \cref{sec:inequality}. In particular, we rewrite the infinite set of inequalities $g(x) \geq 0, ~ \forall x \in \Omega$ in terms of a very sparse set of constraints of the form $g(x_i) = \Phi_i B \Phi_i$, for some points $x_1,\dots, x_n \in \Omega$ and a matrix $B \in \pdm(\R^n)$, with $n$ in the same order of the one required by the equality case. Assume for simplicity that $\Omega$ is contained in the unit ball and the points are uniformly distributed in $\Omega$. From \cref{thm:inequality-scattered-data} we derive that if $B$ exists,
$$ g(x) ~~\geq -C ~n^{-m/d} ~(|g|_{\Omega,m} + \tr(B)),$$
modulo logarithmic factors, where $m$ is the order of smoothness of $g$.
This result is particularly useful for two reasons. First, it recovers the same dependence on $m$, the smoothness of $g$, and $n$ the number of sample points, as in the case of equality constraints. This is particularly convenient when $m \gg d$, e.g. with $m\geq d$ the rate becomes  $O(n^{-1})$, that is independent from $d$ in the exponent (the dependence of $d$ is still present in the hidden constants and it is exponential in the worst case). Second, if used in an optimization problem, the matrix $B$ can be found via a convex formulation, by requiring $u(x_i) = \Phi_i^\top B \Phi_i$ for $i \in [n]$ and penalizing $\tr(B)$ in the functional. This technique allows, for example, to deal with more general optimization problems with infinite constraints than the one considered in this paper, as
$$\min_{\theta \in \Theta} F(\theta) ~~\mbox{ such that }~~ g(\theta, x) \geq 0, ~ \forall x \in \Omega,$$
by translating it as follows
$$\min_{\theta \in \Theta, B \succeq 0} F(\theta) + \la \tr(B) ~~\mbox{ such that } ~~ g(\theta, x_i) = \Phi_i B \Phi_i ~ \forall i \in [n].$$
If $F$ and $u$ are convex in $\theta$ and $\Theta$ a convex set, then the second is a convex problem that has the potential to approximate very efficiently the first, due to \cref{thm:inequality-scattered-data}. From this viewpoint this paper is an application of this principle to \cref{eq:prob-convex-version}.

\paragraph{Duality.} Beyond using duality in \cref{sec:algorithm} for algorithmic purposes, there is also a dual for the infinite-dimensional problem, which can be written as, 
$$  \inf_{p: \Omega \to \R} \int_\Omega p(x) f(x) dx \mbox{ such that } \int_\Omega p(x) dx = 1 \mbox{ and } \int_\Omega p(x) \phi(x)\! \otimes \!\phi(x) dx \succcurlyeq 0.$$ Replacing the constraint
$\displaystyle \int_\Omega p(x) \phi(x)\! \otimes \! \phi(x) dx \succcurlyeq 0$ by $\forall x \in \Omega, \ p(x) \geqslant 0 $ leads to the usual  relaxation of optimization with probability measures. Thus, our formulation corresponds also to a relaxation in the dual formulation.

\paragraph{Comparison with algorithms based on SOS polynomials.} According to recent results on SOS polynomials (see \cite{slot2020near} and references therein), when $f$ is a polynomial, such algorithms can achieve the global mininum with a rate $O(1/r^2)$ via  an SDP problem based on the representation of SOS polynomials of degree $r$ in terms of positive definite matrices. Since the dimension of the corresponding matrix is $n = {d+r \choose r}$ corresponding to $n = O(r^d)$, by expressing the rate with respect to the dimensionality of the matrix, such methods achieve the global minimum with an error that is in the order of $O(n^{-2/d})$. This can be compared with the approach proposed in this paper as \cref{alg:glm}. By sampling $n$ points from the domain of interest, we cast an SDP problem in terms of a $n$-dimensional positive definite matrix, achieving a rate that is $C_{s,d} n^{-s/d + 1/2}$ (see \cref{thm:alg-glm}) modulo logaritmic factors, by using a Sobolev kernel $k_{s+3}$ with $s > d/2$ (see \cref{ex:sobolev-kernel}). Since the polynomials are arbitrarily differentiable, we can choose $s$ arbitrarily large at the cost of a larger constant $C_{s,d}$ completely characterized in \cref{thm:alg-glm}. For example, by choosing $s = 5d/2$ we achieve the global minimum with a rate $O(n^{-2})$ that does not suffer of the curse of dimensionality except in the constants, and that is  faster than the one obtained by SOS polynomial methods especially when $d \gg 1$. It must be noted that our result holds under the sufficient assumption \cref{asm:geom-f:b} that can be relaxed according to \cref{rem:sufficiency-asm-f}, but that it is not required by SOS polynomial methods. It would be of interest to know if such methods can achieve our rates under the same assumption.

\paragraph{Comparison with simpler algorithms.} Similar reasoning can be done with respect to simple algorithms for global optimization. We consider here the algorithm that consists in sampling~$n$ points at random in $\Omega$ and taking the one with minimum value. A simple analysis based on Lipschitzianity of $f$ shows that this method achieves a rate of $O(n^{-1/d})$. So our method is stricly better than taking the minimum $f(x_i)$ for $i \in [n]$ when $f$ is at least $3$-times differentiable (see \cref{sec:low-smoothness}).

\paragraph{Obtaining optimal rates.} Our current analysis, even for functions $f$ in Sobolev spaces, does not lead to the optimal rate of convergence (we obtain an extra term of $2/d$ in the exponents). We conjecture, that this could be removed by a more refined analysis (in particular in the construction of the operator $A_\ast$).

\paragraph{Modelling gradients.} Our current framework only used function values. If gradients are observed, it could be possible to use them to reduce the number of sampled points, using tools from~\cite{zhou2008derivative}.

\paragraph{Efficient kernel approximations.} The current algorithm has a complexity of $O(n^3)$ for $n$ sampled points, partly due to the need to compute inverse of kernel matrices. There is a large literature within machine learning aiming at providing low-rank approximations, either from approximations of $K$ from a subset of its columns (see, e.g.,~\cite{bach2013sharp,rudi2015less} and references therein) or using random feature vectors (see, e.g.,~\cite{rudi2017generalization,bach2017equivalence}  and references therein). This requires to relax the equality constraint on the subset $\widehat{X}$ to an mean square deviations, as allowed by \cref{sec:certificate}.

\paragraph{Constrained optimization.} Following~\cite{lasserre2001global}, we can apply the same algorithmic technique to constrained optimization, by formulating the problem of minimizing $f(x)$ such that $g(x)\geqslant 0$ as maximizing $c$ such that $f(x) = c + p(x) + g(x) q(x)$, and $p, q$ non-negative functions. We can then replace the non-negative constraints by $p(x) = \langle \phi(x), A \phi(x)\rangle$ and $q(x) = \langle \phi(x), B \phi(x)\rangle$ for positive operators $A$ and $B$. We can then subsample and penalize the traces of $A$ and $B$ to obtain an algorithm. A detailed study of the approximation properties of this algorithm remains to be done.

\paragraph{Acknowledgements.}
This work was funded in part by the French government under management
of Agence Nationale de la Recherche as part of the “Investissements d’avenir” program, reference
ANR-19-P3IA-0001 (PRAIRIE 3IA Institute). We also acknowledge support from the European
Research Council (grant SEQUOIA 724063).

%
%

\bibliographystyle{unsrt}      
\bibliography{biblio}   

\appendix

\section{Additional notation and definitions}

We provide here some basic notation that will be used in the rest of the appendices.

\paragraph{Multi-index notation.} Let $\alpha \in \Nz^d$, $x \in \R^d$ and $f$ be an infinitely differentiable function on $\R^d$, we introduce the following notation
$$|\alpha| = \sum_{j \in [d]} \alpha_i, \quad \alpha! = \prod_{j \in [d]} \alpha_j!, \quad x^\alpha = \prod_{j \in [d]} x_j^{\alpha_j}, \quad \partial^\alpha f = \frac{\partial^{|\alpha|} f}{\partial x_1^{\alpha_1}\cdots\partial x_d^{\alpha_d}}.$$

\paragraph{Some useful space of functions.}
Let $\Omega$ be an open set. In this paper we will denote by $C^s(\Omega)$, $s \in \Nz$, the set of $s$-times differentiable functions on $\Omega$ and by $C^s_0(\Omega)$ the set of functions that are differentiable at least $s$ times and that are supported on a compact in $\Omega$. 
Denote by $L^p(\Omega)$ the {\em Lebesgue space} of $p$-integrable functions with respect to the Lebesgue measure and denote by $\|\cdot\|_{L^p(\Omega)}$ the associated norm \cite{adams2003sobolev}. 

\subsection{Fourier Transform.} 
Given two functions $f,g:\Omega \to \R$ on some set $\Omega$, we denote by $f \cdot g$ the function corresponding to {\em pointwise product} of $f, g$, i.e.,
$$(f \cdot g)(x) = f(x)g(x), \quad \forall x \in \Omega.$$
Let $f, g \in L^1(\R^d)$ we denote the {\em convolution} by $f \star g$ 
$$(f \star g)(x) = \int_{\R^d} f(y) g(x-y) dy.$$
Let $f \in L^1(\R^d)$. The Fourier transform of $f$ is denoted by $\tilde{f}$ and is defined as 
$$\tilde{f}(\omega) = (2\pi)^{-\frac{d}{2}}\int_{\R^d} e^{-i \,\omega^\top x} \,f(x)\, dx,$$
We now recall some basic properties, that will be used in the rest of the appendix.
\bp[Basic properties of the Fourier transform \cite{wendland2004scattered}, Chapter 5.2.]\label{prop:fourier}
$ $

\begin{enumprop}
\item\label{prop:fourier:scale} Let $f \in L^1(\R^d)$ and let $r > 0$. Denote by $\tilde{f}$ its Fourier transform and by $f_r$ the function $f_r(x) = f(x/r)$ for all $x \in \R^d$, then
$$\tilde{f}_r(\omega) = r^d \tilde{f}(r\omega).$$ 
\item\label{prop:fourier:product} Let $f, g \in L^1(\R^d)$, then 
$$\widetilde{f \cdot g} = (2\pi)^{d/2}  \tilde{f} \star \tilde{g}.$$
\item\label{prop:fourier:derivative} Let $\alpha \in \N_0^d$, $f: \R^d \to \R$ and $f, \partial^\alpha f \in L^1(\R^d)$, then
$$\widetilde{\partial^\alpha f}\,(\omega) = i^{|\alpha|} \omega^\alpha \tilde{f}(\omega), \quad \forall \omega \in \R^d.$$
\item\label{prop:fourier:Linfty-omega} Let $f \in L^1(\R^d)$, then 
$$\|\tilde{f}\|_{L^\infty(\R^d)} \leq (2\pi)^{-d/2} \|f\|_{L^1(\R^d)}.$$
\item\label{prop:fourier:Linfty-x} Let $f \in L^1(\R^d)$ and assume that $\tilde{f} \in L^1(\R^d)$, then
$$f(x) = (2\pi)^{-\frac{d}{2}}\int_{\R^d} e^{i \,\omega^\top x} \,\tilde{f}(\omega)\, dx, \quad \textrm{and} \quad \|f\|_{L^\infty(\R^d)} \leq (2\pi)^{-d/2} \|\tilde{f}\|_{L^1(\R^d)}.$$
\item\label{prop:fourier:L2} There exists a linear isometry ${\cal F}: L^2(\R^d) \to L^2(\R^d)$ satisfying 
$${\cal F} f = \tilde{f}, \quad f \in L^2(\R^d) \cap L^1(\R^d).$$
The isometry is uniquely determined by the property in the equation above. For any $f \in L^2(\R^d)$ we denote by $\tilde{f}$ the function $\tilde{f} = {\cal F} f$.
\end{enumprop}
\ep

\subsection{Sobolev Spaces}\label{sec:sobolev-spaces} 
For this section we refer to \cite{adams2003sobolev}.
For any $\alpha \in \N_0^d$ we say that $v_\alpha \in L^1_{loc}(\R^d)$ is the $\alpha$-{\em weak derivative} of $u \in L^1_{loc}(\R^d)$ if, for all compactly supported smooth functions $\tau \in C^\infty_0(\R^d)$, we have 
$$\int_{\R^d} v_\alpha(x) \tau(x) dx = (-1)^{|\alpha|}\int_{\R^d} u(x) (\partial^\alpha\tau)(x) dx,$$
and we denote $v_\alpha$ by $D^\alpha u$. Let $\Omega \subseteq \R^d$ be an open set. For $s \in \N, p \in [1,\infty]$ the Sobolev spaces $W^s_p(\Omega)$ are defined as
$$W^s_p(\Omega) = \{f \in L^p(\Omega) ~|~ \|f\|_{W^s_p(\Omega)} < \infty\}, \quad \|f\|_{W^s_p(\Omega)} = \sum_{|\alpha| \leq s} \|D^\alpha f\|_{L^p(\Omega)}.$$
We now recall some basic results about Sobolev spaces that are useful for the proofs in this paper.
First we start by recalling the restriction properties of Sobolev spaces. Let $\Omega \subseteq \Omega' \subseteq \R^d$ be two open sets. Let $\beta \in \N$ and $p \in [1,\infty]$. By definition of the Sobolev norm above we have
$$\|g|_\Omega\|_{W^s_p(\Omega)} \leq \|g\|_{W^s_p(\Omega')},$$
and so $g|_\Omega \in W^s_p(\Omega)$ for any $g \in W^s_p(\Omega')$. Now we recall the extension properties of Sobolev spaces.
\bp[Extension operator, 5.24 in \cite{adams2003sobolev}]\label{prop:extension-sobolev}
Let $\Omega$ be a bounded open subset of $\R^d$ with locally Lipschitz boundary \cite{adams2003sobolev}. Let $\beta \in \N$ and $p \in [1,\infty]$.
There exists a bounded operator $E:W^\beta_p(\Omega) \to W^\beta_p(\R^d)$ and a constants $C_3$ depending only on $\beta, p, \Omega$ such that for any $h \in W^\beta_p(\Omega)$ the following holds 
(a) $h = (Eh)|_\Omega$  (b) $\|Eh\|_{W^\beta_p(\R^d)} \leq C_3 \|h\|_{W^\beta_p(\Omega)}$ with $C_3 = \|E\|_{op}$. 
\ep

\bp[Approximation property of Sobolev spaces, Thm~5.33 in \cite{adams2003sobolev}]\label{prop:approximation-sobolev}
Let $\Omega$ be a bounded open subset of $\R^d$ with locally Lipschitz boundary \cite{adams2003sobolev},  or $\Omega = \R^d$. Let $s,d \in \N, r \geq s$ and $p \in [1, \infty]$. There exists $C_1$ depending only on $s,d,p$ and $C_2$ depending only on $r,s,d,p$ such that for any $\eps \in (0,1]$ and $g \in W^s_p(\Omega)$ there exists $g_\eps \in C^\infty(\Omega)$ satisfying (i) $g_\epsilon$ is the restriction to $\Omega$ of a certain $\tilde{g}_\epsilon 
\in C^\infty(\R^d)$ and (ii) 
$$\|g - g_\eps\|_{L^p(\Omega)} \leq C_1\eps^{s} \|g\|_{W^s_p(\Omega)}, \qquad \|g_\eps\|_{W^r_p(\Omega)} \leq C_2 \eps^{-(r-s)} \|g\|_{W^s_p(\Omega)}.$$
\ep
\bpr
The case $\Omega = \R^d$ is covered explicitly by Thm.~5.33 in \cite{adams2003sobolev}. The result holds also for $W^s_p(\Omega)$ when $\Omega$ has Lipschitz boundaries as discussed in \cite{adams2003sobolev}, above Theorem 5.33. The result is obtained considering that when $\Omega$ has Lipschitz boundaries, then there exists a bounded extension operator between $W^s_p(\Omega)$ and $W^s_p(\R^d)$ \cite{adams2003sobolev}. Here we provide the proof for the sake of completeness. Let $g \in W^s_p(\Omega)$ and let $\eps \in (0,1]$. Then, by \cref{prop:extension-sobolev} since $\Omega$ has Lipschitz boundary, there exists a bounded extension operator $E:W^s_p(\Omega) \to W^s_p(\R^d)$. Denote by $\tilde{g}$ the function $\tilde{g} = Eg$ and note that $\tilde{g} \in W^s_p(\R^d)$. Then, by applying Thm.~5.33 in \cite{adams2003sobolev} we have that there exists $\tilde{g}_\eps \in C^\infty(\R^d)$ such that
$$\|\tilde{g} - \tilde{g}_\eps\|_{L^p(\R^d)} \leq C\eps^{s} \|\tilde{g}\|_{W^s_p(\R^d)}, \qquad \|\tilde{g}_\eps\|_{W^r_p(\R^d)} \leq C' \eps^{-(r-s)} \|\tilde{g}\|_{W^s_p(\R^d)},$$
for some $C$ depending only on $s,p$ and $C'$ depending on $r,s,p$.
Since by \cref{prop:extension-sobolev} we have $\|\tilde{g}\|_{W^s_p(\R^d)} = \|Eg\|_{W^s_p(\R^d)} \leq C_3 \|g\|_{W^s_p(\Omega)}$, so
$$\|g - \tilde{g}_\eps|_\Omega\|_{L^p(\Omega)} \leq \|\tilde{g} - \tilde{g}_\eps\|_{L^p(\R^d)} \leq C \eps^s \|\tilde{g}\|_{W^s_p(\R^d)}   \leq C C_3 \eps^s \|g\|_{W^s_p(\Omega)},$$
and analogously, 
$$\|\tilde{g}_\eps|_\Omega\|_{W^r_p(\Omega)} \leq \|\tilde{g}_\eps\|_{W^r_p(\R^d)} \leq C' \eps^{s-r} \|\tilde{g}\|_{W^s_p(\R^d)} \leq C' C_3 \eps^{s-r}\|g\|_{W^s_p(\Omega)}.
$$
The proof is concluded by taking $g_\eps = \tilde{g}_\eps |_\Omega$ and $C_1 = C C_3, C_2 = C' C_4$.
\epr

In the next proposition we recall some aspects of the more general {\em Sobolev embedding theorem} \cite{adams2003sobolev}.

\bp\label{prop:embedding-sobolev}
Let $\Omega$ be a bounded open set with Lipschitz continouos boundary. Let $r \in \N$ and $1 \leq p \leq q \leq \infty$. Then $W^r_q(\Omega) \subseteq W^r_p(\Omega)$. In particular there exists a constant $C_5$ such that 
$$\|\cdot\|_{W^r_p(\Omega)} \leq C_5 \|\cdot\|_{W^r_q(\Omega)}.$$

Finally, note that for any $f \in C^r(\R^d)$, it holds  $f|_\Omega \in  W^r_\infty(\Omega)$.
\ep
\bpr
The main statement of the proposition is a subcase of the more general Sobolev embedding theorem \cite{adams2003sobolev}.

\noindent Finally, we recall that, since $f \in C^r(\R^d)$ and $\Omega$ is bounded, then $\partial^\alpha f$ is uniformly bounded on $\Omega$, for any $\alpha \in \N^d$ satisfying $|\alpha| \leq r$. This implies that $f|_\Omega \in W^r_\infty(\Omega)$.
\epr

Finally, note that the semi-norm $\|\cdot\|_{\Omega,r}$ defined in \cref{eq:df_norm} and the Sobolev norm $\|\cdot\|_{W^r_\infty}$ are equivalent in the following sense.

\bp\label{prop:equivnorm}
Let $\Omega \subset \Omega^{\prime}$ be two bounded open sets. Let $r \in \Nz$. For any $u \in C^r(\Omega^\prime)$, recall the definition of $\|u\|_{\Omega,r}$ from \cref{eq:df_norm}. There exists an explicit constant $C_6 > 0$ such that
$$ \forall u \in C^r(\Omega^\prime),\tfrac{1}{C_6}\|u|_\Omega\|_{W^r_\infty(\Omega)} \leq  \|u\|_{\Omega,r} \leq C_6 \|u|_\Omega\|_{W^r_\infty(\Omega)}.$$
Note that this inequality holds also when the norms are unbounded, by using the convention $+\infty \leq + \infty$.
\ep
\bpr

  Since by \cref{eq:df_norm},  $\|u\|_{\Omega,r} = \max_{|\alpha| \leq r}\|\partial^\alpha u\|_{L^\infty(\Omega)}$ and $\|u|_\Omega\|_{W^r_\infty(\Omega)} = \sum_{|\alpha| \leq r}\|\partial^\alpha u\|_{L^\infty(\Omega)}$, and $\{|\alpha|\leq r\}$ is of size $1 + d + ... + d^r = \frac{d^{r+1} - 1}{d-1}$ (where this is taken to be equal to $k+1$ in the case where $d = 0$), the result holds for $C_6 = \frac{d^{r+1} - 1}{d-1}$.

\epr

\subsection{Reproducing Kernel Hilbert spaces}\label{sec:rkhs-appedix}

For this section we refer to \cite{aronszajn1950theory,steinwart2008support,paulsen2016introduction}. Let $S$ be a set and $k:S \times S \to \R$ be a p.d. kernel. We denote by $\hh_k(S)$ the reproducing kernel Hilbert space (RKHS) associated to the kernel $k$, and by $\scal{\cdot}{\cdot}_k$ the associated inner product. In particular, we will omit the dependence in $k$ from $\hh$ and $\scal{\cdot}{\cdot}$ when the used kernel is clear from the context. We will omit also the dependence on $S$ when $S = \Omega$, the region we are using in this paper. In particular we will use the following shortcuts $\hh = \hh_k(\Omega)$ and $\hh(\R^d) = \hh_k(\R^d)$.

\paragraph{Concrete constructions and useful characterizations.} 
In the rest of the section we provide other methods to build RKHS and some interesting characterizations of $\hh_k(S)$ and $\scal{\cdot}{\cdot}_k$ that will be useful int the rest of the appendix.

\bp[Construction of RKHS given $S, \phi$, Thm. 4.21 of \cite{steinwart2008support}]
Let $\phi: S \to V$ be a continuous map, where $V$ is separable Hilbert space with inner product $\scal{\cdot}{\cdot}_V$. Let $k(x,x') = \scal{\phi(x)}{\phi(x')}_V$ for any $x,x' \in S$. Then $k$ is a p.d. kernel and the associated RKHS is characterized as follows
$$\hh_k(S) = \{\scal{w}{\phi(\cdot)}_V~|~ w \in V\}, \quad \|f\|_{\hh_k(S)} = \inf_{u \in V} \|u\|_V ~~s.t.~~ f = \scal{u}{\phi(\cdot)}_V.$$  
\ep

\bp[Restriction of a RKHS $\hh_{k_1}(S_1)$ on a subset $S_0 \subset S_1$ \cite{aronszajn1950theory,paulsen2016introduction}]\label{ex:restriction-rkhs}
Let $k_0$ be the restriction on $S_0$ of the kernel $k_1$ defined on $S_1$. Then the following holds
\begin{enumprop}
\item $k_0$ is a p.d. kernel,
\item the RKHS $\hh_{k_0}(S_0)$ is characterized as $\hh_{k_0}(S_0) = \{ f|_{S_0} ~|~ f \in \hh_{k_1}(S_1)\}$,
\item the norm $\|\cdot\|_{\hh_{k_0}(S_0)}$ is characterized by
$$ \|f\|_{\hh_{k_0}(S_0)} = \inf_{g \in \hh_{k_1}(S_1)} \|g\|_{\hh_{k_1}(S_1)}, ~~~s.t.~~~ f(x) = g(x) ~\forall x \in S_0,$$
\item  there exist a linear bounded {\em extension} operator $E:\hh_{k_0}(S_0) \to \hh_{k_1}(S_1)$ such that $(E f)(x) = f(x)$ for any $x \in S_0$ and $f \in \hh_{k_0}(S_0)$ and such that
$$\|f\|_{\hh_{k_0}(S_0)} = \|Ef\|_{\hh_{k_1}(S_1)}, \quad \forall f \in \hh_{k_0}(S_0),$$
\item there exist a linear bounded {\em restriction} operator $R:\hh_{k_1}(S_1) \to \hh_{k_0}(S_0)$ such that $(R f)(x) = f(x)$ for any $x \in S_0$ and $f \in \hh_{k_1}(S_1)$,
\item $R$ and $E$ are partial isometries. In particular $E = R^*$ and $RE$ is the identity on $\hh_{k_0}(S_0)$, while $ER$ is a projection operator on $\hh_{k_1}(S_1)$.
\end{enumprop}
\ep

\bp[Translation invariant kernels on $\R^d$]\label{ex:tr-inv-kernel}
Let $v:\R^d \to \R$ such that its Fourier transform $\tilde{v}$ is integrable and satisfies $\tilde{v} \geq 0$ on $\R^d$. Then 
\begin{enumprop}
\item The function $k:\R^d \times \R^d \to \R$ defined as $k(x,x') = v(x-x')$ for any $x,x' \in \R^d$ is a kernel and is called {\em translation invariant kernel}.
\item The RKHS $\hh_k(\R^d)$ and the norm $\|\cdot\|_{\hh_k(\R^d)}$ are characterized by $$\hh_k(\R^d) = \{f \in L^2(\R^d) ~|~ \|f\|_{\hh_k(\R^d)} < \infty \}, \quad \|f\|^2_{\hh_k(\R^d)} = (2\pi)^{-\tfrac{d}{2}}\int_{\R^d} \frac{|({\cal F}f)(\omega)|^2}{\tilde{v}(\omega)}d\omega,$$
where ${\cal F} f$ is the Fourier transform of $f$ (see \cref{prop:fourier} for more details on ${\cal F}$).
\item The inner product $\scal{\cdot}{\cdot}_k$ is characterized by
$$\scal{f}{g}_k = (2\pi)^{-\tfrac{d}{2}} \int_{\R^d} \frac{({\cal F}f)(\omega)\overline{({\cal F}g)(\omega)}}{\tilde{v}(\omega)} d\omega.$$
\end{enumprop}
\ep

\subsection{Auxiliary results on $C^\infty$ functions}\label{sec:c-infty-results}

\bp\label{prop:Cinfty-extension}
Let $U$ be an open set of $\R^d$ and $K \subset U$ be a compact set. Let $u \in C^\infty(U)$, then there exists $v \in C^\infty_0(\R^d)$ (with compact support), such that $v(x) = u(x)$ for all $x \in K$.
\ep
\bpr
By Thm.~1.4.1, pag.~25 of \cite{hormander2015analysis} there exists $z_{K,U} \in C_0^\infty(U)$, i.e.,  a smooth function with compact support, such that $z_{K,U}(x) \in [0,1]$ for any $x \in U$ and $z(x) = 1$ for any $x \in K$.
Consider now the function $v_{K,U}$ defined as $v_{K,U}(x) = z_{K,U}(x)u(x)$ for all $x \in U$. The function $v_{K,U}$ is in $C^\infty_0(U)$, since it is the product of a $C^\infty_0(U)$ and a $C^\infty(U)$ function, moreover $v_{K,U}(x) = u(x)$ for all $x \in K$. The theorem is concluded by defining $v$ as the extension of $v_{K,U}$ to $\R^d$, i.e.,  the function $v_{K}(x) = z_{K,U}(x)$ for any $x \in U$ and $v_{K}(x) = 0$ for any $x \in \R^d\setminus U$. This is always possible since $v_{K,U}$ is supported on a compact set $K'$ which is contained in the open set $U$, so $v_{K,U}$ is already identically zero in the open set $U \setminus K'$.
\epr

\blm\label{lm:bump}
Given $\zeta \in \R^d$ and $r > 0$, there exists $u \in C^\infty_0(\R^d)$ such that for any $x \in \R^d$, it holds
\begin{enumerate}[(i)]
    \item $u(x) \in [0,1]$;
    \item $\|x\| \geq r \implies u(x) = 0$;
    \item $\|x\| \leq r/2 \implies u(x) = 1$.
\end{enumerate}
\elm

\bpr 
Assume without loss of generality that $\zeta = 0$ and $r = 1$. Consider the following functions : 
\[u_1(x) = \begin{cases} \exp\left(-\frac{1}{1 - \|x\|^2}\right)&\text{if} \,\|x\| <1 \\ 
0 &\text{otherwise}\end{cases},\qquad    u_2(x) = \begin{cases} \exp\left(-\frac{1}{ \|x\|^2 -1/4}\right)&\text{if} \, \|x\| > 1/2 \\ 
0 &\text{otherwise}\end{cases}.\]
Both $u_1$ and $u_2$ belong to $C^{\infty}(\R^d)$ with values in $[0,1]$. Moreover, $u_1 > \alpha_1$ on $B_{3/4}(0)$ and $u_2 \geq \alpha_2$ for some $\alpha_1, \alpha_2 > 0$ on $\R^d \setminus{B_{3/4}(0)}$, which implies that $u_1 + u_2 \in I$ on $\R^d$, where $I = [\min(\alpha_1,\alpha_2), 2]$. Since $(\cdot)^{-1}$ is infinitely differentiable on $(0,\infty)$ we see that $1/(u_1+u_2)$ is well defined on all $\R^d$ and belongs to $C^\infty(\R^d)$, since $I \subset\subset (0,\infty)$.
Consider the function 
$$u_0 = \frac{u_1}{u_1 + u_2}.$$
It is non-negative, bounded by $1$, and infinitely differentiable as a product. Moreover : 
\[\forall x \in B_{1/2}(0),~ u_2(x) = 0 \implies u_0(x) = 1,\qquad \forall x \in \R^d,~ u_1(x) = 0 \Leftrightarrow u_0(x) =0 \Leftrightarrow x \in \R^d\setminus{B_1(0)}. \]
 To conclude the proof, given $r > 0$ and $\zeta \in \R^d$ we will take $u(x) = u_0((x-\zeta)/r)$.
\epr

\blm\label{lm:exists-b2}
Let $N \in \Np$, $\zeta_1,...,\zeta_N \in \R^d$ and $r_1,...,r_N >0$. For $n \in \{1,\dots,N\}$, let $B_n = B_{r_n}(\zeta_n)$ be the open ball centered in $\zeta_n$ of radius $r_n$ and $B^{\prime}_n = B_{r_n /2}(\zeta_n) \subset B_n$ be the open ball centered in $\zeta_n$ of radius $r_n/2$. Then there exists functions $v_0,v_1,...,v_N \in C^\infty(\R^d)$ such that 
\begin{enumerate}[(i)]
    \item $v_0 = v_0\cdot \boldsymbol{1}_{\R^d \setminus{\bigcup_{n =1}^N B^{\prime}_n}}$
    \item $v_n = v_n\cdot \boldsymbol{1}_{B_n},~ \forall n \in \{1,\dots,N\}$
    \item $\sum_{n =0}^N{v_n^2} = 1$.
\end{enumerate}
\elm 

\bpr 
For all $n \in [N]$, take $u_n$ as in \cref{lm:bump} with $r = r_n, \zeta = \zeta_n$ and define $u_0 = \prod_{n=1}^N{(1-u_n)}$. Since $\forall n \in [N],~ u_n \in [0,1]$, we also have $u_0 \in [0,1]$. 
Moreover, let $R = \max_{n \in [N]}{\|\zeta_n\| + r_n}$, then 
\[\forall \|x\| \geq R,~ \forall 1 \leq n \leq N,~u_n(x) = 0 \text{ and } u_0(x) = 1.\]

\noindent{\bf Step 1.} \textit{$u_0\cdot \boldsymbol{1}_{\R^d \setminus{\bigcup_{n \in [N]}B^{\prime}_n}} = u_0$ and for all $n \in [N]$, $u_n \cdot \boldsymbol{1}_{B_n} = u_n$.}

 By point (iii)  of \cref{lm:bump}, $u_n = 1$ on $B^{\prime}_n$ for all $n \in [N]$, which shows that $u_0 = 0$ on $\bigcup_{n = 1}^N{B^{\prime}_n}$ and hence $u_0\cdot \boldsymbol{1}_{\R^d \setminus{\bigcup_{n \in [N]}B^{\prime}_n}} = u_0$. On the other hand, for all $n \in [N]$, point (ii) of \cref{lm:bump} directly implies $u_n \cdot \boldsymbol{1}_{B_n} = u_n$.\\\

\noindent{\bf Step 2.} \textit{The function $\frac{1}{\sqrt{\sum_{n = 0}^N{u_i^2}}}$ is well defined and in $C^\infty(\R^d)$.}

By definition of $u_0$, if $u_0(x) = 0$, then there exists $n \in [N]$ such that $u_n(x) = 1$. Since all the $u_n$ are non-negative, this shows that $s := \sum_{n =0}^N{u_n^2} > 0$. Moreover, consider the closed ball $\bar{B}$ of radius $R$ and centered in $0$. Since $\bar{B}$ is compact, $s$ is continuous and $s(x) > 0$ for any $x \in \bar{B}$, then there exists $0 < m_R \leq M_R < \infty$ such that $s(x) \in [m_R, M_R]$ for any $x \in \bar{B}$. 
Moreover, since for any $\|x\| \geq R,~ u_0(x) = 1 \text{ and }\forall n \in [N],~u_n(x) = 0$, we see that 
\[\forall x \in \R^d \setminus{B_R(0)},~ \sum_{n =0}^N{u_n^2(x)} = 1.\]
Then $s \in [m, M]$ for any $x \in \R^d$, where $m = \min(m_R,1)$ and $M = \max(M_R,1)$.

Since the interval $I = [m,M]$ is a compact set included in the open set $(0,\infty)$ and ${1}/{\sqrt{\cdot}}$ is infinitely differentiable on $(0,\infty)$ then by \cref{prop:Cinfty-extension} there exists $q_I \in C^\infty_0(\R)$ such that $q_I(x) = 1/\sqrt{x}$ for any $x \in I$. Since $s(x) \in I$ for any $x \in \R^d$ we have 
$$\frac{1}{\sqrt{\sum_{n = 0}^N{u_i^2}}} = q_I \circ s.$$
Finally $q_I \circ s \in C^\infty(\R^d)$ since it is the composition of $q_I \in C_0^\infty(\R)$ and $s = \sum_{n =0}^N{u_n^2} \in C^\infty(\R^d)$ (since all the $u_n$ are in $C^\infty(\R^d)$) and $s \in [m, M]$.

\noindent{\bf Step 3.}

Finally, defining $v_n = 
\frac{u_n}{\sqrt{\sum_{n =0}^N{u_n^2}}}$ for all $0 \leq n \leq N$, $v_n \in C^\infty(\R^d)$ since it is the product of two infinitely differentiable functions. Moreover, $\sum_{n=0}^N v_i^2 = 1$ by construction and 
$v_0 = v_0\cdot \boldsymbol{1}_{\R^d \setminus{\bigcup_{n =1}^N B^{\prime}_n}}$ since $u_0$ satisfies the same equality and $v_0$ is the product of $u_0$ by the strictly positive function $1/\sqrt{s}$. Analogously $v_n = v_n\cdot \boldsymbol{1}_{B_n},~ \forall n \in \{1,\dots,N\}$, since $u_n$ satisfy the same equality and $v_n$ is the product of $u_n$ by the strictly positive function $1/\sqrt{s}$.
\epr

\section{Fundamental results on scattered data approximation}\label{sec:scattered-data}

We recall here some fundamental results about local polynomial approximation. In particular, we report here the proofs to track explicitly the constants. The proof techniques are essentially from \cite{narcowich2003refined,wendland2004scattered}. Denote by $\pi_k(\R^d)$ the set of multivariate polynomials of degree at most $k$, with $k \in \N$. In this section $B_r(x) \subset \R^d$ denotes the open ball of radius $r$ and centered in $x$.

\bp[\cite{wendland2004scattered}, Corollary 3.11. Local polynomial reproduction  on a ball]\label{prop:local-poly-repr} 
Let $k \in \Nz,~d,m \in \Np$ and $\delta >0$. Let $B_\delta$ be an open ball of radius $\delta > 0$ in $\R^d$. Let $\widehat{Y}  = \{y_1,\dots,y_m\} \subset B_\delta$ be a non empty finite subset of $B_\delta$. If either $k = 0$ or $h_{\widehat{Y},B_\delta} \leq \frac{\delta}{9k^2}$, there exist $u_j: B_\delta \to \R$ with $j \in [m]$ such that
\begin{enumprop}
\item\label{prop:local-poly-repr:a} $\sum_{j \in [m]} p(y_j) u_j(x) = p(x), \quad \forall x \in B_\delta, p \in \pi_k(\R^d)$
\item\label{prop:local-poly-repr:b} $\sum_{j \in [m]} |u_j(x)| \leq 2, \quad \forall x \in B_\delta$.
\end{enumprop}

\ep

\blm[Bounds on functions with scattered zeros on a small ball \cite{narcowich2003refined,wendland2004scattered}]\label{lm:zeros-scattered-data-small-ball}
Let $k \in \Nz,~d,m \in \Np$ and $\delta > 0$. Let $B_\delta \subset \R^d$ be a ball of radius $\delta$ in $\R^d$. Let $f \in C^{k+1}(B_\delta)$. Let $\widehat{Y}  = \{y_1,\dots,y_m\} \subset B_\delta$ be a non empty finite subset of $B_\delta$. If either $k =0$ or $h_{\widehat{Y}, B_\delta} \leq \frac{\delta}{9 k^2}$, it holds:
$$\sup_{x \in B_\delta} |f(x)| ~\leq~ 3C\delta^{k+1}~+~2\max_{i\in [m]} |f(y_i)|, \qquad  C := \sum_{|\alpha| = k+1} \frac{1}{\alpha!}  \|\partial^\alpha f\|_{L^\infty(B_\delta)}.$$
\elm

\bpr
 Note that since either $k= 0$ or $h_{\widehat{Y}, B_\delta} \leq \frac{\delta}{9 k^2}$, then we can apply \cref{prop:local-poly-repr} obtaining $u_j$ with $j \in [m]$ with the local polynomial reproduction property. Define the function $s_{f,\widehat{Y}} = \sum_{j \in [m]} f(y_j) u_j$ and let $\tau = \max_{i \in [m]} |f(y_i)|$.
Now, by using both \cref{prop:local-poly-repr:a,prop:local-poly-repr:b}, we have that for any $p \in \pi_k(\R^d)$ and any $x \in B_\delta$,
\eqals{
|f(x)| &\leq |f(x) - p(x)| + |p(x) - s_{f,\widehat{Y}}(x)| + |s_{f,\widehat{Y}}(x)| \\
&\leq |f(x) - p(x)| + \sum_{j \in [m]} |p(y_j) - f(y_j)| |u_j(x)| + \max_{j \in [m]} |f(y_j)| \sum_{j \in [m]} |u_j(x)| \\
&\leq \|f - p\|_{L^\infty(B_\delta)}(1 + \sum_{j \in [m]}|u_j(x)|) + \tau \sum_{j \in [m]} |u_j(x)|\\
&\leq 3\|f - p\|_{L^\infty(B_\delta)} + 2\tau.
}
In particular, consider the Taylor expansion of $f$ at the center $x_0$ of $B_\delta$ up to order $k$ (e.g. \cite{brenner2007mathematical} Eq.~4.2.5 pag 95). For any $x \in B_\delta$, it holds

$$f(x) = \sum_{|\alpha| \leq k} \frac{1}{\alpha!} \partial^\alpha f(x_0) (x -x_0)^\alpha ~+~ \sum_{|\alpha| = k+1} \frac{k+1}{\alpha!} (x - x_0)^\alpha \int_0^1 (1-t)^{k} \partial^\alpha f ((1-t)x_0 + tx) dt.$$

By choosing $p(x) = \sum_{|\alpha| \leq k} \frac{1}{\alpha!} \partial^\alpha f(x_0) (x -x_0)^\alpha ~ \in \pi_k(\R^d)$
it holds:
$$\|f - p\|_{L^\infty(B_\delta)} \leq \sum_{|\alpha| = k+1} \frac{\delta^{k+1}}{\alpha!}  \|\partial^\alpha f\|_{L^\infty(B_\delta)} = C\delta^{k+1},$$
where $C = \sum_{|\alpha| = k+1} \frac{1}{\alpha!}  \|\partial^\alpha f\|_{L^\infty(B_\delta)}$ is defined in the lemma. Gathering the previous equations,
$$\sup_{x \in B_\delta} |f(x)| \leq 2\tau + 3C \delta^{k+1}.$$
\epr

\bt[Bounds on functions with scattered zeros \cite{narcowich2003refined,wendland2004scattered}]\label{thm:zeros-scattered-data}
Let $k,m \in \Nz$ s.t. $k \leq m$ and $n,d \in \Np$. Let $r > 0$ and $\Omega$ an open set of $\R^d$ of the form $\Omega = \bigcup_{x \in S} B_r(x)$ for some subset $S$ of $\R^d$. Let $\widehat{X} = \{x_1,\dots,x_n\}$ be a non-empty finite subset of $\Omega$. Let $f \in C^{m+1}(\Omega)$. 
If $h_{\widehat{X}, \Omega} \leq r\max(1,\frac{1}{18k^2})$, then
$$\sup_{x \in \Omega} |f(x)| ~~\leq ~~ C C_f h_{\widehat{X}, \Omega}^{k+1} ~+~ 2 \max_{i \in [n]} |f(x_i)|,$$
where $C = 3\max(1,18~k^2)^{k+1}$ and
$C_f =  \sum_{|\alpha| = k+1} \frac{1}{\alpha!}  \|\partial^\alpha f\|_{L^\infty(\Omega)}.$
\et

\bpr
First, note that the condition that there exists a set $S$ such that $\Omega = \bigcup_{x \in S}B_r(x)$ implies
\[\forall \delta \leq r,~ \Omega = \bigcup_{x_0 \in S_\delta}{B_\delta(x_0)},\qquad S_\delta = \{x^\prime \in \Omega~:~ \exists x \in S,~ \|x-x^\prime\| \leq r-\delta\}.\]
We will now prove the theorem for $k \geq 1$ and then the easier case $k=0$, where we will use essentially only the Lipschitzianity of $f$.

\noindent{\bf Proof of the case $\boldsymbol{k \geq 1}$.} The idea of the proof is to apply \cref{lm:zeros-scattered-data-small-ball} to a collection of balls of radius $\delta$ for a well chosen $\delta \leq r$ and centered in $x_0 \in S_\delta$ defined above. Given $\widehat{X}$, to apply \cref{lm:zeros-scattered-data-small-ball} on a ball of radius $\delta$ we have to restrict the points in $\widehat{X}$ to the subset belonging to that ball, i.e.,  $\widehat{Y}_{x_0,\delta} = \widehat{X} \cap B_\delta(x_0)$, $x_0 \in S_\delta$ and $\delta >0$. The set $\widehat{Y}_{x_0,\delta}$ will have a fill distance $h_{x_0,\delta} = h_{\widehat{Y}_{x_0,\delta},B_\delta(x_0)}$. First we are going to show that $\widehat{Y}_{x_0,\delta}$ is not empty, when $r > \delta > h_{\widehat{X},\Omega}$. To obtain this result we need to study also the ball $B_{\delta'}(x_0)$ with $\delta' = \delta - h_{\widehat{X},\Omega}$.

\noindent{\bf Step 1. Showing that $\widehat{Y}_{x_0,\delta}$ is not empty and for any $y \in B_{\delta'}(x_0)$ there exists $z \in \widehat{Y}_{x_0,\delta}$ satisfying $\|y-z\| \leq h_{\widehat{X},\Omega}$.} 
Let $x_0 \in S_\delta$ and $\delta \leq r$. This implies that $B_\delta(x_0) \subseteq \Omega$ by the characterization of $\Omega$ in terms of $S_\delta$ we gave above. Define now $\delta' = \delta - h_{\widehat{X},\Omega}$ and note that $B_{\delta'}(x_0)$ is non empty, since $\delta^\prime > 0$, and that $B_{\delta'}(x_0) \subset B_\delta(x_0)  \subseteq \Omega$. Now note that by definition of fill distance, for any $y \in B_{\delta'}(x_0)$ there exists a $z \in \widehat{X}$ such that $\|z - y\| \leq h_{\widehat{X},\Omega}$. Moreover note that $z \in B_\delta(x_0)$, since $\|x_0 - z\| \leq \|x_0 - y\| + \|y - z\| < \delta - h_{\widehat{X},\Omega} + h_{\widehat{X},\Omega} = \delta$. Since $z \in \widehat{X}$ and also in $B_\delta(x_0)$, then $z \in \widehat{Y}_{x_0,\delta}$ by definition of $\widehat{Y}_{x_0,\delta}$.

\noindent{\bf Step 2. Showing that $h_{x_0,\delta} \leq 2h_{\widehat{X},\Omega}$.} Let $x \in B_{\delta}(x_0)$.  We have seen in the previous step that the ball $B_{\delta'}(x_0)$ is well defined and non empty, with $\delta' = \delta - h_{\widehat{X},\Omega}$. Now note that also $B_{h_{\widehat{X},\Omega}}(x) \cap B_{\delta'}(x_0)$ is not empty, indeed the distance between the centers $x,x_0$ is strictly smaller than the sum of the two radii, indeed $\|x-x_0\| < \delta = \delta' + h_{\widehat{X},\Omega}$, since $x \in B_{\delta}(x_0)$. Take $w \in B_{h_{\widehat{X},\Omega}}(x) \cap B_{\delta'}(x_0)$. Since $w \in B_{\delta'}(x_0)$ by Step~1 we know that there exists $z \in \widehat{Y}_{x_0,\delta}$ with $\|w-z\| \leq h_{\widehat{X},\Omega}$. Since $w \in B_{h_{\widehat{X},\Omega}}(x)$, then we know that $\|x-w\| < h_{\widehat{X},\Omega}$. So $\|x - z\| \leq \|x-w\| + \|w-z\| < 2 h_{\widehat{X},\Omega}$.

\noindent{\bf Step 3. Applying \cref{lm:zeros-scattered-data-small-ball}.} 
Since, by assumption $h_{\widehat{X},\Omega} \leq r/(18k^2)$ and $k \geq 1$, then the choice $\delta = 18 k^2 h_{\widehat{X},\Omega}$ implies $r \geq \delta > h_{\widehat{X},\Omega}$. So we can use the characterization of $\Omega$ in terms of $S_\delta$ and the results in the previous two steps, obtaining that for any $x_0 \in S_\delta$ the set $B_\delta(x_0) \subseteq \Omega$ and moreover the set $\widehat{Y}_{x_0,\delta}$ is not empty and covers $B_\delta(x_0)$ with a fill distance $h_{x_0,\delta} \leq 2h_{\widehat{X},\Omega}$. Since,
$h_{x_0,\delta} \leq 2h_{\widehat{X},\Omega} \leq \delta/(9k^2)$ then we can apply \cref{lm:zeros-scattered-data-small-ball} to each ball $B_\delta(x_0)$ obtaining 
$$\sup_{x \in B_\delta(x_0)} |f(x)| ~\leq~ 3C_{\delta,x_0}\delta^{k+1}~+~2\max_{z\in \widehat{Y}_{x_0,\delta}} |f(z)|, \qquad  C_{\delta,x_0} := \sum_{|\alpha| = k+1} \frac{1}{\alpha!}  \|\partial^\alpha f\|_{L^\infty(B_\delta(x_0))}.$$
The proof is concluded by noting that $\Omega = \bigcup_{x_0 \in S_\delta} B_\delta(x_0)$ and that for any $x_0 \in S_\delta$ we have $C_{\delta,x_0} \leq C_f$, $\delta^{k+1} \leq (18k^2)^{k+1} h_{\widehat{X},\Omega}^{k+1}$ and moreover that $\max_{z\in \widehat{Y}_{x_0,\delta}} |f(z)| \leq \max_{i \in [n]} |f(x_i)|$, since $\widehat{Y}_{x_0,\delta} \subseteq \widehat{X}$ by construction.

\noindent{\bf Proof of the case $\boldsymbol{k=0}$} 
Since $h_{\widehat{X},\Omega} \leq r$, by assumption, then $\delta = h_{\widehat{X},\Omega}$ implies that $\Omega$ admits a characterization as $\Omega = \bigcup_{x_0 \in S_\delta}B_{\delta}(x_0)$. Now let $x \in \Omega$ and choose $x_0 \in S_\delta$ such that $x \in B_\delta(x_0)$. One the one hand, since the segment $[x_0,x]$ is included in $\Omega$, by Taylor inequality, $|f(x) - f(x_0)| \leq C_f\|x-x_0\| \leq C_f h_{\widehat{X},\Omega}$ and $C_f = \sum_{|\alpha| = 1}{\frac{1}{\alpha !} \|\partial^\alpha f\|_{L^\infty(\Omega)}}$. One the other hand, by definition of $h_{\widehat{X},\Omega}$, there exists $z \in \widehat{X} \subset \Omega$ such that $\|z - x_0\| \leq h_{\widehat{X},\Omega} = \delta$. Since both the open segment $[x_0,z) \subset B_\delta(x_0) \subset \Omega$ and $z \in \Omega$, then the whole segment $[x_0,z] \subset \Omega$ and hence we can apply Taylor inequality to show $\|f(x_0) - f(z)\| \leq C_f\|z-x_0\| \leq C_f h_{\widehat{X},\Omega}$. Then we have 
$$|f(x)| \leq |f(x) - f(x_0)| + |f(x) - f(z)| + |f(z)| \leq  2 C_f h_{\widehat{X},\Omega} + \max_{i \in [n]}|f(x_i)|.$$ 
The proof of the step $k = 0$ is concluded by noting that the previous inequality holds for every $x \in \Omega$.
\epr

\section{Auxiliary results on RKHS}

We recall that the {\em nuclear norm} of a compact linear operator $A$ is defined as $\|A\|_\star = \tr(\sqrt{A^*A})$ or equivalently $\|A\|_\star = \sum_{j \in \N} \sigma_j$, where $(\sigma_j)_{j \in \N}$ are the singular values of $A$ (Chapter 7 of \cite{weidmann1980linear} or \cite{bhatia2013matrix} for the finite dimensional analogue).
\blm\label{lm:H-norm-tr-norm}
Let $\Omega$ be a set, $k$ be a kernel and $\hh $ the associated RKHS. Let $A: \hh  \to \hh $ be a trace class operator. If $\hh $ satisfies \cref{asm:H-rich:a}, then
$$\|r_A\|_{\hh } \leq \Cmul \|A\|_\star, \quad \textrm{where} \quad r_A(x) := \scal{\phi(x)}{A \phi(x)}, ~~ \forall x \in \Omega,$$
and $\|A\|_\star$ is the {\em nuclear norm} of $A$. We recall that if $A \in \pdm(\hh )$ then $\|A\|_\star = \tr(A)$.
\elm
\bpr
Since $A$ is compact, it admits a singular value decomposition $A = \sum_{i \in \N} \sigma_i u_i \otimes v_i$. Here, $(\sigma_j)_{j \in \N}$ is  a non-increasing sequence of non-negative eigenvalues converging to zero, and $(u_j)_{j \in \N}$ and $(v_j)_{j \in \N}$ are two orthonormal families of corresponding eigenvectors, (a family $(e_j)$ is said to be orthonormal if for $i,j \in \N$, $\scal{e_i}{e_j} = 1$ if $i=j$  and $\scal{e_i}{e_j} = 0$ otherwise) \cite{weidmann1980linear}.  Note that we can write $r_A$ using this decomposition as $r_A(x) = \sum_{i \in \N} \sigma_i u_i(x) v_i(x) = \sum_{i \in \N} \sigma_i \, (u_i \cdot v_i)(x)$, for all $x \in \Omega$, where we denote by $\cdot$ the pointwise multiplication between two functions (this equality is justified by the following absolute convergence bound). By \cref{asm:H-rich:a}, the fact that $A$ is trace-class (i.e.,  $\|A\|_\star < \infty$) and the fact that $u_j$, $v_j$ satisfy $\|u_j\|_{\hh }  = \|v_j\|_{\hh } = 1, j \in \N$, the following holds
\eqals{
\|r_A\|_{\hh } & = \|\sum_{j \in \N} \sigma_j (u_j \cdot v_j)\|_{\hh } \leq \sum_{j \in \N} \sigma_j \|u_j \cdot v_j\|_{\hh } \\
& \leq  \Cmul \sum_{j \in \N} \sigma_j \|u_j\|_{\hh }\|v_j\|_{\hh }  \leq \Cmul \sum_{j \in \N} \sigma_j = \Cmul \|A\|_\star.
}
In the case where $A \in \pdm(\hh )$, we have $\|A\|_\star = \tr(\sqrt{A^*A}) = \tr(A)$.
\epr

\subsection{Proof of \cref{prop:properties-V}}\label{proof:prop:properties-V}

Given the kernel $k$, the associated RKHS $\hh $ and the canonical feature map $\phi: \Omega \to \hh $ and a set of distinct points $\widehat{X} = \{x_1,\dots,x_n\}$ define the {\em kernel matrix} $K \in \R^{n \times n}$ as
$K_{i,j} = \scal{\phi(x_i)}{\phi(x_j)} = k(x_i, x_j)$ for all $i,j \in [n]$. Note that, since $k$ is a p.d. kernel, then $K$ is positive semidefinite, moreover when $k$ is universal, then $\phi(x_1), \dots, \phi(x_n)$ are linearly independent, so $K$ is full rank and hence invertible. Universality of $k$ is guaranteed since $\hh$ contains the $C^\infty_0(\Omega)$ functions, by \cref{asm:geom-f:a}, and so can approximate continuous functions over compacts in $\Omega$ \cite{steinwart2008support}. Denote by $R$ the upper triangular matrix corresponding to the Cholesky decomposition of $K$, i.e., $R$ satisfies $K = R^\top R$.
We are ready to start the proof of \cref{prop:properties-V}.

\bpr
Denote by $\widehat{S}:\hh  \to \R^n$ the linear operator that acts as follows 
$$\widehat{S} g ~=~ (\,\scal{\phi(x_1)}{g}\,,~\dots,~\scal{\phi(x_n)}{g}\,) \in \R^n, \qquad \forall g \in \hh .$$
Define $\widehat{S}^*:\R^n \to \hh $, i.e.,  the adjoint of $\widehat{S}$, as $\widehat{S}^*\beta = \sum_{i=1}^n \beta_i \phi(x_i)$ for $\beta \in \R^n$. Note, in particular, that $K = \widehat{S}\widehat{S}^*$ and that $\widehat{S}^* e_j = \phi(x_i)$, where $e_j$ is the $j$-th element of the canonical basis of $\R^n$. We define the operator $V = R^{-\top} \widehat{S}$ and its adjoint $V^* = \widehat{S}^*R^{-1}$. By using the definition of $V$, the fact that $K = R^\top R$ by construction of $R$, and the fact that $K = \widehat{S}\widehat{S}^*$, we derive two facts. 

\noindent On the one hand,
$$VV^* = R^{-\top} \widehat{S}\widehat{S}^* R^{-1} = R^{-\top} K R^{-1} = R^{-\top} R^\top R R^{-1} = I.$$
On the other hand, $P$ is a projection operator, i.e., $P^2 = P$, $P$ is positive definite and its range is $\ran P = \lspan{\phi(x_i)~|~i \in [n]}$, implying $P \phi(x_i) = \phi(x_i)$ for all $i \in [n]$. Indeed, using the equation above, $P^2 =  V^*VV^*V = V^*(VV^*)V  = V^* V = P$, and
the positive-semi-definiteness of $P$ is given by construction since it is the product of an operator and its adjoint. Moreover, the range of $P$ is the same as that of $V^*$ which in turn is the same as that of $S^*$, since $R$ is invertible : $\ran P = \lspan{\phi(x_i)~|~i \in [n]}$.

\noindent Finally, note that since $k(x,x') = \scal{\phi(x)}{\phi(x')}$, for any $x,x' \in \Omega$, then for any $j \in [n]$, $\Phi_j$ is characterized by 
\eqals{
\Phi_j &= R^{-\top}(k(x_1,x_j), \dots, k(x_n,x_j))\\& = R^{-\top}(\scal{\phi(x_1)}{\phi(x_j)}, \dots, \scal{\phi(x_n)}{\phi(x_j)}) = R^{-\top} \widehat{S}\phi(x_j) = V \phi(x_j).
}
\epr

\section{The constants of translation invariant and Sobolev kernels}

\subsection{Results for translation invariant and Sobolev kernels}
\blm\label{lm:pointwise-mult-tr-inv-kernel}
Let $\Omega$ be a set and let $k(x,x') = v(x-x')$ for all $x,x' \in \Omega$, be a translation invariant kernel for some function $v:\R^d \to \R$. Denote by $\tilde{v}$ the Fourier transform of $v$. Let $\hh $ be the associated RKHS. For any $f,g \in \hh $ we have
$$ \|f \cdot g\|_{\hh } \leq C \|f\|_{\hh }\|g\|_{\hh }, \qquad C = (2\pi)^{d/4}\left\|\frac{\tilde{v} \star \tilde{v}}{\tilde{v}}\right\|^{1/2}_{L^\infty(\R^d)}.$$
In particular, if there exists a non-increasing $g:[0,\infty] \to (0,\infty]$ s.t. $\tilde{v}(\omega) \leq g(\|\omega\|)$, then 
$$C \leq \sqrt{2}(2 \pi)^{d/2}v(0)^{1/2}\sup_{\omega \in \R^d} \sqrt{\frac{g(\tfrac{1}{2}\|\omega\|)}{\tilde{v}(\omega)}}.$$
\elm
\bpr
First note that by as recalled in \cref{ex:restriction-rkhs}, there exists an extension operator, i.e.,  a partial isometry $E: \hh  \to \hh(\R^d)$ such that $r = Eu$ satisfies $r(x) = u(x)$ for all $x \in \Omega$ and $\|u\|_{\hh } = \|r\|_{\hh }$, for any $u \in \hh $. Moreover there exists a restriction operator $R:\hh(\R^d)\to \hh $, as recalled in \cref{ex:restriction-rkhs}, such that $RE:\hh  \to \hh $ is the identity operator and $ER:\hh(\R^d) \to \hh(\R^d)$ is a projection operator whose range is $\hh $. Moreover, note that $f \cdot g= R(E f \cdot Eg)$ since for any $x \in \Omega$, $(R(E f \cdot Eg))(x) = (E f)(x)(Eg)(x) = f(x)g(x) = (f \cdot g)(x)$. Since $ER$ is a projection operator, then $\|ER\|_{\rm op} \leq 1$, hence
\eqals{
\|f \cdot g\|_{\hh } &= \|R (Ef \cdot Eg)\|_{\hh } = \|ER (Ef \cdot Eg) \|_{\hh(\R^d)} \\
&\leq \|ER\|_{\rm op} \|Ef \cdot Eg \|_{\hh(\R^d)} \leq \|Ef \cdot Eg \|_{\hh(\R^d)}.
}
Let $a = Ef$ and $b = Eg$. Denote by $\tilde{a},\tilde{b}$ their Fourier transform and by $\widetilde{a \cdot b}$ the Fourier transform of $a \cdot b$ (see \cref{prop:fourier} for more details). By expanding the definition of the Hilbert norm of translation invariant kernel
\eqals{
\|Ef \cdot Eg \|^2_{\hh(\R^d)} &= \| a \cdot b \|^2_{\hh(\R^d)} = (2\pi)^{-d/2}\int_{\R^d} \frac{|\widetilde{a \cdot b}\,(\omega)|^2}{\tilde{v}(\omega)} d\omega.
}
Now we bound $\widetilde{a\cdot b}$. Since $\widetilde{a\cdot b} = (2 \pi)^{d/2}\tilde{a} \star \tilde{b}$ (see \cref{prop:fourier}) where $\star$ corresponds to the convolution, by expanding it and by applying Cauchy-Schwarz we obtain
\eqals{
(2 \pi)^{-d/2}|\widetilde{a \cdot b}(\omega)|^2 &= |(\tilde{a} \star \tilde{b}) (\omega)|^2 = \left(\int_{\R^d} \tilde{a}(\sigma) \tilde{b}(\omega - \sigma) d\sigma\right)^2 \\
&= \left(\int_{\R^d} \frac{\tilde{a}(\sigma)}{\sqrt{\tilde{v}(\sigma)}} \frac{\tilde{b}(\omega - \sigma)}{\sqrt{\tilde{v}(\omega - \sigma)}} ~ \sqrt{\tilde{v}}(\sigma) \sqrt{\tilde{v}}(\omega - \sigma) d\sigma\right)^2\\
& \leq \int_{\R^d} \frac{\tilde{a}^2}{\tilde{v}}(\sigma) \frac{\tilde{b}^2}{\tilde{v}}(\omega - \sigma) d\sigma ~ \int_{\R^d} \tilde{v}(\sigma) \tilde{v}(\omega - \sigma) d\sigma= \left(\frac{\tilde{a}^2}{\tilde{v}} \star \frac{\tilde{b}^2}{\tilde{v}}\right)(\omega)~(\tilde{v} \star \tilde{v})(\omega).
}
By using the bound above together with Hölder inequality and Young inequality for convolutions, we have
\eqals{
(2\pi)^{-d/2}\int_{\R^d} \frac{|\widetilde{a \cdot b}\,(\omega)|^2}{\tilde{v}(\omega)} d\omega & \leq \int_{\R^d} \left(\frac{\tilde{a}^2}{\tilde{v}} \star \frac{\tilde{b}^2}{\tilde{v}}\right)(\omega)\, \frac{(\tilde{v} \star \tilde{v})(\omega)}{\tilde{v}(\omega)}d\omega  ~\leq~ \left\|\frac{\tilde{a}^2}{\tilde{v}} \star \frac{\tilde{b}^2}{\tilde{v}}\right\|_{L^1(\R^d)} \left\|\frac{\tilde{v} \star \tilde{v}}{\tilde{v}}\right\|_{L^\infty(\R^d)} \\
& \leq \left\|\frac{\tilde{a}^2}{\tilde{v}}\right\|_{L^1(\R^d)}  \left\|\frac{\tilde{b}^2}{\tilde{v}}\right\|_{L^1(\R^d)} \left\|\frac{\tilde{v} \star \tilde{v}}{\tilde{v}}\right\|_{L^\infty(\R^d)} \\
&= (2 \pi)^{d/2} \left\|\frac{\tilde{v} \star \tilde{v}}{\tilde{v}}\right\|_{L^\infty(\R^d)}\|a\|^2_{\hh(\R^d)}\|b\|^2_{\hh(\R^d)} = C^2,
}
where in the last step we used the definitions of inner products for translation invariant kernels. The proof is concluded by noting that $\|a\|_{\hh(\R^d)} = \|Ef\|_{\hh(\R^d)} = \|f\|_{\hh }$ and the same holds for $b$, i.e.,  $\|b\|_{\hh(\R^d)} = \|g\|_{\hh }$. A final consideration is that $C$ can be further bounded by applying \cref{prop:conv-ratio} and noting that $v(0) = (2\pi)^{-d/2}\int \tilde{v}(\omega) d\omega = (2\pi)^{-d/2}\|\tilde{v}\|_{L^1(\R^d)}$, via the characterization of $v$ in terms of $\tilde{v}$ in \cref{prop:fourier:Linfty-x}, since $\tilde{v}(\omega) \geq 0$ and integrable.
\epr

\bp\label{prop:conv-ratio}
Let $u \in L^1(\R^d) \cap C(\R^d)$ be $u(x) \geq  0$ for $x \in \R^d$ and such that there exists a non-increasing function $g: [0,\infty) \to (0,\infty)$ satisfying $u(x) \leq g(\|x\|)$ for all $x \in \R^d$.
Then it holds : 
\[\forall x \in \R^d,~ 0 \leq (u \star u)(x) \leq 2 \|u\|_{L^1(\R^d)}g(\tfrac{1}{2}\|x\|).\]
In particular, if $u > 0$, it holds
$$\left\|\frac{u \star u}{u}\right\|_{L^\infty(\R^d)} \leq 2\|u\|_{L^1(\R^d)}\sup_{x \in \R^d} \frac{g(\tfrac{1}{2}\|x\|)}{u(x)}.$$
\ep
\bpr 
For any $x \in \R^d$,
$$(u \star u)(x)= \sup_{x \in \R^d} \int_{\R^d} u(y)u(x-y) dy.$$
Let $S_x = \{y ~|~ \|x - y\| \leq \tfrac{1}{2} \|x\| \}$. Note that, when $y \in \R^d \setminus S_x$, then $\|x-y\| > \tfrac{1}{2}\|x\|$. Instead, when $y \in S_x$, then 
$$\tfrac{1}{2} \|x\| \leq \|x\| - \|x - y\| \leq \|y\|.$$
Since $g$ is non-increasing, for any $x \in \R^d$ we have
\eqals{
\int_{\R^d} u(y)u(x-y) dy &= \int_{S_x} u(y)u(x-y) dy ~+~ \int_{\R^d\setminus S_x} u(y)u(x-y) dy \\
& \leq \int_{S_x} g(\|y\|)u(x-y) dy ~+~ \int_{\R^d\setminus S_x} u(y)g(\|x-y\|)dy \\
& \leq \int_{S_x} g(\tfrac{1}{2}\|x\|)u(x-y)dy ~+~ \int_{\R^d\setminus S_x} u(y)g(\tfrac{1}{2}\|x\|) dy\\
& \leq \int_{\R^d} g(\tfrac{1}{2}\|x\|)u(x-y) dy ~+~ \int_{\R^d} u(y)g(\tfrac{1}{2}\|x\|) dy \\
& = \int_{\R^d} g(\tfrac{1}{2}\|x\|)u(y) dy ~+~ \int_{\R^d} u(y)g(\tfrac{1}{2}\|x\|) dy  = 2 ~g(\tfrac{1}{2}\|x\|) \int_{\R^d} u(y) dy,
}
where: in the first inequality we bounded $u(y)$ with $g(\|y\|)$ and $u(x-y)$ with $g(\|x-y\|)$, in the first and the second integral, respectively; in the second inequality we bounded $g(\|y\|)$ with $g(\tfrac{1}{2}\|x\|)$, since $\|y\| \geq \tfrac{1}{2}\|x\|$ when $y \in S_x$ and we bounded $g(\|x-y\|)$ with $g(\tfrac{1}{2}\|x\|)$, since $\|x-y\| \geq \tfrac{1}{2}\|x\|$ when $y \in \R^d\setminus S_x$; in the third we extended the integration domains to $\R^d$.  
\epr

\subsection{Proof of \cref{prop:sobolev-kernel}}\label{proof:prop:sobolev-kernel}

\bpr
We prove here that the Sobolev kernel satisfies \cref{asm:H-rich}. Let $k = k_s$ from \cref{eq:sobolev-kernel}. As we have seen in \cref{ex:sobolev-kernel} $\hh = W^s_2(\Omega)$ and $\|\cdot\|_{W^s_2(\Omega)}$ is equivalent to $\|\cdot\|_{\hh }$, when $s > d/2$ and $\Omega$ satisfies \cref{asm:geom-f:a} since this assumption implies that $\Omega$ satisfies the {\em cone condition} \cite{wendland2004scattered}. 

\noindent Recall that $k$ is translation invariant, i.e.,   $k(x,x') = v(x-x')$ for any $x,x' \in \R^d$, with $v$ defined in \cref{ex:sobolev-kernel}. The Fourier transform of $v$ is $\tilde{v}(\omega) = C_0(1+\|\omega\|^2)^{-s}$ with $C_0 = \frac{2^{d/2} \Gamma(s)}{\Gamma(s-d/2)}$  \cite{wendland2004scattered}. In the rest of the proof, $C_0$ will always refer to this constant.

\noindent We are going to divide the proof in one step per point of \cref{asm:H-rich}.

\noindent{\bf Proof of \cref{asm:H-rich:d} for the Sobolev kernel.} 
 Let $\alpha \in \Nz^d$,  $m = |\alpha|$. Assume $m < s-d/2$, i.e., $m \in \{1,\dots, \lfloor s - (d+1)/2\rfloor\}$. Since $k$ is translation invariant, then $\partial^\alpha_x \partial^\alpha_y k(x,y) = (-1)^m ~v_{2\alpha}(x-y)$ with $v_{2\alpha}(z) = \partial^{2\alpha}_z v(z)$ for all $z \in \R^d$. So
$$\sup_{x,y \in \Omega}|\partial^\alpha_x \partial^\alpha_y k(x,y)|  = \sup_{x,y \in \Omega}|\partial^\alpha_x \partial^\alpha_y v(x-y)| \leq \sup_{z \in \R^d}|\partial^{2\alpha}_z v(z)| \leq (2\pi)^{-d/2}\|\omega^{2\alpha} \tilde{v}(z)\|_{L^1(\R^d)},$$
where in the last step we used elementary properties of the Fourier transform (in particular the ones recalled in \cref{prop:fourier:derivative,prop:fourier:Linfty-x}).
Let $S_{d-1} = 2\frac{\pi^{d/2}}{\Gamma(d/2)}$ be the area of the $d-1$ dimensional sphere. Since $m < s - d/2$ and $\tilde{v} \geq 0$, 
\eqals{
	\|\omega^{2\alpha} \tilde{v}(z)\|_{L^1(\R^d)} & \leq  \int_{\R^d} \|\omega\|^{2m}\tilde{v}(\omega) d\omega  = C_0 S_{d-1} \int_0^\infty \frac{r^{2m + d - 1}}{(1+r^2)^s} dr \\
	& = C_0 S_{d-1} \int_0^\infty \frac{t^{m +d/2 -1}}{2(1+t)^s} dt = C_0 S_{d-1} \tfrac{\Gamma(m+d/2)\Gamma(s-d/2-m)}{2\Gamma(s)},
}
where we performed a change of variable $r = \sqrt{t}$ and $dr = \frac{dt}{2\sqrt{t}}$ and applied Eq.~5.12.3 pag. 142 of \cite{olver2010nist} to the resulting integral. Thus, \cref{asm:H-rich:d} holds with 
$$\Cdiff^2 = C_0 \frac{\pi^{d/2} \Gamma(m+d/2) \Gamma(s- m -d/2)}{\Gamma(d/2)\Gamma(s)} =  \frac{(2\pi)^{d/2} \Gamma(m+d/2)\Gamma(s-d/2-m)}{\Gamma(s-d/2)\Gamma(d/2)}.$$

\noindent{\bf Proof of \cref{asm:H-rich:a} for the Sobolev kernel.} First, note that $C^\infty(\R^d)|_\Omega \subset W^s_\infty(\Omega) \subset W^s_2(\Omega)$. Indeed, since $\Omega$ is bounded, for any $f \in C^\infty(\R^d)$, $\|\partial^\alpha f|_\Omega\|_{L^\infty(\Omega)} < \infty$ for any $\alpha \in \Nz^d$. This shows that $f|_\Omega \in W^s_\infty(\Omega)$. Moreover $W^s_\infty(\Omega) \subset W^s_2(\Omega)$ since $\|\cdot\|_{L^2(\Omega)} \leq \vol(\Omega)^{1/2} \|\cdot\|_{L^\infty(\Omega)}$ because $\Omega$ is bounded. Second, since $\tilde{v}(\omega) = g_s(\|\omega\|)$ with $g_s(t) = C_0 (1 + t^2)^{-s}$, positive and non-increasing, we can apply \cref{lm:pointwise-mult-tr-inv-kernel}. Therefore, for $C = \sqrt{2}(2\pi)^{d/2}v(0)^{1/2} \sup_{t \geq 0}\big(\tfrac{g_s(t/2)}{g_s(t)}\big)^{1/2}$ it holds $\|f \cdot g\|_{\hh } \leq C \|f\|_{\hh }\|g\|_{\hh }$. In particular we have
$\sup_{t \geq 0}\big(\tfrac{g_s(t/2)}{g_s(t)}\big)^{1/2} \leq 2^{s}$ and $v(0) = 1$, since $\lim_{t\to 0} t^{s-d/2}{\cal K}_{s-d/2}(t) = \Gamma(s-d/2)/2^{1+d/2-s} = 1/C_0$ (\cite{olver2010nist} Eq.~10.30.2 pag. 252) and $v(x) = C_0 t^{s-d/2}{\cal K}_{s-d/2}(t),~t = \|x\|$. Thus, \cref{asm:H-rich:a} holds with constant 
$$\Cmul =\pi^{d/2}  2^{(2s+d+1)/2}.$$ 

\noindent{\bf Proof of \cref{asm:H-rich:b} for the Sobolev kernel.} First we recall from \cite{adams2003sobolev} that for any $s  > d/2$, there exists a constant $C_s$ such that 
$$ \forall h \in W^s_2(\R^d),~ \|h\|_{L^\infty(\R^d)} \leq C_s \|h\|_{ W^s_2(\R^d)}.$$
In particular, this shows that $W^s_2(\R^d) \subset L^\infty(\R^d)$. Fix such a constant $C_s$ in the rest of the proof.

\noindent Let $p \in \N$ and $g \in C^\infty(\R^p)$ with $g(0,0,\dots,0) = 0$. From (i) of Thm.~11 in \cite{sickel1992superposition}, there exists a constant $c_g$ depending only on $g,p,s$ such that for any $h_1,\dots,h_p \in W^s_2(\R^d)\cap L^\infty(\R^d)$, it holds
\[ \|g(h_1,\dots,h_p)\|_{W^s_2(\R^d)} ~~\leq~~ c_g \sup_{i \in [p]}~\|h_i\|_{W^s_2(\R^d)}\left(1  +\|h_i\|_{L^\infty(\R^d)}^{\max(0,s-1)}\right).\]
Since $s > d/2$, the bound above shows, in particular, that for any $h_1,\dots,h_p \in W^s_2(\R^d)$, it holds
\[ \|g(h_1,\dots,h_p)\|_{W^s_2(\R^d)} ~~\leq~~ c^\prime_g \sup_{i \in [p]}~\left(\|h_i\|+\|h_i\|_{W^s_2(\R^d)}^{\max(1,s)}\right),\qquad c^\prime_g = c_g\max\left(1 , C_s^{\max(0,s-1)}\right).\]

Since $W^s_2(\R^d) = \hh(\R^d)$ and $\|\cdot\|_{W^s_2(\R^d)} $ and $ \|\cdot\|_{\hh(\R^d)}$ are equivalent (see \cite{adams2003sobolev}), the previous inequality holds for $\|\cdot\|_{\hh(\R^d)}$ with a certain constant $c^\prime_g$ depending only on $g,p,s,d$. In particular, this implies that $g(h_1,\dots,h_p) \in \hh(\R^d)$ for any $h_1,\dots, h_p \in \hh(\R^d)$. Now we are going to prove the same implication for the restriction on $\Omega$.

\noindent First note that any function in $a \in C^\infty(\R^p)$ can be written as $a(z) = q\,1(z) + g(z)$, $z \in \R^p$ where $q = a(0,0,\cdots,0) \in \R$, $g \in C^\infty(\R^p)$ with $g(0,0,\cdots,0) = 0$ and $1(z) = 1$ for all $z \in \R^p$. Recall the definition and basic results on the extension operator $E:\hh  \to \hh(\R^d)$ from \cref{ex:restriction-rkhs}. For any $f_1,\dots, f_p \in \hh $, note that
$g((Ef_1)(x),\dots,(Ef_p)(x)) = g(f_1(x),\dots,f_p(x))$ for all $x \in \Omega$. We can now apply the results of \cref{ex:restriction-rkhs} to show that $g(f_1,\dots,f_p) \in \hh $ : 
\eqals{
\|g(f_1,\dots,f_p)\|_{\hh } &= \inf_{u} \|u\|_{\hh(\R^d)}~~s.t.~~ u(x) = g(f_1(x),\dots,f_p(x))~\forall x \in \Omega\\
& \leq \|g(Ef_1,\dots,Ef_p)\|_{\hh(\R^d)} \\
& \leq c^\prime_g \sup_{j \in [p]} \|E f_j\|_{\hh(\R^d)} + \|E f_j\|_{\hh(\R^d)}^{\max(1,s)} \\
& = c^\prime_g \sup_{j \in [p]} \|f_j\|_{\hh} +  \| f_j\|_{\hh}^{\max(1,s)} < \infty,
}
where in the last step we used the fact that $\|\cdot\|_{\hh } = \|E \,\cdot\|_{\hh(\R^d)}$.
The proof of this point is concluded by noting that, $a(f_1,\dots,f_p) \in \hh $, since $1 \in \hh $, due to the Point (a) above, and
$$\|a(f_1,\dots,f_p)\|_{\hh } \leq  q\|1\|_{\hh } + \|g(f_1,\dots,f_p)\|_{\hh } < \infty.$$

\noindent{\bf Proof of \cref{asm:H-rich:c} for the Sobolev kernel.} This proof is done in \cref{lm:soboc}, right below.
\epr

Before stating \cref{lm:soboc} we are going to recall some properties. First, recall the Young inequality : 
\[\forall f \in L^2(\R^d),~\forall g \in L^1(\R^d),~ \|f\star g\|_{L^2(\R^d)} \leq \|f\|_{L^2(\R^d)} ~ \|g\|_{L^1(\R^d)}.\]
\noindent Moreover, by definition of the Sobolev kernel, it is a translation-invariant kernel with $v$ defined in \cref{ex:sobolev-kernel}, with Fourier transform $\tilde{v}(\omega) = C_0 (1+\|\omega\|^2)^{-s}$. Let $\hh(\R^d)$ be the reproducing kernel Hilbert space on $\R^d$ associated to the Sobolev kernel $k_s$. As recalled in \cref{ex:tr-inv-kernel}, the $\hh(\R^d)$-norm is characterized by
\eqal{\label{eq:Hfourier}\forall f \in \hh(\R^d),~ \|f\|_{\hh(\R^d)} = (2\pi)^{-d/4} \|\tilde{f}/\sqrt{\tilde{v}}\|_{L^2(\R^d)},}
where $\tilde{f} = {\cal F}(f)$ is the Fourier transform of $f$ (see \cite{adams2003sobolev}).  Then we recall that $\tilde{v} \in L^1(\R^d)$, since $s > d/2$, so for any $f \in \hh(\R^d)$
\eqal{\label{eq:L1F-H-bound}
\|\tilde{f}\|_{L^1(\R^d)} = \|\sqrt{\tilde{v}} \tilde{f}/\sqrt{\tilde{v}}\|_{L^1(\R^d)} \leq \|\sqrt{\tilde{v}}\|_{L^2(\R^d)}\|\tilde{f}/\sqrt{\tilde{v}}\|_{L^2(\R^d)} = C_1 \|f\|_{\hh(\R^d)}.
}
where $C_1 = (2\pi)^{d/4} \|\sqrt{\tilde{v}}\|_{L^2(\R^d)}$. A useful consequence of the inequality above is obtained by considering that $\|f\|_{L^\infty(\R^d)}$ is bounded by the $L^1$ norm of $\tilde{f}$ (see \cref{prop:fourier:Linfty-x}), then
\eqal{\label{eq:Linfty-H-bound}
\|f\|_{L^\infty} \leq (2\pi)^{-d/2} \|\tilde{f}\|_{L^1(\R^d)} \leq  C_2 \|f\|_{\hh(\R^d)},
}
where $C_2 = (2\pi)^{-d/4} \|\sqrt{\tilde{v}}\|_{L^2(\R^d)}$.

\blm[\cref{asm:H-rich:c} for Sobolev Kernels]\label{lm:soboc} Let $\hh$ be the RKHS associated to the translation invariant Sobolev Kernel defined in \cref{ex:sobolev-kernel}, with $s > d/2$. Then \cref{asm:H-rich:c} is satisfied.
\elm 

\bpr
For the rest of the proof we fix $u : \Omega \rightarrow \R$ with $u \in \hh$, $r >0$ and $z \in \R^d$ such that $B_r(z) \subset \Omega$. Let $E_\Omega : \hh \rightarrow \hh(\R^d)$ be the extension operator from $\Omega$ to $\R^d$ (its properties are recalled in \cref{ex:restriction-rkhs}).Let $\chi \in C^\infty_0(\R^d)$ be given by \cref{lm:bump} such that $\chi = 1$ on $B_r(z)$, $\chi = 0$ on $\R^d \setminus{B_{2r}(z)}$ and $\chi \in [0,1]$. Define for any $t \in \R$ and $x \in \R^d$
$$ h_t(x) = \chi(x) w_t(x), \qquad w_t(x) = w((1-t)z + tx), \qquad w = E_\Omega u.$$
In particular we recall that, since $E_\Omega$ is a partial isometry (see \cref{ex:restriction-rkhs}) then $\|w\|_{\hh(\R^d)} = \|u\|_{\hh}$.
\noindent{\bf Step 1. Fourier transform of $w_t$.}
Denote with $\widetilde{w}$ the Fourier transform of $w$ which is well defined since $w \in \hh(\R^d) \subset L^2(\R^d)$ (see \cite{adams2003sobolev}),  with $\tilde{\chi}$ the Fourier transform of $\chi$. Since For any $t \neq 0$, denote with $\widetilde{w_t}$ the Fourier transform of $w_t$ which is well defined using the results of \cref{prop:fourier}, and which satisfies
\[\forall t \neq 0,~\forall \omega \in \R^d,~ \widetilde{w_t}(\omega) = |t|^{-d} e^{i\frac{1-t}{t}z^\top \omega}\tilde{w}(\omega/t).\]

\noindent {\bf Step 2. Separating low and high order derivatives of $h_t$, and bounding the low order terms.} 
For $t \neq 0$, denote with $\widetilde{h_t}$ the Fourier transform of $h_t$ which is well defined since $\chi$ is bounded and $w_t \in L^2(\R^d)$. We will now bound $\|h_t\|_{\hh(\R^d)}$ for all $t \neq 0$, by using the characterization in \cref{eq:Hfourier}.
Since $(x+y)^s \leq 2^{\max(s-1,0)}(x^s + y^s)$ for any $x,y \geq 0$, $s \geq 0$, then $(1 + \|\omega\|^2)^{s/2} \leq c_1(1+\|\omega\|^s)$ for any $\omega \in \R^d$, with $c_1 = 2^{\max(s/2-1,0)}$ so using \cref{eq:Hfourier}, we have
\[\sqrt{C_0}(2\pi)^{d/4}\|h_t\|_{\hh(\R^d)} = \|(1 + \|\cdot\|^2)^{s/2} \widetilde{h_t}\|_{L^2(\R^d)} \leq c_1\,\|\widetilde{h_t}\|_{L^2(\R^d)} + c_1\,\| ~|\cdot|^s_{\R^d} \widetilde{h_t}\|_{L^2(\R^d)}.\]
The first term on the right hand side can easily be bounded using the fact that the Fourier transform is an isometry of $L^2(\R^d)$ (see \cref{prop:fourier} for more details), indeed
\eqals{
\|\widetilde{h_t}\|_{L^2(\R^d)} = \|h_t\|_{L^2(\R^d)} & = \|\chi\cdot w_t\|_{L^2(\R^d)} \leq  \|w_t\|_{L^\infty(\R^d)} \|\chi\|_{L^2(\R^d)} < \infty.
}
since $\chi \in C^\infty_0(\R^d)$ by definition, so it it bounded and has compact support, implying that $\|\chi\|_{L^2(\R^d)} < \infty$, moreover $\|w_t\|_{L^\infty(\R^d)} = \|w\|_{L^\infty(\R^d)}$ and $\|w\|_{L^\infty(\R^d)} \leq C_2 \|w\|_{\hh(\R^d)}$ as recalled in \cref{eq:Linfty-H-bound} (the constant $C_2$ is defined in the same equation).

\noindent {\bf Step 3. Decomposing the high order derivatives of $h_t$.} Note that since $\widetilde{h_t} = \widetilde{\chi \cdot w_t}$, by property of the Fourier transform (see \cref{prop:fourier:product}), $\widetilde{\chi \cdot w_t} = (2\pi)^{d/2} \widetilde{\chi} \star \widetilde{w_t}$. Moreover, since $\|\omega\|^s \leq (\|\omega-\eta\| + \|\eta\|)^s \leq c_s(\|\omega - \eta\|^s + \|\eta\|^s)$ for any $\omega,\eta \in \R^d$, with $c = 2^{\max(s-1,0)}$, then, for all $\omega \in \R^d$ we have
\eqals{
\|\omega\|^s &|\widetilde{h_t}(\omega)| = \|\omega\|^s |\widetilde{\chi \cdot w_t}(\omega)| = \|\omega\|^s (2\pi)^{\frac{d}{2}} |(\tilde{\chi} \star \tilde{w}_t)(\omega)| = (2\pi)^{\frac{d}{2}} |\int_{\R^d} \|\omega\|^s \widetilde{\chi}(\eta) \widetilde{w_t}(\omega - \eta) d \eta |\\ 
& \leq (2\pi)^{\frac{d}{2}} c\int_{\R^d}(|\widetilde{\chi}(\eta)| ~\|\eta\|^s)~|\widetilde{w_t}(\omega - \eta)| ~d\eta  ~+~  (2\pi)^{\frac{d}{2}} c \int_{\R^d}|\widetilde{\chi}(\eta)| ~(|\widetilde{w_t}(\omega - \eta)|~\|\omega - \eta\|^{s}) ~d\eta \\
& = c\, ((J_s |\widetilde{\chi}|) \star |\widetilde{w_t}|)(\omega) ~+~ c\, (|\tilde{\chi}| \star (J_s |\widetilde{w_t}|))(\omega),
}
where we denoted by $J_s$ the function $J_s(\omega) = \|\omega\|^s$ for any $\omega \in \R^d$. Applying Young's inequality, it holds : 
\eqals{
\| J_s \widetilde{h_t}\|_{L^2(\R^d)} &\leq c\, \|(J_s |\widetilde{\chi}|) \star |\widetilde{w_t}|\|_{L^2(\R^d)} ~+~ c\||\tilde{\chi}| \star (J_s |\widetilde{w_t}|)\|_{L^2(\R^d)} \\
&\leq c \|J_s \widetilde{\chi}\|_{L^2(\R^d)} \|\widetilde{w}_t\|_{L^1(\R^d)} + c\|J_s \widetilde{w}_t \|_{L^2(\R^d)}~\|\widetilde{\chi}\|_{L^1(\R^d)}.
}
\noindent{\bf Step 4. Bounding the elements of the decomposition.}
Now we are ready to bound the four terms of the decomposition of $\| J_s \widetilde{h_t}\|_{L^2(\R^d)}$. First term, since $\chi \in C^\infty_0(\R^d) \subset \hh(\R^d)$, and $J_s(\omega) \leq \sqrt{C_0/\tilde{v}(\omega)}$ for any $\omega \in \R^d$, then $\|J_s \widetilde{\chi}\|_{L^2(\R^d)} \leq \sqrt{C_0} \|\widetilde{\chi}/\sqrt{\tilde{v}}\|_{L^2(\R^d)} =(2\pi)^{d/4}\sqrt{C_0}\|\chi\|_{\hh(\R^d)}$, where we used \cref{eq:Hfourier}.
Second term, $\|\widetilde{\chi}\|_{L^1(\R^d)} < \infty$, since $\|\widetilde{\chi}\|_{L^1(\R^d)} \leq C_1 \|\chi\|_{\hh(\R^d)}$, via \cref{eq:L1F-H-bound} (the constant $C_1$ is defined in the same equation) and we have seen already that $\|\chi\|_{\hh(\R^d)}$ is bounded. Third term, by a change of variable $\tau = \omega/t$,
\eqals{
\|\widetilde{w}_t\|_{L^1(\R^d)} &= \int_{\R^d}|\widetilde{w}_t(\omega)| d\omega = \int_{\R^d}|t|^{-d}|\tilde{w}(\omega/t)| d\omega = \int_{\R^d} |\tilde{w}(\tau)| d\tau = \|\tilde{w}\|_{L^1(\R^d)},
}
moreover $\|\tilde{w}\|_{L^1(\R^d)} \leq C_1 \|w\|_{\hh(\R^d)} =  C_1 \|u\|_{\hh}$ via \cref{eq:L1F-H-bound} and the fact that $\|w\|_{\hh(\R^d)} =  \|u\|_{\hh}$ as recalled at the beginning of the proof.
Finally, fourth term, for $t  \in \R \setminus{\{0\}}$,
\eqals{
\|J_s \widetilde{w}_t\|^2_{L^2(\R^d)} &= \int_{\R^d}\|\omega\|^{2s} |\widetilde{w_t}(\omega)|^2~d\omega = t^{-2d}\int_{\R^d}\|\omega\|^{2s} |\widetilde{w}(\omega/t)|^2 d\omega\\
&= t^{2s-d} \int_{\R^d}{\|\tau\|^{2s}|\widetilde{w}(\tau)|^2 d\tau} \leq t^{2s-d} \int_{\R^d}{(1+\|\tau\|^2)^{s}|\widetilde{w}(\tau)|^2 d\tau} \\
& = t^{2s-d} (2\pi)^{d/2}C_0\|w\|^2_{\hh(\R^d)}.
}
where we performed a change of variable $\omega = t\,\tau$, $t^d d\tau = d\omega$ and used the definition in \cref{eq:Hfourier} and the fact that $\|\tau\|^{2s} \leq (1+\|\tau\|^2)^s$ for any $\tau \in \R^d$. The proof of the bound of the fourth term is concluded by recalling that $\|w\|_{\hh(\R^d)} = \|u\|_{\hh}$ as discussed in the proof of the bound for the previous term.

\noindent{\bf Conclusion.} Putting all our bounds together, we get : 
\[\forall t \in \R\setminus{\{0\}},~ \|h_t\|_{\hh(\R^d)} \leq (A + B~t^{s-d/2}) \|\chi\|_{\hh(\R^d)}\|u\|_{\hh},\]
where $A = c_1 C_2 + c c_1 C_1 (2\pi)^{d/4}\sqrt{C_0}$ and $B = c c_1 C_1 (2\pi)^{d/4} \sqrt{C_0}$, where $c = 2^{\max(s-1,0)}$, $c_1 = 2^{\max(s/2-1,0)}$, while $C_1$ is defined in \cref{eq:L1F-H-bound}, $C_2$ in \cref{eq:Linfty-H-bound}.
\noindent Now define
$$ \forall x \in \R^d,~ \overline{g}_{z,r}(x) = \int_0^1 (1-t) h_t(x) dt,$$
and note that, by construction $\overline{g}_{z,r}(x) = \int_0^1 (1-t) u(t z + (1-t)x) dt$ for any $x \in B$ since $u$ and $\chi w$ coincide on $B$.
\noindent Note that the map $t \in (0,1) \mapsto (1-t)\|h_t\|_{\hh(\R^d)}$ is measurable, using the expression in \cref{eq:Hfourier}. 
\noindent Moreover, since for all $t \in (0,1)$, it holds $\|h_t\|_{\hh(\R^d)} \leq (A + B t^{s-d/2})\|\chi\|_{\hh(\R^d)}\|u\|_{\hh} \leq (A + B)\|\chi\|_{\hh(\R^d)}\|u\|_{\hh}$ since $s > d/2$, 
the map $t \mapsto (1-t)h_t$ is in integrable, and thus
\eqals{\|\overline{g}_{z,r}\|_{\hh(\R^d)}  &=  \big\|\int_0^1 (1-t) h_t dt\big\|_{\hh(\R^d)}\leq \int_0^1 |1-t| \|h_t\|_{\hh(\R^d)} dt \leq (A + B)\|\chi\|_{\hh(\R^d)}\|u\|_{\hh}   < \infty,}
\noindent which implies that the function $\overline{g}_{z,r}$ belongs to $\hh(\R^d)$. 
Finally, denote by $R_\Omega:\hh(\R^d) \to \hh$ the restriction operator (see \cref{ex:restriction-rkhs} for more details). By construction $(R_\Omega g)(x) = g(x)$ for any $g \in \hh(\R^d)$ and $x \in \Omega$, defining $g_{z,r} = R_\Omega \overline{g}_{z,r}$ the lemma is proven.
\epr

\section{Proofs for \cref{alg:glm}}

We start with two technical lemmas that will be used by the proofs in this section.
\blm[Technical result]\label{lm:technical}
Let $\alpha \geq 1$, $\beta \geq 2$ and $n \in \Nz$. If $n \geq 2\alpha \log(2\beta \alpha)$, then it holds
\[\frac{\alpha \log(\beta n)}{n} \leq 1.\]
\elm 
\bpr 
  Note that the function $x \mapsto \frac{\log (\beta x)}{x}$ is strictly decreasing on $[\exp(1)/\beta,+\infty]$.
  
  \noindent Moreover, $2\alpha \log(2\beta \alpha) \geq 2\log 4 \geq \exp(1)/2 \geq \exp(1)/\beta$ since $\beta \geq 2$ and $\alpha \geq 1$.
  
  \noindent Now assume $n \geq c \alpha $ with $c = 2\log(2\beta \alpha)$. It holds:
  \[\frac{\alpha \log(\beta n)}{n}  \leq  \frac{\log (\beta c \alpha)}{c} \leq \frac{\log (\tfrac{c}{2}) + \log(2\alpha \beta)}{c} \leq \frac{1}{2} + \frac{1}{2}~\frac{2\log(2\beta \alpha)}{c} \leq 1,\]

\noindent where we used the definition of $c$ and the fact that $\log(c/2) \leq c/2 -1 \leq c/2$.
\epr

\blm\label{lm:cone}
Let $\overrightarrow{u} \in S_{d-1} = \{x \in \R^{d}~|~\|x\| = 1\}$, $\alpha \in [0,\pi/2]$, $x_0 \in \R^d$ and $t > 0$. Define the cone centered at $x_0$, directed by $\overrightarrow{u}$ of radius $t$ with aperture $\alpha$: 
\[C_{x_0,\overrightarrow{u},t}^\alpha = \left\{x \in B_t(x_0)~|~ \tfrac{x-x_0}{\|x-x_0\|} \cdot \overrightarrow{u} \leq \cos(\alpha),~x \neq x_0 \right\},\]
where we denoted by $\cdot$ the scalar product among vectors. 
Then the volume of this cone is lower bounded as 
\[\vol(C_{x_0,\overrightarrow{u},t}^\alpha) \geq \frac{(\sqrt{\pi}\sin(\alpha))^{d-1}(t~\cos \alpha )^d}{d\Gamma((d+1)/2)}. \]

\noindent Moreover, let $x_0 \in \R^d$ and $r >0$. Let $x \in B_r(x_0)$ and $0 < t \leq r$. The intersection $B_t(x) \cap B_r(x_0)$ contains the cone $C_{x,\overrightarrow{u},t}^{\pi/3}$, where $\overrightarrow{u} = \frac{x-x_0}{\|x-x_0\|}$ if $x \neq x_0$ and any unit vector otherwise.
\elm 
\bpr 
\noindent{\bf 1. Bound on the volume of the cone}. Without loss of generality, assume $x_0 = 0$ and $\overrightarrow{u} = e_1$ since the Lebesgue measure is invariant by translations and rotations. A simple change of variable also shows that $\vol(C_{0,\overrightarrow{u},t}^\alpha)  = t^d \vol(C_{0,\overrightarrow{u},1}^\alpha)$. Now note the following inclusion (the proof is trivial): 
\[\widetilde{C} := \left\{x = (x_1,z) \in \R^d = \R \times \R^{d-1}~:~z \leq \cos(\alpha),~\|z\|_{\R^{d-1}} \leq x_1\sin(\alpha)\right\} \subset C_{0,e_1,1}^\alpha.\]
It is possible to compute the volume of the left hand term explicitly : 
\eqals{ 
\vol(\widetilde{C}) &= \int_{\R} {\boldsymbol{1}_{x_1 \leq \cos(\alpha)} \left( \int_{\R^{d-1}}{\boldsymbol{1}_{\|z\| \leq x_1 \sin(\alpha)} dz}\right) dx_1} \\
&= \int_{0}^{\cos(\alpha)}{V_{d-1}(\sin \alpha x_1)^{d-1}~dx} \\
&= V_{d-1}\frac{\sin^{d-1}(\alpha) \cos^d(\alpha)}{d} ,
}
where $V_{d-1} = \pi^{(d-1)/2}/\Gamma((d-1)/2 +1)$ denotes the volume of the $d-1$ dimensional ball.

\noindent{\bf 2. Proof of the second point} The case where $x = x_0$ is trivial since $t \leq r$. Assume therefore $x \neq x_0$ and note that by definition, $C_{x,\overrightarrow{u},t}^{\pi/3} \subset B_t(x)$. We will now show that  $C_{x,\overrightarrow{u},t}^{\pi/3} \subset B_r(x_0)$. Let $y \in C_{x,\overrightarrow{u},t}^{\pi/3}$ and assume $y \neq x$ (if $y=x$ then $y \in B_r(x_0)$). Expanding the dot product

\eqals{ \|y-x_0\|^2 &= \|y - x\|^2 + 2(y-x)\cdot (x-x_0) + \|x-x_0\|^2 \\
& =  \|y - x\|^2 - 2\|y-x\|~\|x_0-x\|~\tfrac{y-x}{\|y-x\|} \cdot \overrightarrow{u} + \|x-x_0\|^2\\
& \leq \|x-y\|^2 - \|x-y\|~\|x-x_0\| + \|x-x_0\|^2,
}
where the last inequality comes from the definition of the cone and $\cos \pi/3 =\frac{1}{2}$. Let us distinguish two cases:

\begin{itemize}
    \item if $t > \|x_0-x\|$, we have $- \|x-y\| \|x_0-x\| \leq - t^2$ and hence $\|y-x_0\|^2 \leq t^2 \leq r^2$;
    \item otherwise $\|x-y\| \leq t \leq \|x_0-x\|$ and thus  $\|y-x_0\|^2 \leq \|x-x_0\|^2 \leq r^2$.
\end{itemize}
In any case, $y \in B_r(x_0)$, which concludes the proof. 
\epr

\subsection{Proof of \cref{thm:fill-distance-random-points}} \label{proof:thm:fill-distance-random-points}

\bpr[Proof of \cref{thm:fill-distance-random-points}] Fix $\Omega$ as in \cref{thm:fill-distance-random-points}.
Let $U$ be the uniform probability over $\Omega$, i.e.,  $U(A) = \tfrac{\vol(A \cap \Omega)}{\vol(\Omega)}$ for any Borel-measurable set $A$. Let $\mathbb{P} = U^{\otimes n}$ over $\Omega^n$. Throughout this proof, we will use the notation $V_d$ to denote the volume of the $d$-dimensional unit ball (recall that $V_d = \tfrac{\pi^{d/2}}{\Gamma(d/2+1)}$).

\noindent{\bf Step 1. Covering $\Omega$.} Let $t > 0$. We say that a subset $\overline{X}$ of $\Omega$ is a $t$ (interior) covering of $\Omega$ if $\Omega \subset \bigcup_{x \in \overline{X}}{B_t(x)}$. Denote with $N_t$ the minimal cardinal $|\overline{X}|$ of a $t$ interior covering of $\Omega$ and fix $\overline{X}_t$ a $t$ interior covering of $\Omega$ whose cardinal is minimum, i.e.,  $|\overline{X}_t| = N_t$. Since the diameter of $\Omega$ is bounded by $2R$, it is known that $N_t \leq (1 + 2R/t)^d$ 

\noindent To prove this fact , one defines a maximal $t/2$-packing of $\Omega$ as a maximal set $\overline{Y}_{t/2} \subset \Omega$ such that the balls $B_{t/2}({\overline{y}})$ are disjoint. It is then easy to check that if $\overline{Y}_{t/2}$ is a maximal $t/2$-packing, then it is also a $t$-covering and hence $N_t \leq |\overline{Y}_{t/2}|$. Finally, since $\Omega$ is included in a ball of radius $B_{2R}(x_0)$ for some $x_0 \in \R^d$ and since $\overline{Y}_{t/2} \subset \Omega$, it holds $\bigcup_{\overline{y} \in \overline{Y}_{t^{\prime}}}B_t(\overline{y}) \subset B_{R+t/2}(x_0)$. Since the $B_t(\overline{y})$ are two by two disjoint, the result follows from the following equation:
$$|\overline{Y}|_{t/2} (t/2)^d V_d = \vol\left(\cup_{\overline{y} \in \overline{Y}_{t^{\prime}}}B_t(\overline{y})\right) \leq \vol(B_{R+t/2}(x_0)) = (R + t/2)^d V_d.$$

\noindent{\bf Step 2. Probabilistic analysis.}  Note that for any $(x_1,...,x_n) \in \Omega^n$, writing $\widehat{X} = \{x_1,..,x_n\}$, it holds:
\eqals{ 
	h_{\widehat{X},\Omega} &= \max_{x \in \Omega} \min_{i \in [n]}\|x-x_i\| = \max_{\overline{x} \in \overline{X}_t} \max_{x \in B_t(\overline{x}) \cap \Omega} \min_{i \in [n]}\|x-x_i\| \\
	& \leq t + \max_{\overline{x} \in \overline{X}_t} \min_{i \in [n]}\|\overline{x}-x_i\|.
}
Define $E$ to be the following event :
$$E= \{(x_1,\dots,x_n) \in \Omega^n ~|~ \max_{j \in [m]} \min_{i \in [n]}\|\overline{x}_j-x_i\| < t \}.$$
The $n$ tuple $(x_1,..,x_n)$ belongs to $E$ if for each $\overline{x} \in \overline{X}_t$ there exists at least one $i \in [n]$ for which $\|\overline{x} - x_i\| < t$. $E$ can therefore be rewritten as follows :
$$E = \bigcap_{\overline{x} \in \overline{X}_t} \bigcup_{i \in [n]} \{(x_1,\dots,x_n) \in \Omega^n ~|~ \|\overline{x}-x_i\| < t\}.$$
In particular, note that 
$$E^c = \Omega^n \setminus E = \bigcup_{\overline{x} \in \overline{X}_t} \bigcap_{i \in [n]} \{(x_1,\dots,x_n) \in \Omega^n ~|~ \|\overline{x}-x_i\| \geq t\} = \bigcup_{\overline{x} \in \overline{X}_t}  ~(\Omega \setminus B_t(\overline{x}))^{ n}.$$
Applying a union bound, we get
\eqals{
\mathbb{P}(E^c) &= \mathbb{P}\big(\bigcup_{\overline{x} \in \overline{X}_t} (\Omega \setminus B_t(\overline{x}))^n \big) \\
& \leq \sum_{\overline{x} \in \overline{X}_t} \mathbb{P}\big((\Omega \setminus B_t(\overline{x}))^{n} \big) = \sum_{j \in [m]} U(\Omega \setminus B_t(\overline{x})))^n,
}
where the last step is due to the fact that $\mathbb{P}$ is a product measure and so $\mathbb{P}(A^n) = U^{\otimes n}(A^n) = U(A)^n$.
Now we need to evaluate $U(\Omega \setminus B_t(\overline{x})) = 1- U(B_t(\overline{x}))$ for $\overline{x} \in \overline{X_t}$. Since $\overline{X}_t \subset \Omega$, it holds
$$\forall \overline{x} \in \overline{X}_t,~U( B_t(\overline{x})) =  \tfrac{\vol(B_t(\overline{x}) \cap \Omega)}{\vol(\Omega)} \geq \tfrac{\min_{x \in \Omega} \vol(B_t(x) \cap \Omega)}{\vol(\Omega)}.$$

\noindent{\bf Step 3. Bounding $\vol(B_t(x) \cap \Omega)$ when $t \leq r$.} Let us now find a lower bound for $\min_{x \in \Omega} \vol(B_t(x) \cap \Omega)$. Recall that since $\Omega$ satisfies \cref{asm:geom-f:a}, $\Omega$ can be written $\Omega = \cup_{z \in S} B_r(z)$. .Let $t \leq r$, $x \in \Omega$. By the previous point, there exists $z \in S$ such that $x \in B_r(z) \subset \Omega$ and hence $B_t(x) \cap B_r(z) \subset B_t(x) \cap \Omega$. Let $C_{x,z,t}$ denote the cone centered in $x$ and directed to $z$ with aperture $\pi/3$. It is easy to see geometrically that $B_r(z) \cap B_t(x)$ contains the cone $C_{x,z,t}$ (this fact is proved in \cref{lm:cone}). Moreover, using the lower bound for the volume of this cone provided in \cref{lm:cone}, it holds:
\eqals{
\operatorname{vol}(\Omega \cap B_t(x)) &\geq  \operatorname{vol}(B_r(z) \cap B_t(x)) \geq  \operatorname{vol}(C_{x,z,t}) \\
&\geq \frac{2V_{d-1}}{\sqrt{3}d} \left(\tfrac{\sqrt{3}}{4}\right)^{d}t^d.
}

\noindent{\bf Step 4. Expressing $t$ with respect to $n$ and $\delta$ and guaranteeing that $t \leq r$.}
To conclude, let $C = \frac{V_{d-1}}{2d \vol(\Omega)} \left(\tfrac{\sqrt{3}}{4}\right)^{d-1}$. Since $N_t \leq (1 + 2R/t)^d$, and $(1 - c)^x \leq e^{-c x}$ for any $x \geq 0$ and $c \in [0,1]$, then
$$\mathbb{P}(E) \geq 1 - N_t \big(1- C t^d\big)^n \geq 1  - e^{-Ct^d n + d\log(1+ 2R/t)} \geq 1 - \delta,$$
where the last step is obtained by setting 
$$t = (Cn)^{-1/d} \left(\log\frac{(1+2R(Cn)^{1/d})^d}{\delta}\right)^{1/d}.$$
Then $h_{\widehat{X},\Omega} \leq 2t$ with probability at least $1-\delta$, when $t \leq r$. The desired result is obtained by further bounding $C$ and $t$ as follows. 

\noindent \textit{Bounding $C$.} It holds $\frac{2V_{d-1}}{\sqrt{3}d V_d} = \left(\tfrac{4}{3 d^2 \pi}\right)^{1/2} \frac{\Gamma(d/2 + 1)}{\Gamma(d/2 + 1/2)}$. Using Gautschi's inequality and the fact that $d \geq 1$, 
$$\left(\tfrac{2}{3 d \pi}\right)^{1/2}\leq \frac{2V_{d-1}}{\sqrt{3}d V_d}\leq \left(\tfrac{2(d+2)}{3 d^2 \pi}\right)^{1/2}\leq 1.$$ 
Since $\left(\tfrac{3d\pi}{2}\right)^{1/2d} \tfrac{4}{\sqrt{3}} \leq 2\sqrt{2\pi} $ for all $d \geq 1$, and since $V_d r^d \leq \vol(\Omega) \leq V_d R^d$, it holds 
$$(2\sqrt{2\pi} R)^{-d} \leq C \leq (4r/\sqrt{3})^{-d} \implies   \frac{n^{1/d}}{2\sqrt{2\pi}R} \leq (Cn)^{1/d} \leq \frac{\sqrt{3}n^{1/d}}{4 r} \leq \frac{n^{1/d}}{2 r} .$$

\noindent\textit{Bounding t.} Since,  $(1 + x)^d \leq (2x)^d $ for any $x \geq 1$ and $2R(Cn)^{1/d}  \leq \frac{R}{r} n^{1/d}$,  and $\frac{R}{r} n^{1/d} \geq 1$, it holds 
$$t \leq 2\sqrt{2\pi} R n^{-1/d} (\log \tfrac{n}{\delta} + d \log\tfrac{2R}{r})^{1/d} .$$

\noindent \textit{Guaranteeing $t \leq r$.} Applying \cref{lm:technical} to $\alpha = (2\pi)^{d/2} (2R/r)^d$ and $\beta = (2R/r)^d/\delta$, it holds that if 
$$n\geq 2\alpha \log(2\alpha\beta) =  2~(2\pi)^{d/2} (2R/r)^d ~\left(\log\frac{2}{\delta} + d/2\log (2\pi) + 2d \log(2R/r)\right),$$ then $\alpha/n \log(\beta n) \leq 1$, so
\[t \leq 2\sqrt{2\pi} R n^{-1/d} (\log \tfrac{n}{\delta} + d \log\tfrac{2R}{r})^{1/d} \leq r (\alpha/n \log(\beta n)^{1/d}  \leq r .\]

\epr

\subsection{Proof of \cref{thm:alg-glm}} \label{proof:thm:alg-glm}

\bpr 

Recall that $s > d/2$ and $m < s-\tfrac{d}{2}$ is a positive integer. Assume that $\Omega$ satisfies \cref{asm:geom-f:a} for a certain $r$ and that the diameter of $\Omega$ is bounded by $2R$. In particular, if $\Omega$ is a ball of radius $R$, then $\Omega$ satisfies \cref{asm:geom-f:a} with $r = R$. In the first step of the proof we guarantee that $n$ is large enough to apply \cref{thm:fill-distance-random-points} and that $h_{\widehat{X},\Omega}$, controlled by \cref{thm:fill-distance-random-points}, satisfies the assumptions of \cref{thm:bound-prob-sampled}. Then we apply \cref{thm:bound-prob-sampled}.

\noindent{\bf Step 1. Guaranteeing $n$ large enough and $h_{\widehat{X},\Omega} \leq r/(18(m-1)^2)$.} Applying \cref{lm:technical} to $\alpha = \left(\frac{2R}{r}\right)^d\max(3,10(m-1))^{2d}$ and $\beta = \tfrac{(2R)^d}{r^d~\delta}$, it holds that if 
$$n \geq 2\alpha \log(2\alpha \beta) = \left(\frac{2R}{r}\right)^d~\max(3,10(m-1))^{2d}\left( 2 \log \frac{2}{\delta} + 4d \log \left(\tfrac{R}{r} \max(6,20(m-1))\right) \right),$$
then $\alpha/n \, \log(\beta n) \leq 1$, which implies
\[n^{-1/d}(\log \frac{n}{\delta} + d\log \beta)^{1/d} \leq \frac{r}{2R\max(3,10(m-1))^2}.\] 
In particular, $n$ satisfying the condition above is large enough to satisfy the requirement of \cref{thm:fill-distance-random-points} (since $r \leq R$). Therefore, by applying \cref{thm:fill-distance-random-points} we have that with probability at least $1-\delta$,
\[h_{\widehat{X},\Omega} ~~\leq~~ 11 R \, n^{-\frac{1}{d}} \, (\log \tfrac{(2R)^d~n}{r^d~\delta} )^{{1}/{d}} \leq \frac{r}{\max(1,18(m-1)^2)}.\]

\noindent{\bf Step 2. Applying \cref{thm:bound-prob-sampled}.} In the previous step we provided a condition on $n$ such that $h_{\widehat{X},\Omega}$ satisfies $h_{\widehat{X},\Omega} \leq \frac{r}{\max(1,18(m-1)^2)}$. By \cref{prop:sobolev-kernel}, \cref{asm:H-rich} holds for the Sobolev kernel with smoothness $s$, for any $m \in \N$ since $m < s - d/2$. 
Then the conditions to apply \cref{thm:bound-prob-sampled} are satisfied. Applying \cref{thm:bound-prob-sampled}  with $\la\geq  2 \eta \max(1, \Cmul\Cdiff)$ and $\eta = \frac{3\max(1,18(m-1)^2)^m~d^m}{m!}  h^m_{\widehat{X},\Omega}$, we have
\eqals{
|\hat{c} - f_*| \leq 2\eta|f|_{\Omega,m} + \la\tr(A_*) \leq 3 \la (|f|_{\Omega,m}+\tr(A_*)),
}
Thus, under this condition, we have with probability at least $1-\delta$,
\[|\hat{c} - f_*| \leq C_{m,s,d} R^m n^{-m/d} (\log \frac{2^d n}{\delta}),\]
where
\[C_{m,s,d} = 6\times 11^m \times \frac{\max(1,18(m-1)^2)^m d^m}{m!} \max(1,\Cmul \Cdiff).\]
\noindent{\bf Step 3. Bounding the constant term $C_{m,s,d}$ in terms of $m, s, d$.} Note that
$$\frac{\Gamma(m+d/2)}{\Gamma(d/2)} = (d/2)...(d/2 + m-1) \leq (d/2 + m-1)^{m-1} $$
and 
$$\frac{\Gamma(s-d/2-m)}{\Gamma(s-d/2)} =\frac{1}{ (s-d/2 - m)....(s-d/2-1)} \leq \left(\frac{1}{s-d/2-m}\right)^{m-1},$$ 
which yields:
$$ \Cdiff \leq (2\pi)^{d/4} \left(\frac{d/2 + m-1}{s-d/2-m}\right)^{(m-1)/2}.$$
Moreover, using the bound on $\Cmul$, we get 
$$ \Cdiff \Cmul \leq 2^{s+1/2} ~ (2\pi)^{3d/4} \left(\frac{d/2 + m-1}{s-d/2-m}\right)^{(m-1)/2} . $$ 
This yields the following bound for $C_{m,s,d}$:
\[C_{m,s,d}\leq   \frac{6 \max(1,18(m-1)^2)^m (11d)^m}{m!} \max \left(1,2^{s+1/2} ~ (2\pi)^{3d/4} \left(\frac{d/2 + m-1}{s-d/2-m}\right)^{(m-1)/2} \right). \]
\epr

\section{Global minimizer. Proofs.}
\subsection{Proof of \cref{rm:existence-beta}}\label{proof:rm:existence-beta}

\bpr 
Since $f$ satisfies both \cref{asm:geom-f:b,asm:unique-minimizer}, denote by $\zeta$ the unique minimizer of $f$ in $\Omega$. Since $\zeta$ is a strict minimum by \cref{asm:geom-f:b}, there exists $\beta_1 > 0$ such that $\nabla^2 f(\zeta) \succeq \beta_1 I$. Thus, since $f \in C^2(\R^d)$, there exists a small radius $t>0$ such that $\nabla^2 f(x) \succeq \tfrac{\beta_1}{2} I$ for all $x \in B_t(\zeta)$ and hence
\eqal{\label{eq:beta1}\forall x \in \Omega \cap B_t(\zeta),~ f(x) - f_* = f(x) -f(\zeta) - \nabla f(\zeta) \geq \tfrac{\beta_1}{4}\|x-\zeta\|^2.}
\noindent Moreover, since $f$ has no minimizer on the boundary of $\Omega$ and since $\zeta$ is the unique minimizer of $f$ on $\Omega$, $f$ has no minimizer on $K = \overline{\Omega} \setminus{B_t(x)}$ which is a compact set. Denote by $m$ the minimum of $f$ on $K$. Since $K$ is compact, this minimum is reached and since $f$ does not reach its global minimum $f_*$ on $K$, we have $m - f_* > 0$. Let $R$ be a radius such that $\overline{\Omega} \subset B_R(\zeta)$, which exists since $\Omega$ is bounded. Then, since for any $x \in \overline{\Omega}$, $\|x - \zeta\| < R$, it holds for any $x \in K$ : 
\eqal{\label{eq:beta2}
f(x) - f_* = f(x) - m + m-f_* \geq m-f_* = \frac{2(m-f_*)}{2R^2}R^2 \geq  \frac{2(m-f_*)}{2R^2}\|x-\zeta\|^2.
}
Thus, taking $\beta = \min(\tfrac{\beta_1}{2},\tfrac{2(m-f_*)}{R^2})$ and combining \cref{eq:beta1,eq:beta2}, it holds
\[\forall x \in  \Omega,~f(x)-f_* \geq \frac{\beta}{2}\|x-\zeta\|^2.\]

\epr

\subsection{Proof of \cref{thm:Hilbert-z-isgood}} \label{proof:thm:Hilbert-z-isgood}

\bpr

Let us divide the proof into four steps.

\noindent{\bf Step 1: Extending the parabola outside of $\Omega$} Since $\Omega$ is an open set containing $\zeta$, there exists $t > 0$ such that the ball $B_t(\zeta) \subset \Omega$. Define $\delta = \tfrac{\beta - \nu}{2}t^2$. It holds :
\eqal{\label{eq:deltabeta} \forall x \in \R^d\setminus{\Omega},~ \frac{\beta}{2}\|x-\zeta\|^2 \geq \frac{\nu}{2}\|x-\zeta\|^2 + \delta.}
\noindent Now define the following open set : 
$$ \widetilde{\Omega} = \left\{x \in \R^d ~:~  f(x) - f_*  - \tfrac{\beta}{2}\|x - \zeta\|^2>  - \delta/2\right\} .$$
It is open since $f$ is continuous. Moreover, it contains the closure of $\Omega$ which we denote with $\overline{\Omega}$ which is compact since it is closed and bounded in $\R^d$. Theorem 1.4.2 in \cite{hormander2015analysis} applied to $X = \widetilde{\Omega}$ and $K = \overline{\Omega}$ shows the existence of $\chi : \R^d \rightarrow \R$ such that $\chi \in C^\infty(\R^d)$, $\chi(x) \in [0,1]$, $\chi = 1$ on $\overline{\Omega}$ and $\chi = 0$ on $\R^d \setminus{\widetilde{\Omega}}$. Finally, define $\overline{p}_{\nu}(x) := \frac{\nu}{2}\|x-\zeta\|^2 \chi(x)$. $\overline{p}_{\nu}$ satisfies the following properties :

\begin{itemize}
    \item $\overline{p}_{\nu} \in C^\infty(\R^d)$;
    \item for all $x \in \overline{\Omega}$, $\overline{p}_{\nu}(x) = \frac{\nu}{2}\|x-\zeta\|^2 \leq \frac{\beta}{2}\|x-\zeta\|^2$;
    \item for all $x \in \R^d \setminus{\widetilde{\Omega}}$, $\overline{p}_{\nu}(x) =0$;
    \item for all $x \in \widetilde{\Omega}\setminus{\Omega}$, $f(x) - f_* - \overline{p}_{\nu}(x) \geq \delta/2$.  
\end{itemize}

\noindent The first, second and third properties are direct consequences of the properties of $\chi$ and the fact that $\nu < \beta$. The last property comes from combining \cref{eq:deltabeta} with the definition of $\widetilde{\Omega}$ and the fact that $\chi \in [0,1]$ :  
\eqals{\forall x \in \widetilde{\Omega}\setminus{\Omega},~f(x) - f_* - \overline{p}_{\nu}(x) &= f(x) - f_* - \chi(x)\tfrac{\nu}{2}\|x-\zeta\|^2 \\
&\geq f(x) - f_* - \tfrac{\nu}{2}\|x-\zeta\|^2\\
& = \left(f(x) - f_* - \tfrac{\beta}{2}\|x-\zeta\|^2\right) + \left(\tfrac{\beta}{2}\|x-\zeta\|^2 - \tfrac{\nu}{2}\|x-\zeta\|^2\right) \\
&\geq -\delta/2 + \delta = \delta/2.
}

\noindent{\bf Step 2: Extending $x 
\mapsto f(x) - \tfrac{\nu}{2}\|x-\zeta\|^2$ outside of $\Omega$.} Define $g(x) = f(x) - \overline{p}_\nu(x)$ on $\R^d$. Then $g$ satisfies \cref{asm:geom-f:b}, $g$ has exactly one minimizer in $\Omega$ which is $\zeta$, and its minimum is $g(\zeta) = f_*$. Indeed, the fact that $g \in C^2(\R^d)$ comes from the fact that $f \in C^2(\R^d)$ by \cref{asm:geom-f:b} on $f$ and the fact that $\overline{p}_\nu \in C^\infty(\R^d)$. Moreover, $g \geq f_*$ on $\R^d$ and $g-f_* \geq \delta/2$ on $\partial \Omega$. Indeed, first note that since $\nu< \beta$, it holds 
$$ \forall x \in \Omega,~g(x) = f(x) - \overline{p}_{\nu}(x) = f(x) - \tfrac{\nu}{2}\|x-\zeta\|^2 \geq f(x) - \tfrac{\beta}{2}\|x-\zeta\|^2 \geq f_*,$$ 
where the last inequality comes from \cref{eq:beta}. Second, since $\overline{p}_{\nu} = 0$ on $\R^d\setminus{\widetilde{\Omega}} $ and since $f_*$ is the minimum of $f$, for any $x \in \R^d \setminus{\widetilde{\Omega}}$, $g(x) - f_* = f(x) - f_* \geq 0$. Finally, by the last point of the previous step, we see that $g(x) \geq f_* + \delta/2 > f_*$ for any $x \in \widetilde{\Omega}\setminus{\Omega}$. In particular, $g(x) \geq f_* + \delta/2$ for any $x \in \partial \Omega$.
Since $g(\zeta) = f(\zeta) = f_*$, we see that $f_*$ is the minimum of $g$ on $\R^d$ and that this minimum is reached at $\zeta$ and is not reached on the boundary of $\Omega$. The fact that $\zeta$ is the unique minimum on $\Omega$ comes from the fact that since $\nu < \beta$ and by \cref{eq:beta} we have that for any $x \in \Omega\setminus{\{\zeta\}}$ the following holds
\eqal{
\begin{split}
\label{eq:useful} g(x) &= f(x) - \overline{p}_{\nu}(x) = f(x) - \tfrac{\nu}{2}\|x-\zeta\|^2 \\
&
> f(x) - \tfrac{\beta}{2}\|x-\zeta\|^2 \geq f_*.
\end{split}
}
The fact that this minimum is not reached on the boundary of $\Omega$ comes from the fact stated above that $g(x) \geq f_* + \delta/2$ for any $x \in \partial \Omega$. Finally, the fact that $\zeta$ is a strict minimum of $g$ also comes from \cref{eq:useful} which implies that $\nabla^2 g(\zeta) \succeq (\beta - \nu) I$ since $g$ reaches a minimum in $\zeta$, $g$ is $C^2$ and $\nu < \beta$. 

\noindent Note that $g$ also satisfies \cref{asm:analytic-f} since $f$ satisfies \cref{asm:analytic-f} and $\overline{p}_\mu \in C^\infty(\R^d) \subset C^2(\R^d)\cap \hh$ by \cref{asm:H-rich:a}.

\noindent{\bf Step 3: Applying \cref{cor:g-good-has-A} to $g$.} The previous point shows that $g$ satisfies \cref{asm:geom-f:b,asm:analytic-f}and that $g$ has a unique minimum in $\Omega$. Moreover, $\hh$ satisfies \cref{asm:H-rich}. Hence, \cref{cor:g-good-has-A} to $g$ and $\hh$, the following holds : there exists $A_* \in \pdm(\hh)$ with $\rank(A_*) \leq d + 1$ such that $g(x) - f^* = \scal{\phi(x)}{A_* \phi(x)}$ for all $x \in \Omega$.

\noindent{\bf Step 4.} Let $p_0$ be the maximum of \cref{eq:prob-conv-z}.
In \cref{lm:conv-z-is-ok} we have seen that the solution of \cref{eq:prob-conv-z} is $p_0 = f_*$.
Since $A \succeq 0$ implies $\scal{\phi(x)}{ A \phi(x)} \geq 0$ for all $x \in \Omega$, the problem in \cref{eq:prob-conv-z} is a relaxation of \cref{eq:prob-intermediate-z}, where the constraint $f(x) - \tfrac{\nu}{2}\|x\|^2 + \nu x^\top z - c = \scal{\phi(x)}{ A \phi(x)}$ is substituted by $f(x) - \tfrac{\nu}{2}\|x\|^2 + \nu x^\top z - c \geq 0, \forall x \in \Omega$.
Then $p_0 \geq p^*$ if a maximum $p^*$ exists for \cref{eq:prob-intermediate-z}. Thus, if there exists $A$ that satisfies the constraints in \cref{eq:prob-intermediate-z} for the value $c_* = f_*+\frac{\nu}{2}\|\zeta\|^2$ and $z_* = \zeta$, then $p_0 = p^*$ and $(c_*, \zeta, A)$ is a minimizer of \cref{eq:prob-intermediate-z}.

\noindent The proof is concluded by noting that indeed there exists $A$ that satisfies the constraints in \cref{eq:prob-intermediate-z} for the value $c_* =  f_*+\frac{\nu}{2}\|\zeta\|^2$ and $z_* = \zeta$ and it is obtained by the previous step.

\epr

\subsection{Proof of \cref{thm:bound-prob-sampled-z}} \label{proof:thm:bound-prob-sampled-z}

\bpr
The proof is a variation of the the one for \cref{thm:bound-prob-sampled}, the main difference is that we take care of the additional term $z - \zeta$.

\noindent{\bf Step 0. The SDP problem in \cref{eq:prob-relax-z} admits a solution} 

\noindent(a) Under the constraints of \cref{eq:prob-relax-z}, $c - \frac{\nu}{2}\|z\|^2$ cannot be larger than $\min_{i \in [n]} f(x_i)$. Indeed, for any $i\in [n]$, since $B \succeq 0$, the $i$-th constraint implies
$$ f(x_i) - \tfrac{\nu}{2}\|x_i - z\|^2 - c + \frac{\nu}{2}\|z\|^2 =  f(x_i) - \tfrac{\nu}{2}\|x_i\|^2 + \nu  x_i^\top z - c =  \Phi_i B \Phi_i \geq 0. $$
Hence, $f(x_i) \geq f(x_i) - \tfrac{\nu}{2}\|x_i - z\|^2 \geq c + \frac{\nu}{2}\|z\|^2$. Thus, since $B \succeq 0$, for any $B,z,c$ satisfying the constraint, $c - \frac{\nu}{2}\|z\|^2 - \lambda \tr(B) \leq \max_{i \in [i]}{f(x_i)}$. 

\noindent (b) There exists an admissible point. Indeed let $(c_*, z_*, A_*)$ be the solution of \cref{eq:prob-intermediate-z} such that $A_*$ has minimum trace norm (by \cref{thm:Hilbert-z-isgood}, we know that this solution exists with $c_* = f_*$ and $z_* = \zeta$, under \cref{asm:geom-f,asm:H-rich,asm:analytic-f,asm:unique-minimizer}). Then, by \cref{lm:exists-B} applied to $g(x) = f(x) - \frac{\nu}{2}\|x\|^2 - \nu x^\top z_* - c_*$ and $A =A_*$, given $\widehat{X} = \{x_1,\dots,x_n\}$ we know that there exists $\overline{B} \in \pdm(\R^n)$ satisfying $\tr(\overline{B}) \leq \tr(A_*)$ s.t. the constraints of \cref{eq:prob-relax-z} are satisfied for $c = c_*$ and $z = z_*$. Then $(c_*,z_*, \overline{B})$ is admissible for the problem in \cref{eq:prob-relax-z}. Since there exists an admissible point for the constraints of \cref{eq:prob-relax-z} and its functional cannot be larger than $\max_{i \in [n]}f(x_i)$, then the SDP problem in \cref{eq:prob-relax-z} admits a solution~\cite{boyd2004convex}.

\noindent{\bf Step 1. Consequences of existence of $A_*$.} 
Let $(\hat{c}, \hat{z}, \hat{B})$ one minimizer of \cref{eq:prob-relax-z}. The existence of the admissible point $(c_*, z_*, \overline{B})$ implies that
\eqal{\label{eq:lb-cc-0-z}
\hat{c} - \tfrac{\nu}{2}\|\hat{z}\|^2 - \la \tr(\hat{B}) \geq c_* - \tfrac{\nu}{2}\|z_*\|^2 - \la\tr(\overline{B}) \geq f_* - \la\tr(A_*).
}
From which we derive,
\eqal{\label{eq:lb-cc-z}
\la\tr(\hat{B}) - \la\tr(A_*) ~~\leq~~ \Delta, \quad \Delta := \hat{c} - \tfrac{\nu}{2}\|\hat{z}\|^2 - f_*.
}

\noindent{\bf Step 2. $L^\infty$ bound due to the scattered zeros.}
Note that the solution $(\hat{c},\hat{z},\hat{B})$ satisfies $\hat{g}(x_i) = \Phi_i^\top \hat{B} \Phi_i$ for $i \in [n]$, where the function $\hat{g}$ is defined as $\hat{g}(x) = f(x) - \frac{\nu}{2}\|x\|^2 + \nu x^\top \hat{z} - \hat{c}$ for $x \in \Omega$, moreover $h_{\widehat{X},\Omega} \leq\frac{r}{\max(1,18(m-1)^2)}=  \frac{r}{18(m-1)^2}$ by assumption, since $m\geq 2$. Then we can apply \cref{thm:inequality-scattered-data} with $g = \hat{g}$, $\tau = 0$ and $B = \hat{B}$ obtaining for all $x \in \Omega$
$$ f(x) - \tfrac{\nu}{2}\|x\|^2 + \nu x^\top \hat{z} - \hat{c} = \hat{g}(x) \geq - \eta (|\hat{g}|_{\Omega,m} + \Cmul\Cdiff\tr(\hat{B})), \quad \eta = C_0 h^m_{\widehat{X},\Omega},$$ 
where $C_0$ is defined in \cref{thm:inequality-scattered-data} and $C_0 = 3\frac{(18d)^m (m-1)^{2m}}{m!}$ since $m \geq 2$. Since the inequality above holds for any $x \in \Omega$, by evaluating it in the global minimizer $\zeta \in \Omega$, we have $f(\zeta) = f_*$ and so
$$-\Delta - \tfrac{\nu}{2}\|\hat{z} - \zeta\|^2 = \hat{g}(\zeta) \geq -\eta (|\hat{g}|_{\Omega,m} + \Cmul\Cdiff\tr(\hat{B})).$$
Now we bound $|\hat{g}|_{\Omega,m}$. Since $\hat{g}(x) = f(x) - p_{\hat{z},\hat{c}}(x)$, where $p_{\hat{z},\hat{c}}$ is a second degree polynomials defined as $p_{\hat{z},\hat{c}} = \tfrac{\nu}{2} \|x\|^2 -  \nu x^\top \hat{z} + \hat{c}$, we have
\eqal{\label{eq:bound-hatg-z}
|\hat{g}|_{\Omega,m} \leq |f|_{\Omega,m} \,+\, |p_{\hat{z},\hat{c}}|_{\Omega,m} \leq |f|_{\Omega,m} + \nu,
}
since for $m = 2$, we have $|p_{\hat{z},\hat{c}}|_{\Omega,2} = \sup_{i,j \in [d], x \in \Omega} |\frac{\partial^2 p_{\hat{z},\hat{c}}(x)}{\partial x_i \partial x_j}| = \nu$ and also $|p_{\hat{z},\hat{c}}|_{\Omega,m} = 0$ for $m > 2$. Then 
\eqal{\label{eq:key-inequality-z}
\Delta \leq \Delta + \tfrac{\nu}{2}\|\hat{z} - \zeta\|^2 \leq \eta |f|_{\Omega,m}  + \eta \Cmul \Cdiff \tr(\hat{B}) + \eta \nu.
}
\noindent{\bf Conclusion.}
Combining \cref{eq:key-inequality-z} with \cref{eq:lb-cc-z}, since $\tfrac{\nu}{2}\|\hat{z} - \zeta\|^2 \geq 0$ and since $\la \geq 2\Cmul \Cdiff \eta$ by assumption, we have
$$\tfrac{\la}{2} \tr(\hat{B}) \leq (\la - \Cmul \Cdiff \eta) \tr(\hat{B}) \leq  \eta |f|_{\Omega,m} + \eta \nu + \la\tr(A_*),$$
from which we obtain \cref{eq:appr-A-z}. Moreover, the inequality \cref{eq:appr-c-z} is derived by bounding $\Delta$ from below as $\Delta \geq -\la \tr(A_*)$ by \cref{eq:lb-cc-z}, since $\tr(\hat{B}) \geq 0$ by construction, and bounding it from above as
$$\Delta  \leq 2\eta |f|_{\Omega,m} + 2\eta \nu + \la\tr(A_*),$$
that is obtained by combining \cref{eq:key-inequality-z} with \cref{eq:appr-A-z} and with the assumption $\Cmul \Cdiff \eta \leq \la/2$. Finally from \cref{eq:key-inequality-z} we obtain
\eqals{
\tfrac{\nu}{2}\|\hat{z} - \zeta\|^2 \leq |\Delta| ~+~ \eta |f|_{\Omega,m}  + \eta \Cmul \Cdiff \tr(\hat{B}) + \eta \nu,
}
from which we derive the bound $\frac{\nu}{2}\|\hat{z} - \zeta\|^2$ in \cref{eq:appr-z}, by bounding $|\Delta|$ and $\tr(\hat{B})$ via \cref{eq:appr-c-z} and \cref{eq:appr-A-z}.
\epr

\section{Proofs for the extensions}

\subsection{Proof of \cref{thm:appr-solution}}\label{proof:thm:appr-solution}

\bpr
Let $(\hat{c}, \hat{B})$ be a minimum trace-norm solution of \cref{eq:prob-relax}. The minimum $p_{\la,n}$ of \cref{eq:prob-relax} then corresponds to $p_{\la,n} = \hat{c} - \la \tr(\hat{B})$. Combining \cref{eq:appr-functional} with \cref{eq:lb-cc} from the proof of \cref{thm:bound-prob-sampled} and the fact that $\theta_2 \leq \la/8$, we have that
\eqal{\label{eq:lb-cc-appr}
\tfrac{7}{8}\la\tr(\tilde{B}) - \la\tr(A_*) - \theta_1 \leq \tilde\Delta, \quad \tilde\Delta := \tilde{c} - f_*.
}
Analogously to Step 3 of the proof of \cref{thm:bound-prob-sampled}, by applying \cref{thm:inequality-scattered-data} to \cref{eq:appr-constraints} with $g(x) = f(x) - \tilde{c}, B = \tilde{B}$ and $\tau = \tau_1 + \tau_2\tr(\tilde{B})$, we obtain for any $x \in \Omega$
$$f(x) - \tilde{c} ~\geq~ -2\tau_1 - 2\tau_2\tr(\tilde{B})  ~-~ \eta(|g|_{\Omega,m} + \Cmul\Cdiff \tr(\tilde{B})),  \qquad \eta = C_0 h^m_{\widehat{X},\Omega},$$
with $C_0$ defined in \cref{thm:inequality-scattered-data}. Now evaluating the inequality above for $x = \zeta$, noting that $|g|_{\Omega,m} = |f|_{\Omega,m}$ since $m \geq 1$, and considering that by assumption $\tau_2 \leq \la/8$ and $\Cmul\Cdiff\eta \leq \la/2$ we derive
\eqal{\label{eq:key-inequality-approx}
\tilde{\Delta} = -(f(\zeta) - \tilde{c}) \leq 2\tau_1 + \tfrac{3}{4}\la\tr(\tilde{B}) + \eta|f|_{\Omega,m}.
}
The desired result is obtained by combining \cref{eq:key-inequality-approx} and \cref{eq:lb-cc-appr} as we did in Step 3 of \cref{thm:bound-prob-sampled}.
\epr

\subsection{Proof of \cref{cor:g-in-Cs-has-A}}\label{proof:cor:g-in-Cs-has-A}

\bpr 
\noindent Define $\hh = \{g \in C^s(\Omega)~:~ \exists f \in C^s(\R^d),~f|_\Omega = g \}$, endowed with the following norm : 
$$ \forall g \in \hh,~ \|g\|_{\hh} = \sup_{|\alpha| \leq s}\sup_{x \in \Omega}{\|\partial^\alpha g(x)\|}.$$
Note that this norm is well defined since for any $g \in \hh$, since there exists $f \in \C^s(\R^d)$ such that $g = f|_\Omega$, since all the derivatives of $f$ are continuous hence bounded on $\Omega$ which is bounded, so are all the derivatives of $g$. 

\noindent Now note that $\hh$ satisfies \cref{asm:H-rich:a,asm:H-rich:b,asm:H-rich:c}. Indeed, given $u,v \in \hh$ the first assumption is satisfied as a simple consequence of the Leibniz formula, since for any $x \in \Omega$,  $\partial^{\alpha}(u\cdot v)(x) = \sum_{\beta \leq \alpha}{ \binom{\alpha}{\beta}\partial^\beta u(x) \partial^{\alpha - \beta} v(x)}$
which in turn implies that for any $|\alpha| \leq s$ and $x \in \Omega$, $\|\partial^{\alpha}(u\cdot v)(x)\| \leq 2^{|\alpha|}~\|u\|_{\hh} ~\|v\|_{\hh}$ and hence $\|u\cdot v\|_{\hh} \leq 2^s \|u\|_{\hh} ~\|v\|_{\hh} $.  \cref{asm:H-rich:b} is trivially satisfied and \cref{asm:H-rich:c} is a simple consequence of the dominated convergence theorem. Indeed, if $u \in \hh$ and $\overline{u} \in C^s(\R^d)$ such that $\overline{u}|_\Omega = u$, define

$$ \forall x,z \in \R^d,~ \overline{v}_z(x) = \int_{0}^1{\overline{u}(z + t(x-z))dt}.$$
$\overline{v}_z$ is in $C^s(\R^d)$ by dominated convergence, and $v_{z} = \overline{v}|_{\Omega}$ satisfies the desired property (in this case, there is no need to depend on $r$ and one can simply take $g_{r,z} = v_z$).

\noindent Moreover, if $f \in C^{s+2}(\R^d)$, then in particular, for any $i,j \in [d]$, $\tfrac{\partial f}{\partial x_i \partial x_j} \in C^s(\R^d) $ and hence its restriction to $\Omega$ is in $\hh$. Moreover, in that case, it is obvious that since $s \geq 0$, $f|_\Omega \in \hh $. This shows that $f$ satisfies \cref{asm:geom-f:b,asm:analytic-f}.

\noindent Therefore, \cref{thm:g-good-has-A} can be applied, and there exist $\tilde{w}_1,\dots,\tilde{w}_p \hh$, $p \in \Np$, such that 
$$ \forall x \in \Omega,~f(x) - f_* = \sum_{j \in [p]} w_j^2(x).$$
By definition of $\hh$, taking $w_1,...,w_p$ such that $w_j|_\Omega = \tilde{w}_j$, the corollary holds.

\epr 

\subsection{Proof of \cref{thm:Aeps-exists}}\label{proof:thm:Aeps-exists}
\bpr
In this proof we will use the results recalled in \cref{sec:sobolev-spaces} about Sobolev spaces. By \cref{cor:g-in-Cs-has-A} we have that there exists $\overline{w}_1,\dots,\overline{w}_p \in C^s(\R^d)$ such that $f(x) = \sum_{j \in [p]} \overline{w}_j^2(x)$ for any $x \in \Omega$. Define $w_j = \overline{w}_j|_{\Omega}$. Note that by \cref{prop:embedding-sobolev}, $w_j \in W^s_\infty(\Omega)$ for $j \in [p]$. 
Now let $\eps \in (0,1]$, for any $j \in [p]$, let $w^\eps_j \in C^\infty(\Omega)$ be the $\eps$ approximation of $w_j$ as defined in \cref{prop:approximation-sobolev}, i.e., $w_j = \tilde{w}_j^\epsilon|_{\Omega}$ where $\tilde{w}_j^\epsilon \in C^{\infty}(\R^d)$ and 
\eqal{\label{eq:weps}\|w_j - w^\eps_j\|_{L^\infty(\Omega)} \leq C_1 \eps^s \|w_j\|_{W^s_\infty(\Omega)}, \qquad \|w^\eps_j\|_{W^r_\infty(\Omega)} \leq C_2 \eps^{s-r}\|w_j\|_{W^s_\infty(\Omega)}.
}
with $C_1, C_2$ depending only on $r,s,d,\Omega$.
Now, since $f(x) - f_* = \sum_{j \in [p]} w_j^2(x)$, for any $x \in \Omega$, we have
\eqals{
\big\|f - f_* - \sum_{j \in [p]} {w^\eps_j}^2 & \big\|_{L^\infty(\Omega)}  = \big\|\sum_{j \in [p]} (w_j - w^\eps_j)(2w_j ~-~  (w_j -w^\eps_j))\big\|_{L^\infty(\Omega)} \\
& \leq \sum_{j \in [p]} \|w_j -  w^\eps_j\|_{L^\infty(\Omega)}(2\|w_j\|_{L^\infty(\Omega)} + \|w_j -  w^\eps_j\|_{L^\infty(\Omega)})\\
& \leq \sum_{j \in [p]} C_1 \eps^s \|w_j\|_{W^s_\infty(\Omega)} (2\|w_j\|_{W^s_\infty(\Omega)} + C_1 \eps^s \|w_j\|_{W^s_\infty(\Omega)}) \\
& \leq \eps^s ~ p~C_1(2+C_1) \max_{j \in [p]} \|w_j\|^2_{W^s_\infty(\Omega)},
}
where we use the first equation of \cref{eq:weps} to go from the second to the third line.

\noindent Recall that $\hh$ is defined to be a RKHS associated to the Sobolev kernel $k_r$ defined in \cref{ex:sobolev-kernel} for a given $r > \max(s,d/2)$. As mentioned in \cref{ex:sobolev-kernel}, we have $\hh = W_2^r(\Omega)$, and $\|\cdot\|_\hh$ is equivalent to $\|\cdot\|_{W^r_2(\Omega)}$, i.e., there exists $C_4$ depending on $\Omega, r, d$ such that $\frac{1}{C_4}\|\cdot\|_{W^r_2(\Omega)} \leq\|\cdot\|_\hh \leq C_4 \|\cdot\|_{W^r_2(\Omega)}$.

\noindent Since $w^\eps_j \in W^r_2(\Omega) = \hh$ by \cref{eq:weps} for all $j \in [p]$, we can define : 
$$A_\eps = \sum_{j \in [p]} w^\eps_j \otimes_{\hh} w^\eps_j.$$
It holds : 
\eqals{
\tr(A_\eps) &\leq p \max_{j \in [p]} \|w^\eps_j\|^2_\hh \leq p C_4^2 \max_{j \in [p]}\|w^\eps_j\|^2_{W^r_2(\Omega)}  \\
&\leq p C_4^2 C_5^2  \max_{j \in [p]}\|w^\eps_j\|^2_{W^r_\infty(\Omega)} \\
&\leq \eps^{2(s-r)} ~ p (C_2C_4 C_5)^2 \max_{j \in [p]}\|w_j\|^2_{W^s_\infty(\Omega)}.
}
where we used \cref{eq:weps} and the fact that there exists $C_5$ such that $\|\cdot\|_{W^r_\infty(\Omega)} \leq C_5 \|\cdot\|_{W^r_2(\Omega)}$ (see  \cref{prop:embedding-sobolev}).
To conclude, we use \cref{prop:equivnorm} to bound $\|\cdot\|_{W^s_\infty(\Omega)} \leq C_6 \|\cdot\|_{\Omega,s}$. 
\epr

\subsection{Proof of \cref{thm:global-min-low-smoothness}}\label{proof:thm:global-min-low-smoothness}
\bpr
The proof of the existence of a minimizer corresponds essentially to the first part of the proof of \cref{thm:bound-prob-sampled} and we skip it. Let $\eps \in (0,1]$, by applying \cref{thm:Aeps-exists} to $f$ we know that there exists $A_\eps \in \pdm(\hh)$ satisfying \cref{eq:Aeps-approx}. Define $f_\eps =  \scal{\phi(x)}{A_\eps \phi(x)}$ for all $x \in \Omega$, by \cref{thm:Aeps-exists} we have
\eqals{
\tr(A_\eps) ~~\leq~~ C_f\epsilon^{-2(r-s)}, \qquad \sup_{x \in \Omega} |f(x) - f_\eps(x)| ~~\leq~~ C_f'\eps^{s}.
}
Now consider the problem in \cref{eq:prob-relax-eps} and denote by $p^\eps_{\la,n}$ its optimum.
Since $f_\eps(x_i) - c = \Phi_i^\top B \Phi_i$ implies $|f(x_i) - c - \Phi_i^\top B \Phi_i| \leq \tau$, since we required $\tau \geq \sup_{x \in \Omega}|f(x) - f_\eps(x)|$. Then in this case \cref{eq:prob-relax-tau} is a relaxation of \cref{eq:prob-relax-eps} and we have that $p^\eps_{\la, n} - \tilde{c} - \la \tr(\tilde{B}) \leq 0$. So, we can apply \cref{thm:appr-solution} to $f_\eps$ with $\theta_1,\theta_2, \tau_2 = 0$ and $\tau_1 = \tau$, obtaining for any $m \in \N$ and $m < s-d/2$
\eqals{|\tilde{c} - f^\eps_*| &~\leq~ 14 \tau + 7\eta |f_\epsilon|_{\Omega,m} + 6 \lambda \tr(A_\epsilon) ,\\
\tr(\tilde{B}) &~\leq~ 8\,\tr(A_{\epsilon}) ~+~ 8\tfrac{\eta}{\la}\,|f_\epsilon|_{\Omega,m} ~+~ 16 \tfrac{\tau}{\lambda}.}

where $f^\eps_*$ is the infimum of $f_\eps$ (see \cref{rem:sufficiency-asm-f}), and satisfies
$$|f_* - f^\eps_*| = |\min_{x \in \Omega} f(x) - \inf_{x \in \Omega} f_\eps(x)| = |\inf_{x \in \Omega} f(x) - \inf_{x \in \Omega} f_\eps(x)| \leq \sup_{x \in \Omega} |f(x) - f_\eps(x)| \leq \tau.$$
By the same reasoning in the proof of \cref{thm:inequality-scattered-data} used to bound $|g|_{\Omega,m}$, we have that
$$|f_\eps|_{\Omega,m} \leq \Cmul \Cdiff \tr(A_\eps) \leq C_f \Cmul \Cdiff \eps^{-2r+2s}.$$
Combining together the inequalities above, with the fact that $\la \geq 2 \Cmul \Cdiff \eta$, we have
\eqals{
|\tilde{c} - f_*| ~\leq~ 10 \la\,C_f\epsilon^{-2(r-s)} ~+~ 15 \tau, \qquad \tr(\tilde{B}) &~\leq~ 12\,C_f\epsilon^{-2(r-s)} ~+~ 16 \tfrac{\tau}{\lambda}.
}
Now we set $\eps$ as large as possible such that $\tau \geq \sup_{x \in \Omega}|f(x) - f_\eps(x)|$ holds. In particular we know that requiring $\tau \geq C'_f \eps^s$ guarantees $\tau \geq \sup_{x \in \Omega}|f(x) - f_\eps(x)|$. Then by setting $\eps = 1$ when $\tau \geq C_f'$, we have
\eqals{
|\tilde{c} - f_*| ~\leq~ 10\, \la\,C_f ~+~ 15\tau \, , \qquad \tr(\tilde{B}) &~\leq~ 12\,C_f ~+~ 16 \tfrac{\tau}{\lambda}.
}
By setting $\eps = (\tau/C_f')^{1/s}$ when $\tau \leq C_f'$, we have
\eqals{
|\tilde{c} - f_*| ~\leq~ 10\, \la\,C_f(C_f')^{2\frac{r-s}{s}}\tau^{-2\frac{r-s}{s}} ~+~ 15\tau \,, \qquad \tr(\tilde{B}) &~\leq~ 12\,C_f( C_f')^{2\frac{r-s}{s}}\tau^{-2\frac{r-s}{s}} ~+~ 16 \tfrac{\tau}{\lambda}.
} 
Selecting $\tau = \la^{\tfrac{s}{2r-s}}$ and combining the inequality for the two cases above, leads to
\eqals{
|\tilde{c} - f_*| ~\leq~ \widetilde{C}_f (\la + \la^{\tfrac{s}{2r-s}}), ~~~\tr(\tilde{B}) &~\leq~ 12\,C_f ~+~ \widetilde{C}^\prime_f \la^{-(1-\frac{s}{2r-s})}.
}
where 
$$ \widetilde{C}_f = \max\left(10C_f (C_f')^{2\frac{r-s}{s}} + 15,10C_f\right),\qquad  \widetilde{C}^\prime_f = 12\,C_f( C_f')^{2\frac{r-s}{s}} + 16$$
\epr

\subsection{Certificate of optimality for the global minimizer candidate of \cref{eq:prob-relax-z}}\label{sec:certificate-prob-relax-z}
\bt[Certificate of optimality for \cref{eq:prob-relax-z}]\label{thm:certificate-prob-relax-z}
Let $\Omega$ satisfy \cref{asm:geom-f:a} for some $r > 0$. Let $k$ be a kernel satisfying \cref{asm:H-rich:a,asm:H-rich:d} for some $m \geq 2$.
Let $\widehat{X} = \{x_1,\dots, x_n\} \subset \Omega$ with $n \in \N$ such that $h_{\widehat{X},\Omega} \leq \frac{r}{18(m-1)^2}$.
Let $f \in C^m(\Omega)$ and let $\hat{c} \in \R, \hat{z} \in \R^d,\hat{B} \in \pdm(\R^n)$ and $\tau \geq 0$ satisfying
\eqal{\label{eq:empirical-inequality-z}
|f(x_i) - \tfrac{\nu}{2}\|x_i\|^2 + \nu x_i^\top \hat z - \hat{c} ~-~ \Phi_i^\top \hat{B} \Phi_i | \leq \tau, \quad i \in [n]
}
where $\Phi_i$ are defined in \cref{sec:sketch}. Let $f_* = \min_{x \in \Omega} f(x)$ and $\hat{f} = \hat{c} - \tfrac{\nu}{2}\|\hat{z}\|^2$. Then,
\eqal{
\label{eq:cert-2:1}
|f(\hat{z}) -f_*| &\leq f(\hat{z}) - \hat{f} +2 \tau + C_1 h_{\widehat{X},\Omega}^m ,  \\
\label{eq:cert-2:2}\tfrac{\nu}{2}\|\zeta -\hat{z}\|^2 & \leq f(\hat{z}) - \hat{f}  + 2\tau +  C_2 h_{\widehat{X},\Omega}^m.
}
and $C_1 = C_0(|f|_{\Omega,m}+\Cmul\Cdiff\tr(\hat{B}) + \Cmul\Cdiff\hat{C})$, $C_2 = C_0(|f|_{\Omega,m}+ \nu+\Cmul\Cdiff\tr(\hat{B}))$, where  $\hat{C} = \tfrac{\nu}{2} \|R^{-\top} (X - 1_n\hat{\zeta}^\top)\|^2$, with $X \in \R^{n\times d}$ the matrix whose $i$-th row corresponds to the point $x_i$ and $1_n \in \R^n$ the vector where each element is $1$.
The constants $C_0$, defined in \cref{thm:inequality-scattered-data}, and $m, \Cmul, \Cdiff$, defined in \cref{asm:H-rich:a,asm:H-rich:d}, do not depend on $n, \widehat{X}, h_{\widehat{X},\Omega}, \hat{c}, \hat{B}$ or $f$.
\et
\bpr
We divide the proof in two steps

\noindent{\bf Step 1.} First note that 
$$\hat{g}(x) := f(x) - \tfrac{\nu}{2}\|x\|^2 + \nu x^\top \hat z - \hat{c} = f(x) - \tfrac{\nu}{2}\|x - \hat{z} \|^2  - \hat{f}.$$
By applying \cref{thm:inequality-scattered-data} with $g = \hat{g}$ and $B = \hat{B}$ we have
that
for any $x \in \Omega$ 
$f(x) - \tfrac{\nu}{2}\|x -\hat{z}\|^2 - \hat{f} = \hat{g}(x) \geq - \eps - 2\tau$, where $\eps = C_0(|\hat{g}|_{\Omega,m} + \Cmul\Cdiff\tr(\hat{B})) h_{\widehat{X},\Omega}^m$ and $C_0$ is defined in \cref{thm:inequality-scattered-data}. In particular this implies that
$$f(\zeta) - \hat{f} - \tfrac{\nu}{2}\|x -\hat{z}\|^2 \geq - \eps - 2\tau,$$
from which \cref{eq:cert-2:2} is obtained by considering that $f(\hat{z}) \geq f(\zeta)$ since $\zeta$ is a minimizer of $f$. To conclude the proof of \cref{eq:cert-2:2} note that $|\hat{g}|_{\Omega,m} \leq |f|_{\Omega,m} + \nu$ since $m \geq 2$.

\noindent{\bf Step 2.} Now to obtain \cref{eq:cert-2:1} we need to do a slightly different construction. Let $u_j(x) = e_j^\top (x - \hat{z})$ for any $x \in \Omega$. Note that since $u_j$ is the restriction to $\Omega$ of a $C^\infty$ function on $\R^d$, by \cref{asm:H-rich:a}, $u_j \in \hh$. Moreover, note that $\tfrac{\nu}{2}\|x-\hat{z}\|^2 = \tfrac{\nu}{2}\sum_{j=1}^d u_j(x)^2$. Take $\hat{u}_j \in \R^n$ defined as $\hat{u}_j = V^*u_j$ and note that 
$$\Phi_i^\top \hat{u}_j = \big<V\phi(x_i),V^*u_j\big> =  \big<V^*V\phi(x_i),u_j\big> = \big<P\phi(x_i),u_j\big> = u_j(x_i).$$
Then, defining $\hat{G} = \tfrac{\nu}{2}\sum_{i=1}^d \hat{u}_j \hat{u}_j^\top \in \pdm(\R^n)$ we have $$\tfrac{\nu}{2}\|x_i-\hat{z}\|^2 = \Phi_i^\top \hat{G} \Phi_i, \qquad \forall i \in [n].$$
Substituting $-\tfrac{\nu}{2}\|x_i\|^2 + \nu x_i^\top \hat{z}$ with $\tfrac{\nu}{2}\|\hat{z}\|^2 - \Phi_i^\top \hat{G} \Phi_i$ in the inequality in \cref{eq:empirical-inequality-z}, we obtain
$$|f(x_i) - \hat{f} ~-~ \Phi_i^\top (\hat{B} + \hat{G}) \Phi_i^\top| \leq \tau , \quad \forall i \in [n].$$
By applying \cref{thm:inequality-scattered-data} with $g(x) = f(x)  - \hat{f}$ and $B = \hat{B} + \hat{G}$ we have that $f(x) - \hat{f} \geq - \eps -2 \tau$ for all $x \in \Omega$, 
where $\eps = C' h^m_{\hat{X},\Omega}$ with $C' = C_0 (|g|_{\Omega,m}+\Cmul\Cdiff\tr(\hat{B} + \hat{G}))$.
In particular, $f(\zeta) - \hat{f} \geq - \eps - 2\tau$, from which \cref{eq:cert-2:1} is obtained considering that $f(\hat{z}) \geq f_*$ since $\zeta$ is a minimizer of $f$.

Finally, note that $|g|_{\Omega,m} \leq |f|_{\Omega,m}$ since $m \geq 1$.
The proof is concluded by noting that using the definition of $V$ we have $\hat{u}_j = R^{-\top}\hat{v}_j$ with $\hat{v}_j \in \R^n$ corresponding to $\hat{v}_j = (u_j(x_1), \dots, u_j(x_n))$ for $j \in [d]$ and that $\tr(\hat{G}) = \tfrac{\nu}{2}\sum_{j\in [d]} \|\hat{u}_j\|^2$. In particular, some basic linear algebra leads to $\tr(\hat{G}) = \tfrac{\nu}{2}\|R^{-\top} (X - 1_n\hat{z}^\top)\|^2$.
\epr

\end{document}

%% file: macro.tex
\newcommand{\Nystrom}[1]{{Nystr\"om}}

\providecommand{\scal}[2]{\left\langle{#1},{#2}\right\rangle}

\providecommand{\ran}[1]{\operatorname{ran}{#1}}
\providecommand{\tr}{\operatorname{Tr}}

\newcommand{\R}{\mathbb R}

\newcommand{\N}{\mathbb N}
\newcommand{\Np}{\mathbb{N}_+}
\newcommand{\Nz}{\mathbb{N}}

\newcommand{\hh}{\mathcal H}

\newcommand{\la}{\lambda}
\newcommand{\eps}{\varepsilon}
\newcommand{\lspan}[1]{\operatorname{span}\{#1\}}

\newcommand{\C}{{C}}

\newcommand{\eqals}[1]{\begin{align*}#1\end{align*}}
\newcommand{\eqal}[1]{\begin{align}#1\end{align}}
\newcommand{\bpr}{\begin{proof}}
\newcommand{\epr}{\end{proof}}
\newcommand{\be}{\begin{equation}}
\newcommand{\ee}{\end{equation}}

\newtheorem{definition}{Definition}
\newcommand{\bd}{\begin{definition}}
\newcommand{\ed}{\end{definition}}
\newlist{enumdef}{enumerate}{1} 
\setlist[enumdef]{label=\upshape(\alph*),ref=\upshape\thedefinition(\alph*)}
\crefalias{enumdefi}{theorem} 

\newcommand{\bi}{\begin{itemize}}
\newcommand{\ei}{\end{itemize}}

\newtheorem{ass}{Assumption}
\newcommand{\ba}{\begin{ass}}
\newcommand{\ea}{\end{ass}}
\newlist{enumasm}{enumerate}{1} 
\setlist[enumasm]{label=\upshape(\alph*),ref=\upshape\theass(\alph*)}
\crefalias{enumasmi}{ass} 

\newtheorem{remark}{Remark}
\newcommand{\br}{\begin{remark}}
\newcommand{\er}{\end{remark}}
\newlist{enumrem}{enumerate}{1} 
\setlist[enumrem]{label=\upshape(\alph*),ref=\upshape\theremark(\alph*)}
\crefalias{enumremi}{theorem} 

\newtheorem{proposition}{Proposition}
\newcommand{\bp}{\begin{proposition}}
\newcommand{\ep}{\end{proposition}}
\newlist{enumprop}{enumerate}{1} 
\setlist[enumprop]{label=\upshape(\alph*),ref=\upshape\theproposition(\alph*)}
\crefalias{enumpropi}{proposition} 

\newtheorem{lemma}{Lemma}
\newcommand{\blm}{\begin{lemma}}
\newcommand{\elm}{\end{lemma}}
\newlist{enumlm}{enumerate}{1} 
\setlist[enumlm]{label=\upshape(\alph*),ref=\upshape\thelemma(\alph*)}
\crefalias{enumlmi}{lemma} 

\newtheorem{theorem}{Theorem}
\newcommand{\bt}{\begin{theorem}}
\newcommand{\et}{\end{theorem}}
\newlist{enumthm}{enumerate}{1} 
\setlist[enumthm]{label=\upshape(\alph*),ref=\upshape\thetheorem(\alph*)}
\crefalias{enumthmi}{theorem} 

\newtheorem{corollary}{Corollary}
\newcommand{\bcor}{\begin{corollary}}
\newcommand{\ecor}{\end{corollary}}
\newlist{enumcor}{enumerate}{1} 
\setlist[enumcor]{label=\upshape(\alph*),ref=\upshape\thecorollary(\alph*)}
\crefalias{enumcori}{corollary} 

\newtheorem{example}{Example}
\newcommand{\bex}{\begin{example}}
\newcommand{\eex}{\end{example}}
\newlist{enumex}{enumerate}{1} 
\setlist[enumex]{label=\upshape(\alph*),ref=\upshape\theexample(\alph*)}
\crefalias{enumexi}{theorem} 

\crefname{ass}{Assumption}{Assumptions}
\crefname{equation}{Eq.}{Eqs.}
\crefname{figure}{Fig.}{Figs.}
\crefname{table}{Table}{Tables}
\crefname{section}{Sec.}{Secs.}
\crefname{theorem}{Theorem}{Theorems}
\crefname{lemma}{Lemma}{Lemmas}
\crefname{corollary}{Cor.}{Cors.}
\crefname{example}{Example}{Examples}
\crefname{appendix}{Appendix}{Appendixes}
\crefname{remark}{Remark}{Remark}
\crefname{page}{page}{page}

\newcommand{\rank}{\operatorname{rank}}